\documentclass[12pt]{amsart}
\usepackage{comment}
\usepackage[utf8]{inputenc}
\usepackage[T1]{fontenc}
\DeclareMathAlphabet{\mathpzc}{OT1}{pzc}{m}{it}
\usepackage{geometry}

\usepackage[dvipsnames]{xcolor}

\usepackage{enumitem}
\usepackage{marginnote}
\usepackage{amsfonts}
\usepackage{amssymb}
\usepackage{amsthm}
\usepackage{amsmath}
\usepackage{amscd}
\usepackage{mathrsfs}

\usepackage{tikz}
\usetikzlibrary{calc,arrows,decorations.pathreplacing}
\usepackage{nicefrac, xfrac}
\usepackage{mathtools,xparse}
\usepackage{stackengine}
\usepackage{tikz}

\usepackage[pagebackref, colorlinks = true,
linkcolor = red,
urlcolor  = blue,
citecolor = blue,
anchorcolor = blue]{hyperref}
\hypersetup{
  pdftitle={Conditional limit theorems in dynamical systems},
  pdfauthor={Jason Atnip, Gary Froyland, Cecilia González-Tokman, Yushi Nakano, Sandro Vaienti}
}

\theoremstyle{plain}

\def\a{\alpha}

\def\osc{\text{osc}}
\newtheorem{theorem}{Theorem}[section]
\newtheorem{corollary}[theorem]{Corollary}
\newtheorem{lemma}[theorem]{Lemma}

\newtheorem{proposition}[theorem]{Proposition}

\theoremstyle{definition}
\newtheorem{definition}[theorem]{Definition}

\newtheorem{remark}[theorem]{Remark}

\theoremstyle{plain}
\newtheorem{thmx}{Theorem}
\newtheorem*{theorem*}{Theorem}

\newcommand{\vp}{\varphi}
\renewcommand{\epsilon}{\varepsilon}

\newcommand{\sub}{\ensuremath{ \subseteq}} % A shortcut for subset notation
\newcommand{\set}[1]{\ensuremath{ \left\{#1\right\} }} % set notation
\newcommand{\norm}[1]{\ensuremath{ \|#1\| }} % norm

\newcommand{\vertiii}[1]{{\left\vert\kern-0.25ex\left\vert\kern-0.25ex\left\vert #1
		\right\vert\kern-0.25ex\right\vert\kern-0.25ex\right\vert}}

\newcommand{\bs}{\ensuremath{\backslash}}

\newcommand{\mas}{\ensuremath{\text{$m$--a.s.}}}
\newcommand{\mae}{\ensuremath{\text{$m$--a.e.}}}
\newcommand{\lt}{\ensuremath{\left}}
\newcommand{\rt}{\ensuremath{\right}}

\newcommand{\Om}{\ensuremath{\Omega}}

\newcommand{\om}{\ensuremath{\omega}}
\newcommand{\lm}{\ensuremath{\lambda}}
\newcommand{\gm}{\ensuremath{\gamma}}
\newcommand{\al}{\ensuremath{\alpha}}

\newcommand{\sg}{\ensuremath{\sigma}}

\newcommand{\eps}{\ensuremath{\epsilon}}
\newcommand{\dl}{\ensuremath{\delta}}
\newcommand{\kp}{\ensuremath{\kappa}}
\newcommand{\zt}{\ensuremath{\zeta}}

\newcommand{\vta}{\ensuremath{\vartheta}}

\newcommand{\vrho}{\ensuremath{\varrho}}

\DeclareSymbolFont{bbold}{U}{bbold}{m}{n}
\DeclareSymbolFontAlphabet{\mathbbold}{bbold}

\newcommand{\ind}{\ensuremath{\mathbbold{1}}}
\DeclareMathOperator*{\esssup}{ess\sup}
\DeclareMathOperator*{\essinf}{ess\inf}

\newcommand{\cB}{\ensuremath{\mathcal{B}}}

\newcommand{\cL}{\ensuremath{\mathcal{L}}}
\newcommand{\cM}{\ensuremath{\mathcal{M}}}

\newcommand{\cP}{\ensuremath{\mathcal{P}}}

\newcommand{\sB}{\ensuremath{\mathscr{B}}}

\newcommand{\sF}{\ensuremath{\mathscr{F}}}

\newcommand{\NN}{\ensuremath{\mathbb N}}%Natural Numbers

\newcommand{\RR}{\ensuremath{\mathbb R}}% Real Numbers

\def\var{\text{{\rm var}}}

\def\BV{\text{{\rm BV}}}
\def\Leb{\text{{\rm Leb}}}

\def\lt{\left}
\def\rt{\right}

\providecommand{\phantomsection}{}% In case hyperref is not loaded
\AtBeginDocument{\let\textlabel\label}% http://tex.stackexchange.com/q/9939/5764
\makeatletter
\newcommand{\mylabel}[2]{\raisebox{.7\normalbaselineskip}{\phantomsection}(#1)%
	\def\@currentlabel{#1}\textlabel{#2}}
\makeatother

\makeatletter
\newcommand\xlabel[2][]{\phantomsection\def\@currentlabelname{#1}\label{#2}}
\makeatother

\ExplSyntaxOn
\NewDocumentCommand{\mathlist}{ O{,} m m }
 {
  \egreg_mathlist:nnn { #1 } { #2 } { #3 }
 }

\seq_new:N \l__egreg_mathlist_seq
\cs_new_protected:Npn \egreg_mathlist:nnn #1 #2 #3
 {
  \seq_set_split:Nnn \l__egreg_mathlist_seq { #1 } { #3 }
  \seq_use:Nnnn \l__egreg_mathlist_seq { #2 } { #2 } { #2 }
 }
\ExplSyntaxOff

\newcommand{\oio}{\ensuremath{\omega\in\Omega}}

\allowdisplaybreaks

\newcommand{\forevery}{for every }

\let\oldtocsection=\tocsection
\let\oldtocsubsection=\tocsubsection
\renewcommand{\tocsection}[2]{\hspace{0em}\bf\oldtocsection{#1}{#2}}
\renewcommand{\tocsubsection}[2]{\hspace{1.8em}\oldtocsubsection{#1}{#2}}
\let\oldtocsubsubsection=\tocsubsubsection
\renewcommand{\tocsubsubsection}[2]{\hspace{4.2em}\oldtocsubsubsection{#1}{#2}}

\numberwithin{equation}{section}
\title[Conditional limit theorems in dynamical systems]{Conditional limit theorems in dynamical systems}
\date{\today}

\author[J.Atnip]{Jason Atnip}
\address{School of Mathematics and Statistics, University of New South Wales, Sydney, NSW 2052, Australia}
\email{\href{j.atnip@unsw.edu.au}{j.atnip@unsw.edu.au} }

\author[G. Froyland]{Gary Froyland}
\address{School of Mathematics and Statistics, University of New South Wales, Sydney, NSW 2052, Australia}
\email{\href{g.froyland@unsw.edu.au}{g.froyland@unsw.edu.au} }

\author[C. Gonz\'alez-Tokman]{Cecilia Gonz\'alez-Tokman}
\address{School of Mathematics and Physics, The University of Queensland, Brisbane, QLD 4072, Australia}
\email{\href{cecilia.gt@uq.edu.au}{cecilia.gt@uq.edu.au} }

\author[Y. Nakano]{Yushi Nakano}
\address{Faculty of Science, Hokkaido University, Kita 10, Nishi 8,
Kita-Ku, Sapporo, Hokkaido, 060-0810, JAPAN}
\email{\href{yushi.nakano@math.sci.hokudai.ac.jp}{yushi.nakano@math.sci.hokudai.ac.jp} }

\author[S. Vaienti]{Sandro Vaienti}
\address{Universit\'e de Toulon, Aix Marseille Universit\'e, CNRS, CPT, 13009 Marseille, France and Great Bay University, Dongguan, China}
\email{\href{vaienti@cpt.univ-mrs.fr}{vaienti@cpt.univ-mrs.fr} }

\begin{document}

\begin{abstract}
This paper develops a comprehensive theory of conditional limit theorems for open dynamical systems.
We study the statistical properties of the dynamics within transient regimes by establishing a suite of limit theorems describing the behavior of trajectories under the condition that they remain outside a ``hole'', or region of no return.
Under a spectral gap assumption for the associated open transfer operator, we first establish conditioning and Yaglom limits and identify the corresponding conditionally invariant and conditionally ergodic measures.
We then establish conditional exponential decay of correlations and the Strong Law of Large Numbers for both deterministic and quenched random systems.
We also obtain a conditional Central Limit Theorem, Large Deviation Principle, and Berry--Esseen theorem, with additional assumptions in the quenched random case.
We also derive standard (non-conditional) limit theorems on the surviving set of points that never enter the hole.
\end{abstract}

	\maketitle

\tableofcontents

\section{Introduction}

Open dynamical systems, in which trajectories are removed once they enter a designated ``hole'' or region of no return, provide a natural framework for studying transient and metastable behaviour in deterministic and random dynamics.
The time-asymptotic behaviour of such systems has been investigated extensively through escape rates, conditionally invariant measures, quasi-stationary behaviour, and Yaglom limits \cite{Yaglom47,PianigianiYorke,ColletMS,Colletetal00,LMD,DemersYoung06,AFGTV-TFPF}.

In this work we turn to complementary statistical questions that remain much less developed, namely \textit{conditional limit laws} for observables when the underlying law is conditioned on survival up to a large but finite time.
We are partly motivated by the fact that the loss of mass in open dynamical systems means that observed ensembles depend on the observation horizon.
In long transient simulations only those trajectories that have survived up to the time of observation are available for analysis.
This point of view is naturally phrased in terms of \(N\)-conditioned measures, obtained by conditioning an initial distribution on survival up to time \(N\).
These measures describe the statistics seen in long transient experiments, and they are the primary objects of this paper.
The resulting questions include conditional ergodicity, laws of large numbers, central limit theorems, large deviations, and Berry--Esseen estimates.

To the best of our knowledge, rigorous results addressing these questions remain remarkably scarce.
In deterministic open dynamics, an early result is the conditional ergodicity and central limit theorem of Collet--Mart\'inez \cite{ColletMartinez99} for subshifts of finite type with respect to the distinguished conformal measure.
More recent results on conditional limit laws for open dynamical systems and absorbing Markov chains include \cite{LNN-open,BassolsCastro2026,WangSu2026}.
The present paper develops a systematic abstract transfer operator framework for conditioning limits and conditional statistical laws.
Our setting accommodates both deterministic systems and quenched random dynamical systems in the standard cocycle framework \cite{arnold_random_1998,crauel_random_2002}, and allows a broad class of initial distributions that are absolutely continuous with respect to a conformal reference measure.
Our abstract assumptions are verified for several standard classes of dynamics, including random piecewise expanding interval maps with random holes \cite{AFGTV20,AFGTV-IVC,AFGTV-TFPF}, multidimensional expanding maps with random holes \cite{KLG,saussol00}, and subshifts of finite type with holes \cite{pollicott_open_2017}.

Under a spectral gap assumption \eqref{G} for the normalized open transfer-operator cocycle, we first prove convergence to conditioning limits (Theorem~\ref{main thm qstat and qerg measures}), Yaglom limits (Theorem~\ref{main thm yaglom_limit}(1)--(2)), uniqueness of the conditionally invariant measure (Theorem~\ref{main thm yaglom_limit}(3)), and conditional exponential mixing (Theorem~\ref{thm:c}).
Building on these estimates, we establish conditional ergodic behaviour, a conditional weak law of large numbers, and a conditional strong law of large numbers (Theorem~\ref{thm:e1}).
In particular, the conditional strong law in Theorem~\ref{thm:e1}(3) seems to be new already in its formulation.
Under the stronger limit law condition \eqref{LL}, and using a perturbative transfer-operator approach adapted from the literature for random closed systems \cite{DFGTV18b,DFGTV20,CNW,AFGTV-TFPF,DH-Limit20}, we further prove a conditional central limit theorem (Theorem~\ref{clt}), a conditional large deviation principle (Theorem~\ref{ldp}), and a
conditional Berry--Esseen theorem (Theorem~\ref{be}).
The corresponding dynamical variance is identified via a Green--Kubo type formula (Proposition~\ref{prop:0625a}).
%\marginpar{\tiny small change for this sentence?}
Under suitable continuity and boundary assumptions, these estimates also yield standard dynamical properties and limit laws on the surviving set (Proposition~\ref{cor:0201a} and Corollary~\ref{cor:CLTLDPXinfty}).
Furthermore, the method is not tied to these particular statements: with the appropriate perturbative or rare event estimates, the same conditional transfer operator framework should be useful for conditional versions of other statistical laws, such as functional central limit theorems, laws of the iterated logarithm, extreme value laws, Erd\H{o}s--R\'enyi laws, and almost sure invariance principles.

It is important to distinguish these conditional limit laws from the ordinary limit laws for the restricted dynamics on $X_\infty$, consisting of points whose forward orbits never enter the hole.
Invariant measures on $X_\infty$ describe distributions of points that survive forever, whereas the conditioned laws we introduce describe the statistics of trajectories selected by survival up to a large finite time.
This distinction is important in applications: the basin of the invariant measure on $X_\infty$ is frequently negligible with respect to natural initial distributions, such as Lebesgue measure, and therefore is only informative for a negligible set of initial conditions.
In contrast, our finite-time conditioned measures and our conditional laws capture the statistics of non-negligible sets of initial conditions: those that are  visible in long transient observations.

The paper is organized as follows. Section~\ref{sec: rds} introduces (closed) random dynamical systems, and Section \ref{sec: ords} defines open random dynamical systems and the associated open random transfer operator cocycle. 
In Section \ref{s:main_defs} we define the conditioning limit, Yaglom limit, and conditionally ergodic limit for random (and deterministic) probability measures.
Section~\ref{s:results} contains statements of all of our main results. 
Section~\ref{s:example} verifies our assumptions for several example classes of deterministic and random dynamics.
Section~\ref{sec: proof of A,B,C} proves Theorems \ref{main thm qstat and qerg measures}--\ref{thm:c}, Section~\ref{ss:variance} justifies the existence of a limiting conditional dynamical variance, and Section~\ref{s:thmD} proves Theorem \ref{thm:e1}.
Section~\ref{s:taylor} contains crucial derivative information on the leading Lyapunov exponent of the twisted transfer operator cocycle, and Section~\ref{sec: proof of E, F, G} contains the proofs of Theorems \ref{clt}--\ref{be}.

%\bigskip

\section{Random dynamical systems}
\label{sec: rds}
Let $(\Om,\sF,m)$ be a probability space and $\sg:\Om\to\Om$ an invertible, ergodic map which preserves the measure $m$. 
Let $(\Xi ,d)$ be a complete, 
%locally 
compact metric space equipped with the Borel $\sg$-algebra $\sB$ and let $\{ X_\om \}_{\omega\in\Omega}$ be a closed random subset of $\Xi$ 
(i.e.~$X_\om$ is a closed subset of $\Xi$ for each $\oio$ and $\om\mapsto d(x,X_\om)$ is measurable \forevery $x\in \Xi$),
in particular, 
\[
X_\om \subseteq \Xi \quad \text{ for each $\oio$.}
\]
We let 
\[
X:=\bigcup_{\oio}\{\om\}\times X_\om \subset \Om \times \Xi .
\]
Further, we consider a family of surjective maps $\{T_\om: X_\om \to X_{\sg\om}\}_{\omega\in\Omega}$ and define the induced skew-product map $T: X \to X $ given by
\[
T(\omega,x)= (\sg\om,T_\om(x)).
\]
We make the following assumption:
\begin{enumerate}[align=left,leftmargin=*,labelsep=\parindent]
\item[\mylabel{M1}{M1}] The map $T:X\to X$ is measurable with respect to product $\sg$-algebra $\sF\otimes\sB$.
\end{enumerate} 
Random iteration of our maps $T_\om$ is denoted by $T_\om^n: X_\om \to X_{\sg^n\om} $ and is given by 
\begin{align}\label{eq:20240613c}
    T_\om^n:=T_{\sg^{n-1}\om}\circ T_{\sg^{n-2}\om}\circ\cdots\circ T_\om
\end{align}
for each $n\geq 1$. Notice that $T^n(\om,x)=(\sg^n\om,T_\om^n(x))$.

A probability measure $\zeta$ on $(\Om\times \Xi ,\sF\otimes\sB)$ is said to be a \textit{random probability measure} relative to $m$ if it has marginal $m$ (i.e.~$\zeta\circ\pi^{-1}_1=m$, where $\pi_1(\omega,x)=\omega$).
This property holds if and only if the disintegrations $\set{\zeta_\om}_{\om\in\Om}$ of $\zeta$ with respect to the partition $\{\{\om\}\times \Xi \}_{\om\in\Om}$ satisfy the following properties:
		\begin{enumerate}
			\item For every $B\in\sB$, the map $\Om\ni\om\longmapsto\zeta_\om(B)\in [0,1]$ is measurable, 
			\item For $m$-a.e.~$\om\in\Om$, the map $\sB\ni B\longmapsto\zeta_\om(B)\in [0,1]$ is a Borel probability measure.
		\end{enumerate}
         Refer to 
         \cite{arnold_random_1998,crauel_random_2002}.
         Thus, we often identify a random probability measure $\zeta$ with $\set{\zeta_\om}_{\omega\in\Om}$ by writing $\zeta=\set{\zeta_\om}_{\omega\in\Om}$.
Given a measurable $A\sub X$, we let $\cM_1(A)$ denote the set of all random probability measures supported in $A$. If $A=X$, we simply write $\cM_1$ for $\cM_1(X)$.
  We also define a random measure 
  $\zeta =\set{\zeta_\om}_{\omega\in\Om}$
   by (1) and (2) with ``measure'' instead of ``probability measure''.

\subsection{Closed random transfer operator}
We consider a Banach space $(\cB,\Vert \cdot \Vert )$ 
of complex-valued Borel measurable functions on $\Xi$ and
a random weight function $g: X \to (0,\infty)$ satisfying the following properties.
We define the (closed) random transfer operator $\cL_{\om,c}:\cB\to\cB$ associated with the weight $g_\om :=g(\omega,\cdot )$  as follows:
for each $\oio$, $h\in\cB$, and $x \in X_{\sg\om}$, 
\begin{align}\label{eq: L_c def}
    \cL_{\om,c}(h)(x) := \sum_{y \in T_\om^{-1}(x)} g_\om(y)h(y).
\end{align}
The same formula defines $\cL_{\om,c}(h)$ for every bounded Borel measurable function $h:\Xi\to\mathbb C$.
We assume that the weight $g$ is chosen so that $\cL_c$ is a well-defined, bounded linear operator. In fact we will assume that there is $C_1>0$ such that $\cL_{\om,c}h\in\cB$ and 
\begin{align}\label{C0}\tag{C0}
    \|\cL_{\om,c}h\|\le C_1\|h\|
\end{align}
for almost every $\oio$ and $h\in\cB$,
and that 
\begin{enumerate}[align=left,leftmargin=*,labelsep=\parindent]
%\label{eq:measurabilityofLc}
\item[\mylabel{M2}{M2}] the map $(\om, x) \mapsto (\cL_c h)_\om(x)$
 is measurable \forevery $\mathcal B$-valued random variable $h$, where
 \[
 (\cL_c h)_\om:=\cL_{\sg^{-1}\om,c}h_{\sg^{-1}\om} \qquad \text{for $\oio$}.
 \]
\end{enumerate}
\begin{remark}
    Note that if the random maps $T_\om$ are finite-to-one and the weight functions $g_\om$ are bounded, then the transfer operator $\cL_{\om,c}$ is always well-defined. On the other hand if the maps $T_\om$ are countable-to-one, then we require that the weight functions are summable. See \cite{AFGTV20} for details. 
\end{remark}

We make the following assumption on the closed system. Moreover, for a measurable set $A\subset \Xi$, let $\hat A$ denote its indicator function. For a probability measure $\mu$ on $\Xi$ and a $\mu$-integrable measurable function $h:\Xi\to\mathbb C$, we use the shorthand $\mu(h):=\int h \,d\mu$.
\begin{enumerate}[align=left,leftmargin=*,labelsep=\parindent]
\item[\mylabel{C1}{C1}] 
There exists a random non-atomic Borel probability $\nu_c$ such that
\begin{align}\label{eq:20240618}
\nu_{\sg\om,c}(\tilde\cL_{\om,c} h) = \nu_{\om,c} (h)
\end{align}
for $m$-a.e.~$\oio$ and every bounded Borel function $h:\Xi\to\mathbb C$, where
\[
\tilde{\mathcal L}_{\om,c}h:= \lambda_{\om,c}^{-1} \mathcal L_{\om,c}h, \quad\text{ and }\quad \lambda_{\om,c}:=\nu_{\sg\om,c}(\mathcal L_{\om,c}\hat X_\om).
\]
\end{enumerate}

Notice that 
\[
\mathcal L_{\om,c}(f\circ T_\om \cdot h) =f\cdot \mathcal L_{\om,c}(h)
\]
\forevery  pair of functions $f,h\in\cB$ by the definition of $\mathcal L_c$ in \eqref{eq: L_c def}. 
Therefore, we get the following important consequence of  \eqref{eq:20240618} regarding duality:
\begin{equation}\label{eq:duality}
\nu_{\om,c}(f\circ T_\om \cdot h)= \nu_{\sg\om,c}(f\cdot \tilde{\mathcal L}_{\om,c} h).
\end{equation}	
In fact, below we will only use \eqref{eq:duality} instead of \eqref{eq:20240618}.

\begin{remark}
    The random probability measure $\nu_c$ is known as the random (closed) conformal measure. Such measures have been shown to exist in many random settings, e.g.  \cite{Bogenschutz_RuelleTransferOperator_1995a,Bruck_Generalizediteration_2003,mayer_distance_2011,mayer_countable_2015,mayer_random_2018,AFGTV20,AFGTV-IVC}.
\end{remark}

\section{Open random dynamical systems}
\label{sec: ords}
Let $H$ be a measurable set in $ X $.
The fibers $H_\om$ are determined by the condition that $\{\om\}\times H_\om = H \cap (\{ \om \} \times X_\om)$.
For each $\oio$ and $n\geq 1$ we define
	\begin{align}\label{def: X_n surv}
		X_{\om,n}:&
		=\bigcap_{j=0}^n T_\om^{-j}\left(X_{\sg^j\om,0}\right)
	\end{align}
to be the set of points in $X_{\om,0}:=X_\om \bs H_\om$ whose trajectories avoid the hole $H$ for $n$ iterates, where $T^{-j}_\om :=(T_\om^j)^{-1}$.
Note that $X_{\om,n}$ can be expressed as $X_{\om,n}=\{\tau_\om>n\}$, where $\tau_\om$ is the first hitting time to $H_\om$ (i.e.~$\tau_\omega(x)=\inf\{j\ge0:T_\omega^j(x)\in H_{\sigma^j\om}\}$).
 We call $X_{\om,n}$ the \textit{$n$th level random surviving set}. 
We naturally define
	\begin{align}\label{def: X_infty surv}
		X_{\om,\infty}:=\bigcap_{n=0}^\infty X_{\om,n}=\bigcap_{n=0}^\infty T_\om^{-n}(X_{\sg^n\om,0})
	\end{align}
to be the set of points whose trajectories will never land in the hole under iteration of the maps $T^n_\om$ \forevery $n\geq 0$.
We call $X_{\om,\infty}$ the \textit{random surviving set}. Note that the random surviving set 
satisfies the following fiberwise invariance properties: 
\begin{align}
    	    T_\om(X_{\om,n})\sub X_{\sg\om,n-1}
	    \qquad\text{ and }\qquad
	    T_\om(X_{\om,\infty})\sub X_{\sg\om,\infty}  
	\label{surv set forw inv}
	\end{align}
for every $n\ge 1$ and a.e. $\om$.

\subsection{Random open transfer operator}
Now we define the random open (perturbed) transfer operators $\cL_\om:\cB\to\cB$ by
\begin{align*}
\cL_\om h:=\cL_{\om,c}\left(h\cdot\hat X_{\om,0}\right),
\qquad h\in\cB, \quad \oio
\end{align*}
(recall that $\hat A$ is the indicator function of $A$).
We let
\[
(\mathcal Lh)_\om := \mathcal L_{\sg^{-1}\om}h_{\sg^{-1}\om} 
\]
for each $\mathcal B$-valued random variable $h$ and $\oio$.
Note that the measurability of the hole $H$ and \eqref{M2} imply that for every $\cB$-valued random variable $h$ the map $(\om,x)\mapsto (\cL h)_\om (x)$ is also measurable.
Iterates of the open operator $\cL_\om^n:\cB \to\cB $ are given by
\begin{align*}
\cL_\om^n:=\cL_{\sg^{n-1}\om}\circ\cL_{\sg^{n-2}\om}\circ\dots\circ\cL_\om,
\end{align*}
which, using induction, we may write as
\begin{align}\label{eq: iterate pert tr op}
\cL_\om^n h=\cL_{\om,c}^n\left(h\cdot\hat{X}_{\om,n-1}\right), \qquad h\in\cB.
\end{align}
In the global notation, we have
\[
(\mathcal L^nh)_\om = \mathcal L^n_{\sg^{-n}\om}h_{\sg^{-n}\om}
\]
for each $\mathcal B$-valued random variable $h$ and $\oio$.
We make the following assumptions on the Banach space $\cB$:
\begin{enumerate}
    \item[\mylabel{B1}{B1}] $\cB$ is dense in $L^1(\nu_{\om,c})$ for every $\omega\in\Omega$ and the set 
    $\cB\cap C(\Xi)$ 
    is dense in $C(\Xi)$ the space of complex-valued continuous functions on $\Xi$.
    \item[\mylabel{B2}{B2}] $\widehat\Xi\in\cB$, and there is $C_0>0$ such that
        $$
        \|h_1h_2\|\le C_0\|h_1\|\cdot\|h_2\|\qquad\mbox{and}\qquad
        \|h\|_\infty:=\sup_{x\in \Xi}|h(x)| \le C_0\Vert h\Vert
        $$
    \forevery $h_1,h_2,h\in\mathcal B$.
     \item[\mylabel{B3}{B3}] $\hat X_{\om}, \hat X_{\om,n} \in\cB
    $ for $\mae$ $\oio$ and every $n\in \NN$.
    Moreover, the random variable
$
\omega\mapsto\|\hat X_{\omega,0}\|
$
is tempered with respect to $\sigma$.
\end{enumerate}

We are now ready to state our main technical assumptions. We assume the following:
\begin{enumerate}[align=left,leftmargin=*,labelsep=\parindent]
\item[\mylabel{G}{G}] \label{hyp:G}
\textit{Spectral gap condition:}
\begin{enumerate}
\item[\mylabel{G0}{G0}]
There exists a measurable, $m$-a.s. finite map $\om\mapsto K_\om$ such that for all $f\in\cB$
\begin{align*}
    \|\cL_\om f\|\leq K_\om\|f\|.
\end{align*}
\item[\mylabel{G1}{G1}]
There exists a $\cB^*$-valued random variable $\nu_\infty$ such that 
 \[
 \nu_{\sg\om,\infty}( \tilde\cL_\om h)=\nu_{\om,\infty}(h),
 \]
 \forevery $h\in\cB$ and $\mae~\oio$, where 
\begin{align*}
\tilde\cL_\om:=\lm_\om^{-1}\cL_\om \quad \text{and} \quad
\lm_\om:=\nu_{\sg\om,\infty}(\mathcal L_\om\hat X)
\end{align*} 	
for $\log\lm\in L^1(m)$.
\item[\mylabel{G2}{G2}] There exists a $\mathcal B$-valued random variable $\phi$ such that $\nu_c(\phi)$
is tempered\footnote{Recall that a random variable $C:\Om\to\RR$ is said to be tempered if $\lim_{n\to\infty}\frac{1}{n}\log |C\circ\sg^n| =0$ $\mas$, and that a convenient sufficient condition for $C$ to be tempered is $\log C\in L^1(m)$.} and
\begin{align*}
\tilde\cL_\om\phi_\om =\phi_{\sg\om},  \quad \text{and} \quad \nu_{\om,\infty}(\phi_\om)=1
\end{align*}
for $m$-a.e. $\oio$. 
\item[\mylabel{G3}{G3}]
There exists a tempered measurable function $D:\Omega\to(0,\infty)$ and $\kp\in(0,1)$ such that for each $n\in\NN$, each $h\in \mathcal B$, and for $m$-a.e. $\oio$ we have 
\begin{align}\label{eq:stronggap}
   \norm{\tilde\cL_\om^n h - \nu_{\om,\infty}(h)\cdot\phi_{\sg^n\om}}\leq D_{\sg^n\om}\norm{h}\kp^n .
\end{align}
\end{enumerate}
\end{enumerate}

\begin{remark}\label{rmk:0411}
Assume in addition that $\log \lambda_c\in L^1(m)$.
    Assuming the limit exists, one can define the fiberwise \textit{escape rate} of a random measure $\zt$ by
    \begin{align*}
        R(\zt_\om):=-\lim_{n\to\infty}\frac{1}{n}\log\zt_\om(X_{\om,n}).
    \end{align*}
  Under \eqref{C1} and \eqref{G}, the limits exist for $\nu_c$ and for
$\eta_0$ defined in \eqref{d;etai}, and
\[
 R(\nu_{\omega,c})=R(\eta_{\omega,0})
 =\int_\Omega\log\frac{\lambda_{\omega,c}}{\lambda_\omega}\,dm
 \qquad\mas
\]
Indeed, \eqref{C1}, \eqref{eq: iterate pert tr op}, and \eqref{G2} give
\[
 \eta_{\omega,0}(X_{\omega,n})
 =\frac{\nu_{\sigma^{n+1}\omega,c}(\phi_{\sigma^{n+1}\omega})}
        {\nu_{\omega,c}(\phi_\omega\hat X_{\omega,0})}
   \prod_{j=0}^{n}\frac{\lambda_{\sigma^j\omega}}
                          {\lambda_{\sigma^j\omega,c}}.
\]
The terminal scalar is tempered by \eqref{G2}, while the denominator
is fixed. The analogous asymptotic formula for $\nu_c$ follows from
\eqref{G3}, and Birkhoff's theorem gives the displayed rates.
Compare Propositions~1.12.4 and~2.3.4 in \cite{AFGTV-TFPF}
\end{remark}

\begin{remark}
The open conformal measure $\nu_\infty$ also appears naturally in spherical averages for iterated function systems generated by contractions. Consider the deterministic map $T(x)=3x \bmod 1$ with hole $H=[\frac13,\frac23]$, and let
$S_0(x)=\frac{x}{3}$, $S_1(x)=\frac{x+2}{3}$.
For $a=(a_1,\ldots,a_n)\in\{0,1\}^n$, write $S_a:=S_{a_n}\circ\cdots\circ S_{a_1}$ and define the so-called spherical average
\[
A_n f(x):=\frac{1}{2^n}\sum_{a\in\{0,1\}^n} f(S_a x).
\]
The associated open transfer operator is
$
\mathcal Lf(x)=\frac13\bigl(f(S_0x)+f(S_1x)\bigr)
$
modulo Lebesgue-null sets,
so $\lambda=\frac{2}{3}$ and therefore
$\widetilde{\mathcal L} f(x)=\lambda^{-1}\mathcal Lf(x)=\frac12\bigl(f(S_0x)+f(S_1x)\bigr).$
Hence 
\[
A_n f=\widetilde{\mathcal L}^n f
\]
for every $n\geq 1$. In this example $\phi\equiv 1$, so \eqref{eq:stronggap} gives
\[
A_n f\to \nu_\infty(f)
\]
uniformly for $f\in\mathcal B$. 
More generally, whenever the normalised open transfer operator can be written as the branch-average operator of a finite family of contractions, $\nu_\infty$ is the limiting functional and the limiting profile is $\phi$.
\end{remark}

\section{Random conditional measures}\label{s:main_defs}
In this section, we introduce the definitions that will be used in the main theorems below. We stress that these definitions are independent of \eqref{G}.
Let $\mathcal M_1^H$ 
be the set of random measures $\zeta \in \mathcal M_1 $ such that 
\begin{equation}\label{eq:20240612}
\zt_\om(X_{\om,n}) >0 \quad \mas \qquad \mbox{for every $n\ge 0$.}
\end{equation}
In all examples in Section \ref{s:example}, the reference measure $\nu_c$ satisfies \eqref{eq:20240612}; in particular, $\mathcal M_1^H \neq \emptyset$.
For each random (not necessarily probability) measure $\zt$ satisfying the condition \eqref{eq:20240612} and $n\geq 0$, we define the \emph{conditioned measure} 
$\zt_n$ of $\zt$, conditioned on the event $X_n$, by 
\begin{align}\label{eq:con}
     \zt_{\om,n}(A):=\frac{\zt_\om(A\cap X_{\om,n})}{\zt_\om(X_{\om,n})}
\end{align}
for each $\oio$ and measurable set $A\subset X_\om$.
Furthermore, we denote by $T_*(\zt)=\{T_*(\zt)_\om\}_{\oio}$ the pushforward of $\zt$ by $T:X\to X$, so that
\[
T_*(\zeta)_\om(A)=\zt_{\sg^{-1}\om}(T_{\sg^{-1}\om}^{-1}(A))
\]
for each $\oio$ and measurable set $A\subset X_\om$.
In the following definition, convergence of probability measures on $\Xi$ is understood in the weak* topology.
For a sequence of random probability measures $(\mu^{(n)})_{n\ge 1}$ and a random probability measure $\mu$, we write $\lim_{n\to\infty}\mu^{(n)}=\mu$ to mean that $\mu_{\omega}^{(n)}$ converges weakly
to $\mu_\omega$ for $m$-a.e.~$\oio$;  that is $\lim_{n\to\infty} \mu_\omega^{(n)}(f)=\mu_\omega(f)$ for all $f\in C(\Xi)$ and a.e.\ $\omega\in\Omega$.

\begin{definition}\label{def qstat}
  We say that a random probability measure $\vrho \in \mathcal M_1$ is 
    \begin{itemize}
        \item the \textit{conditioning limit} 
        of a random probability measure $\zt\in \mathcal M_1^H$ 
    if 
    \begin{align*}
   \lim_{n\to\infty} \zeta_{n} =\vrho.
    \end{align*}
    (namely,   $\lim_{n\to\infty} \zeta _{\om,n} =\vrho_\om$ for $m$-a.e.~$\oio$).
      \item the \textit{Yaglom limit} 
      of a random probability measure $\zt\in \mathcal M_1^H$ 
    if 
        \begin{align*}
   \lim_{n\to\infty} T_*^n (\zeta _n) =\vrho.
    \end{align*}
    \item \textit{conditionally invariant} 
    if $\vrho \in \mathcal M_1^H$ and\footnote{This is equivalent to the condition $T_*(\varrho_1)=\varrho$: note that $T_*^n (\varrho _n)=\frac{T_* (\varrho _1)(T^{-n+1}(\cdot )\cap X_{n-1})}{T_*(\varrho_1)(X_{n-1})}$.}
    \begin{align*}
T_*^n (\varrho _n) = \varrho \qquad \mbox{for every $n\ge 0$.}
    \end{align*}
    \item the \textit{conditionally ergodic limit} 
    of a random probability measure $\zt\in \mathcal M_1^H$ 
    if
    \begin{align*}
   \lim_{n\to\infty} \frac{1}{n}\sum_{j=0}^{n-1} T_*^j (\zeta_{n-1})= 
   \varrho .
    \end{align*}
  \item  \textit{conditionally ergodic} 
 over a 
 subset $\mathcal D$ of $\mathcal M_1^H$  if $\vrho$ is the conditionally ergodic limit  of $\zt$ \forevery $\zt\in\mathcal D$.
    \end{itemize}
\end{definition}

Analogous terms to the definitions above, with the exception of the conditioning limit, exist in the context of Markov chains, namely quasi-limiting distribution (or Yaglom limit), quasi-stationary distribution, and quasi-ergodic limit, which correspond to Yaglom limit, conditionally invariant measure, and conditionally ergodic limit, respectively (see, e.g., \cite{ColletMS,vanDoornPollett-review,MeleardVillemonais}).
However, as will be discussed in Remark \ref{rmk:1229}, the conditionally ergodic limit and the quasi-ergodic limit are not identical concepts, primarily because the target class of measures, $\mathcal D$, is different (typically restricted to Dirac initial measures, under which the scope of the conditionally ergodic limit is substantially narrowed).
Indeed, the definition of the Yaglom limit in a large part of the Markov chain literature differs for the same reason (the target class $\mathcal D$ for Markov chains is basically the set of Dirac initial measures), but our definition aligns with that found in the monograph by Collet, Mart\' inez, San Mart\' in \cite{ColletMS}.
Furthermore, we adopt the term ``conditionally invariant measure'' following the convention established by Pianigiani and Yorke \cite{PianigianiYorke}.

\begin{remark}\label{rmk:1229}
We will consider (dense subsets of)
\[
 \mathcal P(\vrho):=\{\zt \in \mathcal M_1^H : \text{$\zeta _\om$ is absolutely continuous with respect to $\varrho_\om$ for $m$-a.e.~$\oio$}\}
 \]
 as the target class $\mathcal D$ for Yaglom and conditionally ergodic limits due to the following reason.
 
 When determining absolute continuity of $\zeta$ with respect to $\rho$, we compute Radon-Nikodym derivatives $\omega$-wise:  $d\zeta  / d\vrho:=\{d\zeta _\omega / d\vrho_\om\}_{\oio}$. 
Recall that
 a random probability measure $\vrho$ on $ X $ is called invariant (or equivariant) if $T_*(\varrho )=\varrho$ (so that $T_*^n(\varrho )=\varrho$ for every $n\ge 0$);  note the difference with conditional invariance.
Furthermore, it is not difficult to see that a $T$-invariant random probability measure $\varrho$ is ergodic (i.e.~$\varrho (A)=0$ or $1$ for any $T$-invariant $A\in\mathscr F\times \mathscr B$),
or equivalently by Birkhoff ergodic theorem it holds that
\begin{equation}\label{eq:20240614}
\lim_{n\to\infty} \frac{1}{n}\sum_{j=0}^{n-1}  f  \circ T^j(\om,x) =\lim_{n\to\infty} \frac{1}{n}\sum_{j=0}^{n-1}  f_{\sigma^{j}\omega} \circ T_{\omega}^j(x) 
=
\int  \vrho_\om (f_\omega)\, dm(\omega)
\end{equation}
\forevery $f\in L^1(\vrho)$ and $\vrho$-a.e.~$(\omega,x)$,
if and only if
\begin{equation}\label{eq:20240626}
\lim_{n\to\infty}\frac1n\sum_{j=0}^{n-1}
(T_*^j\zeta)_{\sigma^j\omega}(f_{\sigma^j\omega})
=\int_\Omega\varrho_\omega(f_\omega)\,dm(\omega)
\qquad\text{for $m$-a.e.~$\omega$},
\end{equation}
for every bounded measurable $f:X\to\mathbb R$ and every random probability measure $\zeta$ such that $\zeta_\omega$ is absolutely continuous with respect to $\varrho_\omega$ for $m$-a.e.~$\omega$. See Appendix~\ref{app:B} for the proof.
Compare \eqref{eq:20240626} with the definition of conditional ergodicity.

A related property is
\begin{align}\label{eq:equivalence_mixing}
\lim_{n\to\infty}T_*^n(\zeta)=\varrho
\end{align}
in the weak topology for every random probability measure
$\zeta$ such that $\zeta_\omega$ is absolutely continuous with respect to $\varrho_\omega$ for $m$-a.e.~$\omega$. 
For comparison, the forward fibrewise mixing property is
\begin{equation}\label{eq:equivalence_mixing_2}
\int (f\circ T_\omega^n)u\,d\varrho_\omega
-\int f\,d\varrho_{\sigma^n\omega}\int u\,d\varrho_\omega
\longrightarrow0\qquad(n\to\infty)
\end{equation}
for every $f\in C(\Xi)$, $u\in L^1(\varrho_\omega)$ and
$m$-a.e.~$\omega$. The former is a pullback statement,
whereas the latter fixes the initial fibre; they coincide in
the deterministic case. See Appendix~\ref{app:B}.
Compare \eqref{eq:equivalence_mixing} with the definition of Yaglom limit.

We remark that, in the definition of conditionally ergodic or Yaglom limit for Markov processes, the convergence is usually with respect to the total variation, but according to the above observation, the weak convergence (that we employed here) might be more natural in the context of dynamical systems.

We also remark that there are definitions of conditionally invariant measures and Yaglom limits for deterministic dynamical systems (cf.~\cite{ColletMS}), and that of conditionally invariant measures for quenched random dynamical systems (cf.~\cite{AFGTV-TFPF}).
\end{remark}

\section{Main results}\label{s:results}
In this section we introduce our main results. 
\subsection{Convergence to conditional measures, conditional exponential mixing, and conditional weak law of large numbers}

We introduce 
 \[
\mathcal P_{\mathcal B,+}(\nu_c):=\left\{\zt \in\mathcal P(\nu_c) : \text{$\frac{d\zeta_\om   }{ d\nu_{\om,c}} \in \mathcal B$ and $ \nu_{\om,\infty}\left(\frac{d\zeta  _\om }{ d\nu_{\om,c}} \right) >0$ for $m$-a.e.~$\oio$}\right\}\subset\mathcal{M}_{1}^H.
\]
Furthermore, for simple notation we use
 \[
 (\psi \cdot \nu_c )_\om (B):=\nu _{\om,c} (\hat B \cdot \psi _\om) 
 \]
 for each $\psi \cdot \nu_c$, measurable set $B$ and $\oio$.
 Given a random $\mas$ positive measure $\zeta$ (i.e.~$\zeta (X)>0$ $\mas$),~we denote by $\widetilde\zeta $ its normalization, that is, 
 \[
 \widetilde \zeta_\omega (B):=\frac{\zeta_\omega (B)}{\zeta_\omega (X_\omega)} 
 \]
  for each measurable set $B$ and $\oio$.

Referring to the objects $\nu_\infty$ and $\phi$ in assumption \eqref{G} and the conditioning \eqref{eq:con}, we define the random probability measure $\eta_{0}$ and $\cB^*$-valued random variable $\eta_\infty$ by
 \begin{equation}\label{d;etai}
 \eta_{\om,0}:= (\phi \cdot \nu_c)_{\om,0}
 \qquad\mbox{and}\qquad
\eta_{\om,\infty}:=(\phi \cdot \nu_\infty)_\om\qquad \text{for $m$-a.e.~$\oio$.}
\end{equation}
Finally, for each $\zt \in\mathcal P_{\mathcal B,+}(\nu_c)$, we define
$\cB^*$-valued random variable $\zeta_\infty$ by
    \[
    \zeta_{\om,\infty}(f):=(\widetilde{\psi \cdot \nu_{\infty} })_\om (f):=\frac{\nu _{\om,\infty} (f \cdot \psi_\om)}{\nu _{\om,\infty} (\psi_\om)},\qquad f\in\cB,
    \]
    where $\psi$ is the random density of $\zeta$.

We introduce  a random variable with values in $(0,\infty)$
\begin{equation}\label{eq:0820}
C_{\om,\phi}:= \max\lt\{1+ \frac{C_0\Vert \hat X_{\om,0}  \Vert D_{\sg^k\om}}{\nu_{\sg^k\om,c}(\phi_{\sg^k\om})} \kappa ^k:k\ge 0\rt\},
\end{equation}
which is tempered by \eqref{B3} and \eqref{G}, as shown in Lemma~\ref{eq:20240618a}.
Consequently, 
for any $\zt\in \mathcal P_{\mathcal B,+}(\nu_c)$ with random density $\psi$, 
\begin{equation}\label{eq:nomzt}
n_{\om,\zt}:=\inf\left\{n_0 \in \mathbb N : C_0D_{\sg^n\omega} ||\psi_\om||\kappa ^n <\frac{\nu_{\om,\infty}(\psi_\om)}{2C_{\sg^n\om,\phi}} \quad \text{for all $n\ge n_0$}\right\}
\end{equation}
is finite $\mas$ since $DC_\phi$ is tempered.

\begin{thmx}[Conditioning Limit]\label{main thm qstat and qerg measures}
Assume that the spectral gap condition \eqref{G} holds.  Then, 
the following holds for every $\zt \in\mathcal P_{\mathcal B,+}(\nu_c)$:
\begin{enumerate}
\item     
For every $f\in\mathcal B$, $n\ge n_{\om,\zt}$ and $\mas~\oio$,
\begin{equation}
    \label{eq:zeta conv}
\left\vert \zt _{\om,n}(f) -  \zeta_{\om,\infty} (f) \right\vert  \le \frac{2C_0  (1+C_0\Vert f\Vert  )^2
\Vert \psi_\om\Vert }{ \nu_{\om,\infty}( \psi _\om)  }(D_{\sg^n\om}C_{\sg^n\om,\phi})
\kappa ^n;
\end{equation}
    \item For $m$-a.e.~$\oio$, 
  $\zeta_{\om,\infty}$ is positive, $\zeta_{\om,\infty}(\hat \Xi)=1$ and  $ |\zeta_{\omega,\infty}(f)|\le \|f\|_\infty$
    for every $f\in \mathcal B\cap C(\Xi)$, so that there exists a unique Borel probability measure $\widehat \zeta_{\om,\infty}$ satisfying
    \[
    \zeta_{\om,\infty}(f)=\widehat\zeta_{\om,\infty}(f),
    \qquad f\in \cB\cap C(\Xi).
    \]
    \item   If $\zeta=\nu_c$, then $\zeta_\infty=\nu_\infty$;
if $\zeta=\eta_0$, then $\zeta_\infty=\eta_\infty$.
We denote the corresponding representing measures given in (2)
by $\widehat\nu_\infty$ and $\widehat\eta_\infty$, respectively;
\item 
$\widehat\zeta_\infty$ is the conditioning limit of $\zeta$.
Furthermore, for $m$-a.e.~$\oio$, $\widehat\zeta_{\om,\infty}$ is supported on $\bigcap_{n\ge0}\overline{X_{\omega,n}}$ and $\widehat\nu_{\om,\infty}$-absolutely continuous.
\end{enumerate}
\end{thmx}
Throughout, we use a widehat for a probability measure representing a functional initially defined on $\mathcal B$ (e.g.~$\widehat\zeta_\infty$), but not for probability measures defined directly as such (e.g.~$\zeta_0$ and $\zeta_n$).

Under the continuity hypothesis of Proposition~\ref{cor:0201a}, the representation in Theorem~\ref{main thm qstat and qerg measures}(2) holds on all of $\mathcal B$.
If, in addition,
\begin{equation}\label{eq:survivor-boundary-null}
 \widehat\nu_{\omega,\infty}
 \bigl(\overline{X_{\omega,n}}\setminus X_{\omega,n}\bigr)=0
 \qquad\text{for every $n\ge0$ and $m$-a.e.~$\omega$},
\end{equation}
then $\widehat\nu_{\omega,\infty}(X_{\omega,\infty})=1$ by continuity
from above.
The same is true of every $\widehat\zeta_{\omega,\infty}$
in Theorem~\ref{main thm qstat and qerg measures}, by the absolute continuity of Theorem~\ref{main thm qstat and qerg measures}(4).
Condition~\eqref{eq:survivor-boundary-null} holds if each $X_{\omega,n}$
is closed, or if $\widehat\nu_{\omega,\infty}(\partial X_{\omega,n})=0$
for every $n$; in particular, it holds under the continuity hypothesis
of Proposition~\ref{cor:0201a}.
Furthermore, under this continuity hypothesis, $\widehat\eta_\infty=\phi\widehat\nu_\infty$ is the unique ergodic $T$-invariant probability measure on $X_\infty$ that is absolutely continuous with respect to $\widehat\nu_\infty$; see Proposition~\ref{cor:0201a}(2).

For random expanding interval maps with holes, a convergence estimate corresponding to Theorem \ref{main thm qstat and qerg measures}(1) in the case \(\zeta=\eta_0\) was established in the final display equation of \cite[Theorem B]{AFGTV-TFPF}, while the invariant measure was identified as the unique relative equilibrium state for the random weight \(g\) in \cite{AFGTV-IVC}.
For other results on random systems, see \cite{nan_2026}; for deterministic open systems, see \cite{LMD,Colletetal00,ColletMS,pollicott_open_2017}.
 \begin{remark}%\label{rmk:1230}
The restriction `$n\ge n_{\om,\zt}$' in Theorem \ref{main thm qstat and qerg measures}, which will repeatedly appear in the following theorems, is not essential: Note that, by the supremum bounds for $\zt_{\om,n}$ and $\zt_{\om,\infty}$, the left-hand side of the above inequality is bounded by
$2\|f\|_\infty \kappa^{-n_{\om,\zt}} \kappa^{n}$ for all $0\le n\le n_{\om,\zt}$. 
Hence, for every $n\ge 0$, we have
\[
\left\vert \zt _{\om,n}(f) -\zeta_{\om,\infty} (f) \right\vert  \le 
\max\left\{\frac{2C_0  (1+C_0\Vert f\Vert  )^2
\Vert \psi_\om\Vert }{ \nu_{\om,\infty}( \psi _\om)  },2\|f\|_\infty \kappa^{-n_{\om,\zt}}\right\}
\max\{1,D_{\sg^n\om}C_{\sg^n\om,\phi}\}\kappa^n.
\]
Notice also that the right-hand side of the above inequality converges to zero $\mas$ since both $D$ and $C_{\phi}$ are tempered.
\end{remark}

Put
\[
\widetilde C_{\omega,\phi}:=\sup_{j\ge0}
D_{\sigma^j\omega}C_{\sigma^j\omega,\phi}\kappa^j.
\]
This random variable is finite and tempered, since $DC_\phi$ is tempered.

\begin{thmx}[Yaglom Limit and Conditionally Invariant Measure]\label{main thm yaglom_limit}
 Assume that the spectral gap condition \eqref{G} holds.  Then, 
 \begin{enumerate}
     \item  
     For every $\zeta \in \mathcal P_{\mathcal B,+}(\nu_c)$ with density $\psi$, $f\in\mathcal B$,  
$n\ge n_{\om,\zt}$ and $\mas~\oio$,
\[
\left\vert T_*^n(\zt_n)_\om(f) - \eta_{\om,0}(f)\right\vert  \le 
 \frac{2C_0(1+\Vert f \Vert _\infty )^2
 \Vert \psi _{\sg^{-n}\om}\Vert }{\nu_{\sg^{-n}\om,\infty} (\psi _{\sg^{-n}\om})} \widetilde C_{\omega,\phi}\kappa ^n;
\]
\item $\eta_0$ is the Yaglom limit of every $\zeta\in\mathcal P_{\mathcal B,+}(\nu_c)$ with random density $\psi$ such that $\|\psi\|$ and $\nu_\infty(\psi)^{-1}$ are tempered with respect to $(\sigma^{-1},m)$;
\item $\eta_0$
 is the unique conditionally invariant random probability measure in $ \mathcal P_{\mathcal B,+}(\nu_c)$.
\end{enumerate}
\end{thmx}

For random expanding interval maps with holes, a Yaglom convergence
estimate corresponding to Theorem~\ref{main thm yaglom_limit}(1)
for $\zeta=\nu_c$ was established in the second displayed equation
of Theorem B of \cite{AFGTV-TFPF}, while Theorem A(4) of the same
paper established the uniqueness of $\eta_0$ among random
conditionally invariant probability measures with bounded-variation
densities.\footnote{Some deterministic open systems admit uncountably
many absolutely continuous conditionally invariant measures;
see \cite{DemersYoung06}.}

Given   $\zeta \in \mathcal P_{\mathcal B,+}(\nu_c)$ with random density $\psi$ and 
    $f,v\in \mathcal B$, we define
    \[
C_{\omega, \zeta,f,v}:=
\frac{2C_0(1+\Vert  f\Vert _{\infty}\Vert  v\Vert _{\infty})\max\{\Vert v\psi_\omega\Vert ||f||_{\infty},\Vert \psi_\om\Vert\} }{\nu_{\om,\infty}(\psi_\om)}.
\]
\begin{thmx}[Conditional Exponential Mixing]\label{thm:c}
Assume that the spectral gap condition \eqref{G} holds.  
Then, 
    \begin{align*}
    \lt|
\zt _{\om,N}\lt( f\circ T^n_\om \cdot v \rt)
	-
	\eta_{\sg^n\om,N-n}(f)   \zt_{\om,\infty}(v)
    \rt|
    \leq 
 C_{\omega, \zeta,f,v} (D_{\sg^n\om} C_{\sg^n\om,\phi})
  \kp^n
    \end{align*}
    for every
$\zeta\in\mathcal P_{\mathcal B,+}(\nu_c)$,
$f,v\in\mathcal B$,
integers $N,n$ satisfying
$n_{\omega,\zeta}\le n\le N$,
and $m$-a.e.~$\omega$.
\end{thmx}
The proofs of Theorems \ref{main thm qstat and qerg measures}, \ref{main thm yaglom_limit}, and \ref{thm:c} appear in Section~\ref{sec: proof of A,B,C}.

\begin{thmx}[Conditional Ergodicity and Conditional Laws of Large Numbers]\label{thm:e1}
Assume that the spectral gap condition \eqref{G} holds.
Then, \begin{enumerate}
\item $\widehat\eta_\infty$
is the unique conditionally ergodic limit over $\cP_{\cB,+}(\nu_c)$.
Moreover, for every $\zeta\in\mathcal P_{\mathcal B,+}(\nu_c)$ and $m$-a.e.~$\oio$, we have
\begin{equation}\label{eq:20261111}
\lim_{n\to\infty}\frac1n\sum_{j=0}^{n-1}
\bigl(T_*^j(\zeta_{n-1})\bigr)_{\sigma^j\omega}
=\overline\eta_\infty,
\end{equation}
where $\overline\eta_\infty$ is the non-random probability measure defined by
\[
\overline\eta_\infty:=\int_\Omega\widehat\eta_{\omega,\infty}\,dm(\omega).
\]
\item For every $\zeta\in\mathcal P_{\mathcal B,+}(\nu_c)$,
every real-valued $f\in L^\infty(\Omega,\cB)$,
every $\delta>0$ and $m$-a.e.~$\oio$, it holds that
\begin{align*}
   \lim_{n\to\infty}\zeta_{\om,n}\left( 
  \Bigl|  \frac{1}{n} \sum_{j=0}^{n-1}f_{\sg^j\om} \circ T^j_\om
    - \int_\Omega\eta_{\omega,\infty}(f_\omega)\,dm(\omega)
    \Bigr| >\delta\right)=0;
\end{align*}
\item For every $\zeta\in\mathcal P_{\mathcal B,+}(\nu_c)$, every real-valued $f\in L^\infty(\Omega,\cB)$,
 $m$-a.e.~$\oio$ and $(\prod_{n=0}^\infty \zeta_{\om,n})$-a.e.~$(x_n)_{n\ge 0} \in \prod_{n=0}^\infty X _{\omega,n}$,
we have 
\begin{align*}
    \lim_{n\to\infty}
\frac{1}{n}\sum_{j=0}^{n-1}f_{\sg^j\om}\circ T^j_\om(x_n)=\int_\Omega\eta_{\omega,\infty}(f_\omega)\,dm(\omega).
\end{align*}
\end{enumerate}
\end{thmx}
Equation \eqref{eq:20261111} is the base-forward version of conditional ergodicity. In the deterministic case, it coincides with the conditionally ergodic limit above.
The proof of Theorem \ref{thm:e1} is given in Section \ref{sec: proof of D}

\begin{remark}
We note that -- similar to closed random dynamics -- with Theorem \ref{thm:e1}(3) one has a notion of a non-random ``physical'' measure.
For comparison, suppose the closed random system admits an ergodic invariant probability measure $\eta_c$. 
Birkhoff applied to $\eta_c$ yields
$$\lim_{n\to\infty}\frac{1}{n}\sum_{j=0}^{n-1} f_{\sigma^j\omega}\circ T_\omega^j (x)=\int_{X} f(\om,x)\ d\eta_c(\om,x)$$
for $\eta_c$-a.e.\ initial condition $(\omega,x)$.
In applications where $\sigma$ denotes evolution of a random environment configuration $\omega$, it is common to only have access to observations that depend on a physical phase space $X$.
This motivates considering the special situation where $f(\om,x)=f(x)$ is independent of $\omega$.
The Birkhoff result becomes 
$$\lim_{n\to\infty}\frac{1}{n}\sum_{j=0}^{n-1} f\circ T_\omega^j(x)=\int_X f(x)\ d\bar{\eta}_c(x),$$
where $\bar{\eta}_c:=\int_\Omega \eta_{\om,c}\, dm(\omega)$,
again for $\eta_c$-a.e.\ initial condition. 
One may think of $\bar{\eta}_c$ as a  measure that is physically observable;  it is the measure exhibited on $X$ by random trajectories starting at $\eta_c$ a.e.\ initial condition.
See \cite{buzzi_absolutely_nodate,FGTM19} for further discussion.

In the open random setting, under the continuity hypothesis of
Proposition~\ref{cor:0201a}, define
$\bar\eta_\infty:=\int_\Omega\widehat\eta_{\omega,\infty}\,dm(\omega)$.
Then Theorem~\ref{thm:e1}(3) gives, for $f\in\mathcal B$,
$$  \lim_{n\to\infty}
    \frac{1}{n} \sum_{j=0}^{n-1}f\circ T_\om^j(x_n)= \int_X f(x) \, d\bar{\eta}_{\infty}(x),$$
for $m$-a.e.\ $\omega$ and $(\prod_{n=0}^\infty \zt_{\om,n})$-a.e.~$(x_n)_{n\ge 0} \in \prod_{n=0}^\infty X _{\omega,n}$.
The probability measure $\bar{\eta}_\infty$ is the measure exhibited on $X$ by long, simultaneously conditioned random trajectories in the sense given above, and the one that is relevant in experiments where only observations depending on $X$ are possible.
\end{remark}

\begin{remark}\label{rmk:integrable-observables}
The restriction to $\cB$-valued observables in Theorem~\ref{thm:e1}(2)--(3)
cannot be replaced by a naive integrability assumption. Consider an example in which
$\widehat\eta_{\omega,\infty}(X_{\omega,\infty})=1$ and
$\nu_{\omega,c}(X_{\omega,\infty})=0$ for $m$-a.e.~$\omega$ (e.g.~examples in Section \ref{s:example}).
If one works with arbitrary
$\widehat\eta_{\omega,\infty}$-integrable functions, the conclusion may fail:
the function
$
    f_\omega=\hat X_{\omega,0} -\hat X_{\omega,\infty}
$
satisfies $\widehat\eta_{\omega,\infty}(f_\omega)=0$, whereas its conditioned Birkhoff
averages under $\zeta_{\omega,n}$ are equal to $1$ for $(\prod_{n=0}^\infty\zeta_{\omega,n})$-a.e.\
point.  On the other hand, if one works with $L^1(\widehat\eta_{\om,\infty})$-classes, the
left-hand sides are not well defined: the same function represents the class
$[0]$ in $L^1(\widehat\eta_{\om,\infty})$, but gives different finite-time conditioned
averages from the representative $0$. 
\end{remark}

\begin{remark}\label{rmk:1230b}
Although we assumed exponential decay in \eqref{eq:stronggap}, as one can see from the proofs, weaker decay rates suffice for the conclusions of
Theorem~\ref{thm:e1}.
In the deterministic case, Theorem~\ref{thm:e1}(1) remains valid if the error in Theorem~\ref{thm:c} decays to zero at any rate, uniformly over $N\ge n$.
In particular, polynomial decay of any positive order suffices.

On the other hand, the proof of Theorem~\ref{thm:e1}(3) relies on the four-point correlation estimate \eqref{eq:qslln-4pt-bound}.
It is enough that the corresponding normalized fourth moments be summable in $N$.
For instance, the same argument works if the product type exponential decay in \eqref{eq:qslln-4pt-bound} is replaced by polynomial decay of order $p>1/2$, with corresponding control of the remaining error term.
Theorem~\ref{thm:e1}(2) then follows from (3).

In the random case, however, the random prefactors in these estimates require additional care.
For Theorem~\ref{thm:e1}(1), the random prefactors in the mixing estimate can still be controlled at the level of Ces\` aro averages by truncation and Birkhoff's ergodic theorem, provided that the corresponding conditioning-limit convergence is available.
For Theorem~\ref{thm:e1}(3), the prefactors must be controlled so that the normalized fourth moments remain summable.

Furthermore, \eqref{eq:stronggap} can be replaced by another type of exponential decay (which also naturally appears in examples, see \cite{AFGTV-TFPF}): there is $\kp\in(0,1)$ such that
\forevery $\mathcal B$-valued random variable $f$, one can find a tempered measurable function $D_f:\Omega\to(0,\infty)$ such that
\begin{align}\label{eq:strongspectralgapGinfinf}
    \norm{ \tilde\cL ^n_\om f  - \nu_{\om,\infty}(f) \phi _{\sg^n\om}}_\infty\leq D_f\kp^n 
\end{align}
for each $n\in\NN$ and $\mas~\oio$.
If one employs \eqref{eq:strongspectralgapGinfinf} instead of \eqref{eq:stronggap}, the estimates in Theorems~\ref{main thm qstat and qerg measures}, \ref{main thm yaglom_limit}, and~\ref{thm:c} remain valid, with the corresponding constants expressed in terms of $D_\psi$ and $D_{v\psi}$.
These estimates also suffice to obtain Theorem~\ref{thm:e1}(1).
 \end{remark}

\subsection{Conditional limit laws}

The second part of our main results presents conditional limit laws for open random dynamical systems.
Limit theorems, such as central limit theorems, large deviation principles, and Berry--Esseen estimates, for closed random dynamical systems have been developed in several settings; see \cite{HH01,DFGTV18b,DFGTV20,CNW,DH-Limit20,AFGTV-TFPF} and references therein.
For open systems the corresponding conditional theory is much less developed.
To the best of our knowledge, the results below provide the first systematic abstract transfer operator treatment of conditional measures and statistical limit laws in a unified setting for deterministic and quenched random open systems.

For the conditional limit theorems below, we additionally require a \emph{uniform} spectral gap condition \eqref{LL}.
Since we will employ a Taylor expansion of an eigenvalue $\lambda(\theta)$ of $\mathcal L(e^{\theta} \cdot )$, called a Nagaev--Guivarc'h type approach (cf.~\cite{Gouezel2015}), unlike the previous part (recall Remark \ref{rmk:1230b}), we need the spectral gap.
This spectral gap is also utilized for applying Crimmins' perturbation theorem for examples in Section \ref{s:example}.

\mylabel{LL}{LL}
\textit{Limit law condition:} 
\begin{enumerate}
\item[\mylabel{LL1}{LL1}] The condition \eqref{G} holds, the constant \(D\) in \eqref{G3} can be chosen uniformly in \(\omega\), and, 
$
    \esssup_{\omega\in\Omega}\lambda_\omega^{-1}K_\omega<\infty 
$ with \(K\) and \(\lambda\) as in \eqref{G0} and \eqref{G1}.
\item[\mylabel{LL2}{LL2}] Either (a) $\Omega$ is a Polish space, $\sg$ is a homeomorphism and $\cL :\Omega \to L(\cB )$ is measurable, where $L(\cB )$ is the space of bounded operators on $\cB$ equipped with the Borel $\sg$-algebra induced by the operator norm, or (b) $\Omega$ is a Lebesgue--Rokhlin probability space, $\cB$ is separable, and the map $\omega \mapsto \cL_\omega h$ is measurable for each $h\in \cB$.
\item[\mylabel{LL3}{LL3}] 
There exist a norm $\vert\cdot\vert$ on $\cB$, a positive integer $n_0$, and positive random variables $\alpha,\beta$ such that $\{h\in\cB:\Vert h\Vert\le1\}$ is compact in $(\cB,\vert\cdot\vert)$, $\int\log\alpha\,dm<0$, $\int\beta\,dm<\infty$, and
    \[
    \Vert \tilde{\mathcal L}^{n_0}_\om h\Vert \le \alpha _\om\Vert h\Vert +\beta _\om\vert h\vert 
    \]
for any $h\in\cB$ and $\mas ~\oio$.
\item[\mylabel{LL4}{LL4}] $\essinf_\om \nu_{\om,c}(\phi_\om)>0$ and $\esssup_\om \| \hat X_{\om,0}\| <\infty$.
 \end{enumerate}

\begin{remark}\label{rmk:5.6}
The assumptions \eqref{LL2} and \eqref{LL3} are included only to justify the spectral perturbation argument for small values of the twist parameter in the proof of Theorem~\ref{conj:1}.
In the deterministic case, this argument is not needed: the spectral gap in \eqref{G}, together with the usual analytic perturbation theorem, gives the required decomposition directly.
Furthermore, in the deterministic case, \eqref{LL1} and \eqref{LL4} also follow from \eqref{B3} and \eqref{G}.

We emphasize that the conditional limit laws below are new already for deterministic open systems.
\end{remark}

In the following, we have two time parameters $N$ and $n$.
For comparison with previous theorems, it is more natural to only consider the case $n=N$, but this formulation is important to obtain the corollaries in Section \ref{s:cors}.

Let $f$
be a $\cB$-valued random variable such that $f_\om$ is real-valued and $\eta_{\om,\infty}(f_\om) = 0$ $m$-a.s.
If the zero $\eta_\infty$-mean condition is not satisfied, one can subtract the random $\eta_\infty$ means (namely, consider $\tilde f_\omega:=f_\om -\eta_{\om,\infty}(f_\om)$).
We define the random dynamical variance of $f$ as  
\begin{equation}\label{eq:20260705c}
\Sigma^2(f)=\int_\Omega\left(\eta_{\omega,\infty}(f_\omega^2)+2\sum_{k=1}^{\infty}\nu_{\sg^k\omega,\infty}\bigl(f_{\sg^k\omega}\, \tilde{\mathcal L}_\omega^k(\phi_\om f_{\omega})\bigr)\right)\,dm(\omega).
\end{equation}
In Proposition \ref{prop:0625a} we prove  that this random variance
exists and is finite under \eqref{LL}.
\begin{thmx}[Conditional Central Limit Theorem]\label{clt}
Assume that the limit law condition \eqref{LL} holds.  
Let $f\in L^\infty(\Omega,\cB)$ be real-valued
with $\eta_{\om,\infty}(f_\om) =0$ $m$-a.s.~and $\Sigma ^2(f) >0$.
With the above notation, \forevery $\zeta\in\mathcal P_{\mathcal B,+}(\nu_c)$, $a\in \mathbb R$ and $m$-a.e. $\oio$,
\begin{equation*}
\lim _{n\to \infty} \sup _{N\ge n}\left|\zeta_{\om,N}\left(\frac{1}{\sqrt n}\sum_{j=0}^{n-1}f_{\sg^j\om} \circ T^j_\om \le a\right) -
\frac{1}{\sqrt{2\pi \Sigma ^2(f)}} \int ^a _{-\infty} \exp\left(-\frac{z^2}{2\Sigma ^2(f)}\right) dz\right|=0
\end{equation*}
\end{thmx}

\begin{thmx}[Conditional Large Deviation Principle]\label{ldp}
Assume that the limit law condition \eqref{LL} holds.  
Let $f\in L^\infty(\Omega,\cB)$ be real-valued
with $\eta_{\om,\infty}(f_\om) =0$ $m$-a.s.~and $\Sigma ^2(f) >0$.
Then, there exist $\delta_0>0$ and {non-random}, nonnegative, continuous, strictly convex function $c:(-\delta _0,\delta_0)\to \mathbb R$ with $c(0)=0$ such that \forevery $\zeta\in\mathcal P_{\mathcal B,+}(\nu_c)$, $\delta \in [0,\delta_0)$ and $\mae~\oio$,
\[
\lim_{n\to\infty}\sup_{N\ge n}
\left|\frac1n\log\zeta_{\omega,N}\left(    \frac1n \sum_{j=0}^{n-1}f_{\sg^j\om} \circ T^j_\om>\delta\right)+c(\delta)\right|
=0.
\]
\end{thmx}

\begin{thmx}[Conditional Berry--Esseen Theorem]\label{be}
Assume that the limit law condition \eqref{LL} holds.  
Let $f\in L^\infty(\Omega,\cB)$ be real-valued
with $\eta_{\om,\infty}(f_\om) =0$ $m$-a.s.~and $\Sigma ^2(f) >0$.
Then, \forevery $\zeta\in\mathcal P_{\mathcal B,+}(\nu_c)$,
there are a 
random variable $C=C_{f,\zt}>0$ and (self-normalising) random variables $(\sigma_n^2(f))_{n\ge 0}$ (defined explicitly in \eqref{eq:0805b}) such that $\frac{\sg_{\om,n}^2(f)}{n}$ converges to $\Sigma^2(f)$ $\mas$ and that for $m$-a.e.~$\oio$,
\[
\sup_{a\in\mathbb{R},N\ge n}\left| \zeta_{\om,N} \left(   \frac{1}{\sqrt{\sg_{\om,n}^2(f)}}\sum_{j=0}^{n-1}f_{\sg^j\om} \circ T^j_\om \leq a \right) - \frac{1}{\sqrt{2\pi }} \int^a_{-\infty} \exp\left(-\frac{z^2}{2}\right)\ dz
\right| 
\leq \frac{C_\om}{\sqrt n}.
\]
\end{thmx}
Since $\sigma_{\omega,n}^2(f)/n\to\Sigma^2(f)$, the central limit
theorem with normalisation $\sigma_{\omega,n}(f)$ is equivalent to
Theorem~\ref{clt}; it also follows immediately from Theorem~\ref{be}.

Theorems \ref{clt}, \ref{ldp}, and \ref{be} are proven in Section \ref{sec: proof of E, F, G}.

\begin{remark}\label{rmk:cLDPtocSLLN}
The conditional Large Deviation Principle (Theorem~\ref{ldp}) immediately implies a conditional Strong Law of Large Numbers in the sense of Theorem~\ref{thm:e1}(3).
More precisely, under the assumptions of Theorem~\ref{ldp} (in particular, $\eta_{\om,\infty}(f_\om)=0$ $\mas$ and $\Sigma^2(f)>0$), 
\forevery $\zt\in\cP_{\cB,+}(\nu_c)$,
$\mae~\oio$ and 
$\bigl(\prod_{n=0}^\infty \zt_{\om,n}\bigr)$-a.e.~$(x_n)_{n\ge 0}\in\prod_{n=0}^\infty X_{\om,n}$,
we have
\[
  \lim_{n\to\infty}
  \frac{1}{n}\sum_{j=0}^{n-1}f_{\sg^j\om}\circ T^j_\om(x_n)=0.
\]

In fact, fix $\omega$ in the full-measure set on which Theorem~\ref{ldp} holds, and let
$
  Z_\om:=\bigotimes_{n\ge 0}\zt_{\om,n}
$
be the product probability measure on $\prod_{n\ge 0}X_{\om,n}$.
For $n,k\in\mathbb N$, define the event
\[
  E_{n,k}:=\left\{(x_j)_{j\ge 0}\in \prod_{j\ge 0}X_{\om,j}\,:\,
  \left|\frac{1}{n}\sum_{j=0}^{n-1}f_{\sg^j\om}\circ T^j_\om(x_n)\right|>\frac{1}{k}\right\}.
\]
Since the $n$th marginal of $Z_\omega$ is $\zeta_{\omega,n}$, we have
\[
  Z_\om(E_{n,k})
  =\zt_{\om,n}\!\left(\left|\frac{1}{n}\sum_{j=0}^{n-1}f_{\sg^j\om}\circ T^j_\om\right|>\frac1k\right).
\]
Applying Theorem~\ref{ldp} to $f$ and to $-f$ (note $\Sigma^2(-f)=\Sigma^2(f)>0$), 
choose $\delta_0>0$ smaller than both admissible deviation ranges, let $c_f$ and $c_{-f}$ be their rate functions, and set $\hat{c}(u):=\min\{c_f(u),\,c_{-f}(u)\}>0$ for $u\in(0,\delta_0)$.
By Theorem~\ref{ldp} (with $N=n$), for any $k$ with $k^{-1}\in(0,\delta_0)$,
\[
  \limsup_{n\to\infty}\frac{1}{n}\log Z_\om(E_{n,k})\le -\hat{c}(k^{-1})<0.
\]
Therefore, $\sum_{n\ge 1}Z_\om(E_{n,k})<\infty$ for every integer $k> \delta_0^{-1}$. Hence, it follows from the Borel--Cantelli lemma that for every $k> \delta_0^{-1}$,
\[
Z_\om\big(\limsup_{n\to\infty}E_{n,k}\big)=0.
\]
Take $E:=\bigcap_{k> \delta_0^{-1}}(\limsup_{n\to\infty}E_{n,k})^c$, then $Z_\omega(E)=1$ because $E$ is a countable intersection of the $Z_\omega$-full measure sets.
For every $(x_n)_{n\ge 0}\in E$, one has
$\bigl|\frac{1}{n}\sum_{j=0}^{n-1}f_{\sg^j\om}\circ T^j_\om(x_n)\bigr|\le \frac1k$ for all $k> \delta_0^{-1}$ and all $n\ge n_k$ with some integers $n_k$.
Hence we conclude
\[
  \frac{1}{n}\sum_{j=0}^{n-1}f_{\sg^j\om}\circ T^j_\om(x_n)\to 0
  \qquad\text{$Z_\om$-a.s.}
\]

To remove the assumption $\eta_{\om,\infty}(f_\om)=0$, 
one replaces $f$ by $\tilde{f}_\om:=f_\om-\eta_{\om,\infty}(f_\om)$ (refer to the proof of Theorem~\ref{thm:e1}(3)):
since $(\sg,m)$ is ergodic, Birkhoff's ergodic theorem gives
$\frac1n\sum_{j=0}^{n-1} \eta_{\sg^j\om,\infty}(f_{\sg^j\om}) \to\int_\Omega\eta_{\omega,\infty}(f_\omega)\,dm(\omega)$ for $\mae~\oio$.
Under the other hypotheses of Theorem~\ref{ldp} for $\tilde f$, including $\Sigma^2(\tilde f)>0$, this gives the conclusion of Theorem~\ref{thm:e1}(3).
\end{remark}

\subsection{Limit laws on the surviving set \texorpdfstring{$X_\infty$}{X-infinity}}\label{s:cors}

In this subsection, we record conventional limit theorems on the surviving set $X_\infty$.
A route to such results is to consider the coded Markov dynamics on survival (called Q-process in the context of Markov processes \cite{BreyerRoberts99,zhang_quasi-stationarity_2014, he_quasi-ergodic_2019}), and then to establish a spectral gap for the closed transfer operator and apply the (classical) Nagaev--Guivarc'h method.
See the monograph \cite{pollicott_open_2017} and references therein.
In the present paper, however, we do \emph{not} assume a spectral gap for the closed transfer operator describing the restricted random dynamics $T:X_\infty\to X_\infty$ itself.
Instead, we work under the assumption \eqref{LL} for the open transfer operators, together with the boundary and continuity conditions stated below, and obtain limit theorems on $X_\infty$ via 
the uniform-in-time $n$ estimates for the $N$-conditioned measures in the theorems of the previous sections.
Although limit theorems on $X_\infty$ are not the main focus of the present paper,\footnote{
For example, under the hypotheses of Proposition~\ref{cor:0201a}, $\widehat\eta_\infty$ is ergodic.
For each fixed $f\in L^1(\widehat\eta_\infty)$, Birkhoff's ergodic theorem gives
\begin{equation}\label{eq:0104c}
\lim_{n\to\infty}\frac1n\sum_{j=0}^{n-1}
f\circ T^j(\om,x)
=\int f\,d\widehat\eta_\infty
\end{equation}
for $\widehat\eta_\infty$-a.e.~$(\omega,x)$.
However, in our examples (see Section~\ref{s:example}), the basin of attraction of $\widehat\eta_\infty$ (i.e.~the set of points $(\omega,x)$ satisfying \eqref{eq:0104c} for every $f\in C(\Xi)$) may have $\zeta$-measure zero for every $\zeta\in\mathcal P_{\mathcal B,+}(\nu_c)$ (especially when $\nu_{\omega,c}$ is Lebesgue measure).
In contrast, Theorem~\ref{thm:e1}(3) gives convergence for the admissible observables under the conditioned measures
$\zeta_{\omega,n}$.
This observation is particularly meaningful when
$\zeta_\om(X_{\om,n})$ does not decay rapidly,
which indeed occurs in our examples.}
we include these results since the mechanism may be of independent interest.

The following proposition may be regarded as a generalization of our previous results in \cite{AFGTV-TFPF} (Theorem A(3) and Theorem C in the ``large holes'' case, and Theorem J for the ``small holes'' case) to an abstract setting, with a slightly different proof.

\begin{proposition}[Invariance and Exponential Decay of Correlations on $(X_\infty,\widehat\eta_\infty)$]\label{cor:0201a}
 Assume that \eqref{G} holds and that every function in $\cB$ is $\widehat\nu_{\om,\infty}$-a.e.~continuous for $m$-a.e.~$\oio$.
 Then, the following holds.
 \begin{itemize}
    \item[(1)] For every
$\zeta\in\mathcal P_{\mathcal B,+}(\nu_c)$ and every $f\in\mathcal B$,
\begin{equation}\label{eq:zeta-functional-measure-identification}
\zeta_{\omega,\infty}(f)
=
\widehat\zeta_{\omega,\infty}(f)
\qquad\text{for $m$-a.e.~$\oio$.}
\end{equation}
        \item[(2)]  $\widehat\eta_\infty$ is the unique ergodic, 
 $T$-invariant probability measure on $X_\infty$ being absolutely continuous with respect to $\widehat\nu_\infty$.
 \item[(3)]
For  every $\zeta \in \mathcal P_{\mathcal B,+}(\nu_c)$, every $f,v\in \cB$,  
every integer $n\ge n_{\om,\zt}$ and $\mas ~\oio$, 
    \begin{align*}
 \lt|
\widehat\zt _{\om,\infty}\lt( f\circ T^n_\om \cdot v \rt)
	-
	\widehat\eta_{\sg^n\om,\infty}(f)   \widehat\zt_{\om,\infty}(v)
    \rt|
    \leq 
C_{\omega, \zeta,f,v} (D_{\sg^n\om} C_{\sg^n\om,\phi})
  \kp^n.
    \end{align*}
 \end{itemize}
\end{proposition}
The proof of Proposition~\ref{cor:0201a} is postponed to Appendix~\ref{app:survivor-invariance}.

In our previous work \cite{AFGTV-TFPF}, limit theorems for the conditioned dynamics on the surviving set were not obtained.
See Section~2.6 of \cite{AFGTV-TFPF}, where the authors explicitly note that they could not extend their closed-system limit theorems to the open systems considered there.

\begin{corollary}[Large Deviation Principle and Berry--Esseen Theorem on $(X_\infty,\widehat\eta_\infty)$]\label{cor:CLTLDPXinfty}
Assume that the limit law condition \eqref{LL} holds.  
Let $f\in L^\infty(\Omega,\cB)$ be real-valued and satisfy $\eta_{\om,\infty}(f_\om)=0$ $m$-a.s.\ and $\Sigma^2(f)>0$.
We also assume that \eqref{eq:survivor-boundary-null} holds and $f_{\sg^n\om}\circ T_\om^n$ is $\widehat\nu_{\om,\infty}$-a.e.~continuous for every $n\ge 0$ and $m$-a.e.~$\oio$.
Then, the following hold for every $\zeta\in\mathcal P_{\mathcal B,+}(\nu_c)$.
\begin{itemize}
    \item[(1)] There is a random variable $C>0$ such that for every $n\ge n_{\om,\zt}$ and $\mas ~\oio$,
\[
\sup_{a\in\mathbb{R}}\left| \widehat\zeta_{\om,\infty}\left(   \frac{1}{\sqrt{\sg_{\om,n}^2(f)}}\sum_{j=0}^{n-1}f_{\sg^j\om} \circ T^j_\om \leq a \right) - \frac{1}{\sqrt{2\pi }} \int^a_{-\infty} \exp\left(-\frac{z^2}{2}\right) dz
\right| 
\leq \frac{C_\om}{\sqrt n}. 
\]
    \item[(2)] There exist $\delta_0>0$ and nonnegative, continuous, strictly convex function $c:(-\delta _0,\delta_0)\to \mathbb R$ with $c(0)=0$ such that \forevery $\zeta\in\mathcal P_{\mathcal B,+}(\nu_c)$, $\delta \in [0,\delta_0)$ and $\mas~\oio$,
\[
\lim _{n\to \infty}  \frac{1}{n} \log \widehat\zeta_{\om,\infty}\left( \frac{1}{n}\sum_{j=0}^{n-1}f_{\sg^j\om} \circ T^j_\om  > \delta  \right) =  -c(\delta ).
\]
\end{itemize} 
\end{corollary}
\begin{proof}
Fix $\omega$ in a full-measure set on which Theorems~\ref{main thm qstat and qerg measures}, \ref{ldp}, and~\ref{be} hold.
By the absolute continuity of Theorem~\ref{main thm qstat and qerg measures}(4),
the assumed almost everywhere continuity of $f$ also holds with respect to $\widehat\zeta_{\omega,\infty}$.
Furthermore, condition~\eqref{eq:survivor-boundary-null} ensures that $\widehat\zeta_{\omega,\infty}(X_{\omega,\infty})=1$.

We first prove item (1). Fix \(n\ge1\) so large that
\(\sigma_{\omega,n}^2(f)>0\).
By Theorem~\ref{main thm qstat and qerg measures} the measures \(\zeta_{\omega,N}\) converge weakly to \(\widehat\zeta_{\omega,\infty}\), and  by the assumption 
\(S_{\omega,n}f\) is \(\widehat\zeta_{\omega,\infty}\)-a.e.~continuous.
Hence, Billingsley's mapping theorem\footnote{If a sequence $\varrho_n$ of probability measures converges weakly to a probability measure $\varrho_\infty$ and the set of discontinuity points of a measurable map $h$ is a $\varrho_\infty$-null set, then the sequence $h_*\varrho_n$ of the pushforward measures converges weakly to the pushforward measure $h_*\varrho_\infty$. Refer to \cite[Theorem 2.7]{Billingsley1999}. \label{f:0706}} gives weak convergence of the pushforward measures
under
$(\sigma_{\omega,n}^2(f))^{-1/2}\sum_{j=0}^{n-1}f_{\sg^j\om} \circ T^j_\om.$
Consequently, for every continuity point \(b\) of the distribution of
$(\sigma_{\omega,n}^2(f))^{-1/2}\sum_{j=0}^{n-1}f_{\sg^j\om} \circ T^j_\om$
with respect to \(\widehat\zeta_{\omega,\infty}\), we have
\begin{equation}\label{eq:BE-Xinfty-cont-point}
\zeta_{\omega,N}
\left(A_{n,b}\right)
\to
\widehat\zeta_{\omega,\infty}
\left(A_{n,b}\right),\qquad \text{where}\quad A_{n,b}:= \left\{  \frac{1}{\sqrt{\sg_{\om,n}^2(f)}}\sum_{j=0}^{n-1}f_{\sg^j\om} \circ T^j_\om \leq b\right\}.
\end{equation}

Let \(a\in\mathbb R\) and let \(\varepsilon>0\).  Since a distribution
function has at most countably many discontinuities, we may choose
continuity points \(a_-<a<a_+\) of the limiting distribution such that
\[
\frac{1}{\sqrt{2\pi}}
\int_{a_-}^{a_+}e^{-z^2/2}\,dz
<\varepsilon .
\]
Using \eqref{eq:BE-Xinfty-cont-point} at \(b=a_+\) and Theorem~\ref{be}, we obtain
\begin{align*}
\widehat\zeta_{\omega,\infty}\left(A_{n,a}\right)
\le
\widehat\zeta_{\omega,\infty}
\left(A_{n,a_+}\right)
&=
\lim_{N\to\infty}
\zeta_{\omega,N}
\left(A_{n,a_+}\right)\\
&\le
\frac{1}{\sqrt{2\pi}}\int_{-\infty}^{a_+}e^{-z^2/2}\,dz
+
\frac{C_\omega}{\sqrt n}
\le
\frac{1}{\sqrt{2\pi}}\int_{-\infty}^{a}e^{-z^2/2}\,dz
+
\varepsilon
+
\frac{C_\omega}{\sqrt n}.
\end{align*}
Similarly, using \eqref{eq:BE-Xinfty-cont-point} at \(b=a_-\) and Theorem~\ref{be}, we obtain
\begin{align*}
\widehat\zeta_{\omega,\infty}\left(A_{n,a}\right)
\ge
\widehat\zeta_{\omega,\infty}\left(A_{n,a_-}\right)
&=
\lim_{N\to\infty}\zeta_{\omega,N}\left(A_{n,a_-}\right)\\
&\ge
\frac{1}{\sqrt{2\pi}}\int_{-\infty}^{a_-}e^{-z^2/2}\,dz
-
\frac{C_\omega}{\sqrt n}
\ge
\frac{1}{\sqrt{2\pi}}\int_{-\infty}^{a}e^{-z^2/2}\,dz
-
\varepsilon
-
\frac{C_\omega}{\sqrt n}.
\end{align*}
Letting \(\varepsilon\to 0\), we get the conclusion of item (1).

We next show the upper bound in item (2). Fix \(\delta\in[0,\delta_0)\), \(\epsilon>0\) and $n$.
As in the proof of item (1), by Theorem~\ref{main thm qstat and qerg measures} and the assumption,
we can apply Billingsley's mapping theorem
to give the weak convergence of
the pushforward measures of $\zeta_{\omega,N}$ under $\frac{1}{n}\sum_{j=0}^{n-1}f_{\sg^j\om} \circ T^j_\om$.
Since \((\delta,\infty)\) is open, the Portmanteau theorem gives
\[
\widehat\zeta_{\omega,\infty}\left(B_{n,0}\right)
\le
\liminf_{N\to\infty}\zeta_{\omega,N}\left(B_{n,0}\right),
\qquad \text{where}\quad B_{n,0}:= \left\{  \frac{1}{n}\sum_{j=0}^{n-1}f_{\sg^j\om} \circ T^j_\om >\delta \right\}.
\]
Therefore, by Theorem \ref{ldp}, there is an integer $n_\varepsilon$ such that if $n\ge n_\varepsilon$, then we obtain
\[
\widehat\zeta_{\omega,\infty}\left(B_{n,0}\right)
\le
e^{-(c(\delta)-\epsilon)n}.
\]
This estimate gives the desired upper bound after logarithmic normalization, passage to the limsup,
and letting \(\epsilon\to 0\).

We now prove the lower bound of item (2).  Fix
\(\delta'\in(\delta,\delta_0)\). 
Since \([\delta',\infty)\) is closed, the Portmanteau theorem gives
\[
\widehat\zeta_{\omega,\infty}\left(B_{n,1}\right)
\ge
\limsup_{N\to\infty}\zeta_{\omega,N}\left(B_{n,1}\right),
\qquad
\text{where}\quad B_{n,1}:= \left\{  \frac{1}{n}\sum_{j=0}^{n-1}f_{\sg^j\om} \circ T^j_\om \ge \delta '\right\}.
\]
Let $B_{n,2}:=\{\frac{1}{n}\sum_{j=0}^{n-1}f_{\sg^j\om} \circ T^j_\om >\delta'\}$.
For all sufficiently large $n$, Theorem~\ref{ldp} gives
\begin{align*}
\widehat\zeta_{\omega,\infty}(B_{n,0})
\ge\widehat\zeta_{\omega,\infty}(B_{n,1})
 \ge\limsup_{N\to\infty}\zeta_{\omega,N}(B_{n,1})
\ge\limsup_{N\to\infty}\zeta_{\omega,N}
 (B_{n,2})
 \ge e^{-(c(\delta')+\epsilon)n}.
\end{align*}
Divide logarithms by $n$, then let $\epsilon\downarrow0$ and
$\delta'\downarrow\delta$. Continuity of $c$ gives the desired lower
bound. Together with the upper bound, this proves item~(2).
\end{proof}
We emphasize that, as seen in the above proof, conditioning limits (Theorem \ref{main thm qstat and qerg measures}) play a crucial role in our approach to show limit theorems on $X_\infty$.

\begin{remark}\label{rem:P-examples}
The $\widehat\nu_{\omega,\infty}$-a.e.~continuity hypothesis in Corollary \ref{cor:CLTLDPXinfty} is automatic in several standard situations.  
Suppose that, for $f\in L^\infty(\Omega,\cB)$, both $T_\omega$ and $f_\omega$ are continuous for $m$-a.e.~$\omega$.
Then, for every $n\ge0$, $f_{\sigma^n\omega}\circ T_\omega^n$ is continuous for $m$-a.e.~$\omega$.
This includes, for example, the finite subshifts of finite type with Lipschitz observables in Section~\ref{s:6.4}. In this setting, the finite-time survivor sets are clopen, so \eqref{eq:survivor-boundary-null} also holds.

For the one-dimensional expanding maps of Section~\ref{s:6.2}, with $\mathcal B=\BV(I)$, the discontinuity set of (the regulated representative of) $f_{\sigma^n\omega}$ is at most countable. 
Moreover, $T_\omega^n$ has at most countably many branch boundaries and is continuous and injective on each of at most countably many branches.
The only possible discontinuities of $f_{\sigma^n\omega}\circ T_\omega^n$ are therefore the branch boundaries of $T_\omega^n$ and points whose image under $T_\omega^n$ is a discontinuity of $f_{\sigma^n\omega}$. 
The first set is countable. 
On each branch, injectivity implies that the inverse image of a countable set is countable; since there are at most countably many branches, the second set is also countable. 
Therefore, if $\widehat\nu_{\omega,\infty}$ is non-atomic, then the stated continuity hypotheses hold. 
Non-atomicity holds when the hole is sufficiently small; see \cite[Lemma 1.10.11]{AFGTV-TFPF} and Section~\ref{s:6.2}.

For the multidimensional expanding setting, almost everywhere continuity is established under a very mild condition in Appendix~\ref{app:multidimensional-ae-continuity}, in particular Proposition~\ref{prop:ValphaIdentification}.
\end{remark}

\section{Examples}\label{s:example}

Our Assumptions \eqref{B1}--\eqref{B2}, \eqref{C0}--\eqref{C1}, and \eqref{M1}--\eqref{M2} 
have been shown to hold for several large classes of random dynamical systems, for example random interval maps \cite{AFGTV20,AFGTV-IVC}, random subshifts \cite{Bogenschutz_RuelleTransferOperator_1995a,mayer_countable_2015}, random distance expanding maps \cite{mayer_distance_2011}, random polynomial systems \cite{Bruck_Generalizediteration_2003}, and random transcendental maps \cite{mayer_random_2018}. 
In the subsections below, we verify these assumptions, along with \eqref{B1}, \eqref{G}, and \eqref{LL}, for specific classes of random map cocycles.

\begin{remark}\label{rem:eqclass}
In the abstract results above, the Banach space $\cB$ consists of functions on $\Xi$.
In Sections~\ref{s:6.2} and~\ref{s:multidimensional-expanding}, however,
$\BV([0,1])$ and $V_\alpha(X)$ are spaces modulo Lebesgue-a.e.\ equality;
in those sections we continue to write $f$ for a class $[f]$.
The operator-theoretic statements are understood on these quotient spaces.
In particular, nonsingularity with respect to Lebesgue measure makes the class of $f\circ T_\omega$ independent of the representative of $f$.

In these two examples,
$\nu_{\omega,\infty}=\widehat\nu_{\omega,\infty}$ on $\cB\cap C(\Xi)$ means
\[
   \nu_{\omega,\infty}([g])
   =\widehat\nu_{\omega,\infty}(g)
   =\int_\Xi g\,d\widehat\nu_{\omega,\infty}
\]
for every $g\in C(\Xi)$ such that $[g]\in\cB$.  This is precisely
Theorem~\ref{main thm qstat and qerg measures}(2) with $\zeta=\nu_c$.
An identity on all of $\cB$ additionally requires a convention for
integrating an equivalence class against the representing measure.

For $\cB=\BV([0,1])$, we use representatives $g$ satisfying $\var(g)<\infty$, where $\var(g)$ is the usual variation.
Every $\BV([0,1])$ class has such a representative by the definition of the quotient variation (cf.~\cite{rychlik_bounded_1983}).  Moreover, $\var(g)<\infty$ implies that $g$ is regulated, that is, it has finite one-sided limits at every point, and hence that $g$ has at most countably many discontinuities.
If $g$ and $h$ are two such representatives of the same class, then $g-h=0$ Lebesgue-a.e.\ and $\var(g-h)<\infty$.  Since $g-h$ vanishes at every continuity point, $\{g\ne h\}$ is at most countable.  
Consequently, when $\widehat\nu_{\omega,\infty}$ is non-atomic,
$  \int g\,d\widehat\nu_{\omega,\infty} =\int h\,d\widehat\nu_{\omega,\infty}$.
Thus, the identity $\nu_{\omega,\infty}=\widehat\nu_{\omega,\infty}$
on $\mathrm{BV}([0,1])$ is understood as
\[
\nu_{\omega,\infty}([f])
=\int g\,d\widehat\nu_{\omega,\infty}
\]
for every $[f]\in\mathrm{BV}([0,1])$ and one, equivalently every, representative
$g\in[f]$ with $\operatorname{var}(g)<\infty$.

For $\cB=V_\alpha(X)$, finite seminorm does not give an analogous convention. 
Indeed, if $N\subset X$ is Lebesgue-null, then $[\hat N]=[0]$ and
\[
   |\hat N|_\alpha=0,
\]
because $\operatorname{osc}(\hat N,B_\varepsilon(x))=0$ for every
$x\in X$ and $0<\varepsilon\le\varepsilon_0$ (refer to Section \ref{s:multidimensional-expanding} for these notations). 
Nevertheless, $\widehat\nu_{\omega,\infty}(\hat N) =\widehat\nu_{\omega,\infty}(N)$ may be positive.
Thus $0$ and $\hat N$ are representatives of the same $V_\alpha(X)$ class, both with finite seminorm, but their integrals against $\widehat\nu_{\omega,\infty}$ may differ.
We therefore fix a canonical representative: for a real-valued class
$[f]\in V_\alpha(X)$, choose any measurable representative $f$ and put
\[
   f_+(x):=\lim_{r\downarrow0}\esssup_{B(x,r)}f,
   \qquad
   f_-(x):=\lim_{r\downarrow0}\essinf_{B(x,r)}f,
   \qquad
   f^\#(x):=\frac{f_+(x)+f_-(x)}2.
\]
For complex-valued $f$, we apply this construction to the real and imaginary parts separately.
The functions $f_+$ and $f_-$ depend only on $[f]$,
and $f^\#$ is a Borel representative of this class.
Thus, the identity $\nu_{\omega,\infty}=\widehat\nu_{\omega,\infty}$
on $V_\alpha(X)$ is understood as
\[
\nu_{\omega,\infty}([f])
=\int f^\#\,d\widehat\nu_{\omega,\infty}
\]
for every $[f]\in V_\alpha(X)$.
The midpoint is only a convenient fixed choice.  More generally, if
$\widetilde f$ is any Borel function satisfying
\[
   f_-(x)\le \widetilde f(x)\le f_+(x),\qquad x\in X,
\]
then $\widetilde f=f$ Lebesgue-a.e.\ and $\widetilde f$ is continuous at
every $x$ for which
\[
   \lim_{r\downarrow0}\operatorname{osc}(f,B(x,r))=0.
\]
Indeed, when $y$ is sufficiently close to $x$, both $\widetilde f(x)$ and $\widetilde f(y)$ lie between the essential infimum and supremum of $f$ on the same small ball centred at $x$. 
Hence the discontinuity estimates in Appendix~\ref{app:multidimensional-ae-continuity} apply to every such choice.

The same interpretation applies to $\zeta_{\omega,\infty}$ and $\widehat\zeta_{\omega,\infty}$.
\end{remark}

\subsection{Illustrative example:  A deterministic map with a hole}
\begin{figure}
    \centering
\includegraphics[width=0.5\linewidth]{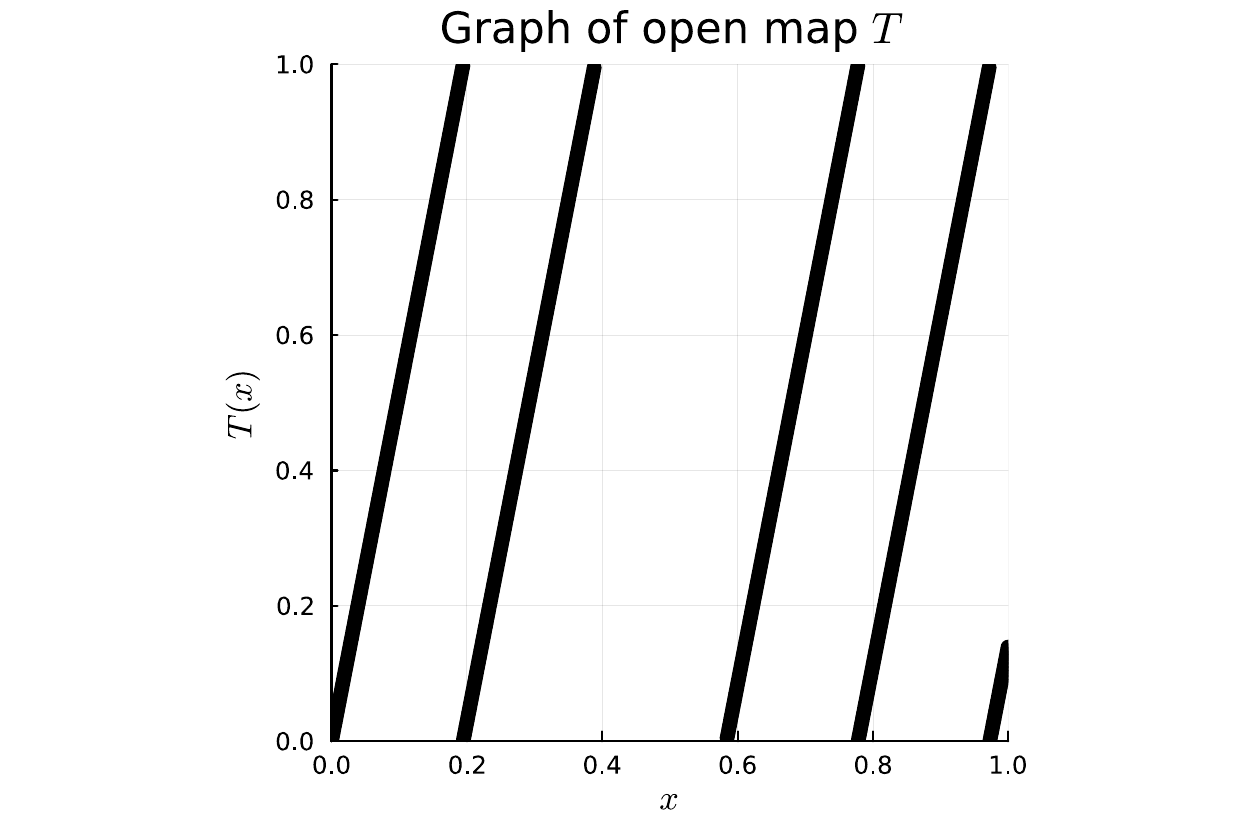}
\caption{Graph of $T(x)=\beta x\pmod 1$, where $\beta=\pi+2$. The hole $H$ is the interval $[2/\beta,3/\beta]$.}
\label{fig:openmap}
\end{figure}

We illustrate the various conditional measures with a deterministic map.
We consider the $\beta$-map shown in Figure \ref{fig:openmap}.
We use the geometric weight $g=1/|T'|=1/\beta$ and take
$\nu_c=\mathrm{Leb}$, so that $\lambda_c=1$.
Conditions \eqref{M1} and \eqref{M2} are trivial, conditions \eqref{B1}, \eqref{B2}, and \eqref{B3} hold for $\mathcal{B}=\BV([0,1])$, and \eqref{C0} and \eqref{C1} are standard.
Condition \eqref{G} follows from Theorem A \cite{LMD}, after checking the conditions 0, 1, and 2 in that paper.

Condition 0 is trivial. 
Because $T$ contains some full branches, Condition~1 holds.
For condition 2, we refer to the sufficient condition ($\widehat{Q1}$) in \cite{AFGTV-TFPF}. 
We seek upper bounds of the two terms on the LHS of ($\widehat{Q1}$) and a lower bound for the RHS of ($\widehat{Q1}$), and then ask that the sum of these upper bounds is less than the lower bound, which implies ($\widehat{Q1}$).
The first term in ($\widehat{Q1}$) is $\log 1/\beta$, while by Proposition 1.15.3 \cite{AFGTV-TFPF} $\lim_n \frac{1}{n}\log\zeta^{(n)}\le \log(\zeta^{(1)}+2)=\log(1+2)$.
A lower bound for $\lim_n \frac{1}{n}\log\rho^n$ is $\log\inf\mathcal{L}\widehat{[0,1]}=\log 4 + \log(1/\beta)$.
Thus to satisfy ($\widehat{Q1}$) we require $\log(1/\beta)+\log(3)<\log(4)+\log(1/\beta)$, which holds.
As explained in Remark~\ref{rmk:5.6}, conditions \eqref{G} and \eqref{B3} suffice for the conditional limit laws in this deterministic setting.

\begin{figure}
\centering
\includegraphics[width=\linewidth]{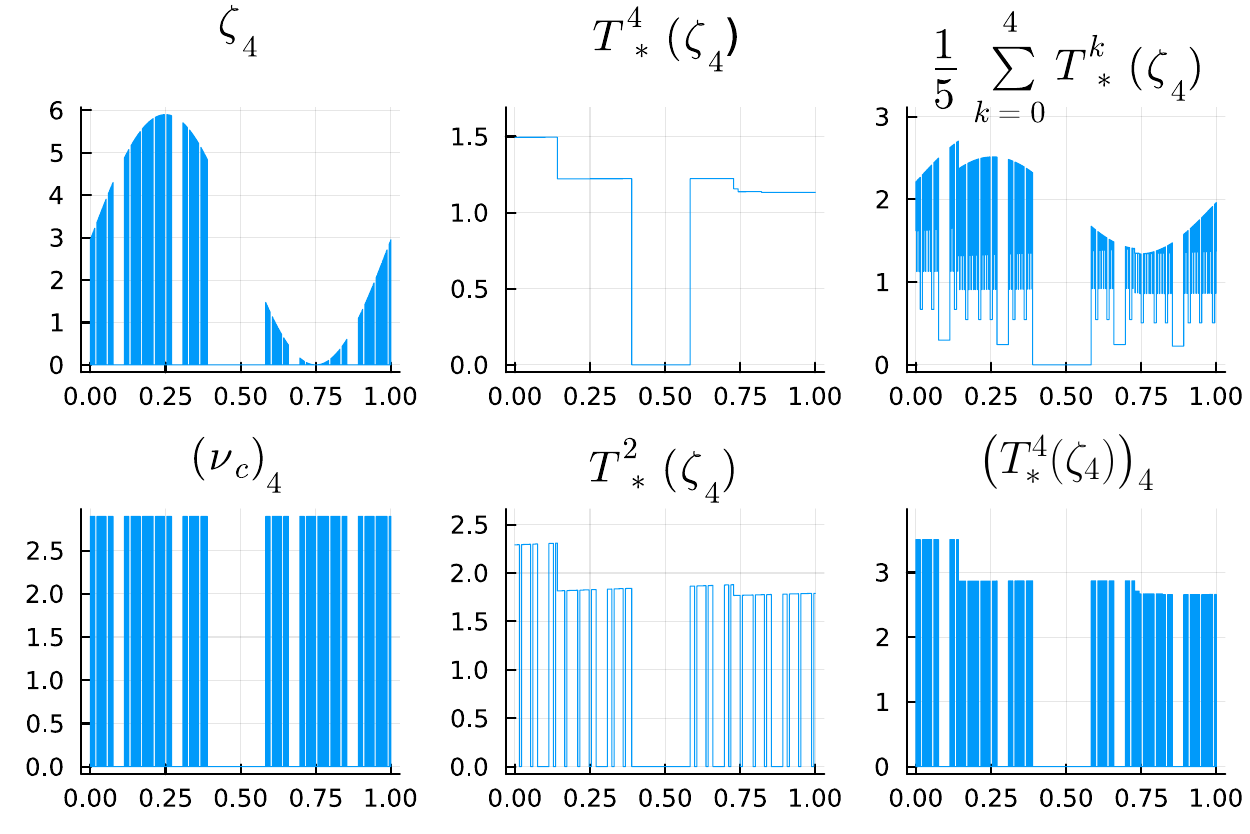}
    \caption{Densities of various conditioned measures arising from an initial probability measure $\zeta$ with density with respect to $\nu_c={\rm Leb}$ given by  $\frac{d\zeta}{d\nu_c}(x)=\sin(2\pi x)+1$.
    \textit{Upper left:} The conditioned measure $\zeta_4$ ($\zeta$ conditioned to the 4-survivor set of $T$).
    By Theorem \ref{main thm qstat and qerg measures}, as we increase the number 4, we converge exponentially quickly with rate $\kappa$ to the conditioning limit $\zeta_\infty$.
     \textit{Lower left:} The same conditioning process applied to the closed conformal measure $\nu_c={\rm Leb}$. The conditioning limit is $\nu_\infty$, the open conformal measure.
    \textit{Upper centre:} Partial Yaglom limit, conditioning $\zeta$ to the 4-survivor set and pushing forward 4 steps with $T_*$. By Theorem \ref{main thm yaglom_limit}, as the number 4 increases, we converge exponentially quickly with rate $\kappa$ to the Yaglom limit $\eta_0$.
    \textit{Lower centre:} Partial conditional exponential mixing. By Theorems~\ref{main thm qstat and qerg measures}(1) and~\ref{thm:c}, by pushing half as much as we condition, we converge exponentially quickly with rate $\sqrt{\kappa}$ to the conditionally ergodic measure $\eta_\infty$.
    \textit{Upper right:} Partial conditional ergodic limit of $\zeta$. By Theorem \ref{thm:e1}(1) as we increase the number 4, we converge to the conditionally ergodic measure $\eta_\infty$.
        \textit{Lower right:} A combination of the partial Yaglom limit of $\zeta$ and the partial conditioning limit of $\nu_c$. By Theorems \ref{main thm qstat and qerg measures} and \ref{main thm yaglom_limit} we converge exponentially with rate $\kappa$ to the conditionally ergodic measure $\eta_\infty$.}\label{fig:conditioned}
\end{figure}

The various conditioned measures associated with our conditioning limit theorems are illustrated in Figure \ref{fig:conditioned}.
In this deterministic example, Theorems~\ref{main thm qstat and qerg measures}(1) and~\ref{thm:c} give, for each fixed $f\in\mathcal B$,
\[
\left|(T_*^{\lfloor N/2\rfloor}\zeta_N)(f)-\eta_\infty(f)\right|
=O\bigl((\sqrt\kappa)^N\bigr).
\]
The lower-centre panel of Figure \ref{fig:conditioned} shows the case $N=4$, namely $T_*^2(\zeta_4)$.


\subsection{Random interval maps with small random holes}\label{s:6.2}

We show that a class of random open interval maps satisfy our general hypotheses.
In the setting of \cite{AFGTV20}, where $\mathcal{B}=\BV([0,1])$, the space of complex-valued functions of bounded variation on the unit interval, \eqref{B2} and \eqref{B3} hold.
Since \eqref{C1} implies that $\BV$ is dense in $L^1(\nu_{\omega,c})$ for a.e.\ $\omega$, we obtain \eqref{B1}. 
The standard property of denseness of $\BV$ in $C([0,1])$ provides \eqref{B1}.
Conditions \eqref{C0} and \eqref{M1} are satisfied by construction.
To verify \eqref{M2} and \eqref{LL2}(a), assume that $\Omega$ is Polish, $\sigma:\Omega\to\Omega$ is a homeomorphism, and the triples $(T_\omega,g_\omega,H_\omega)$ take at most countably many values, with measurable level sets.
Then both the closed and open transfer operator cocycles are operator-norm Borel measurable.
If we further assume we are in the class of random weighted covering systems (Definition 2.18 \cite{AFGTV20}, which asks for nonuniform random expansion, nonuniform random covering, and nonuniformly random contracting weight), then Theorem 2.19 \cite{AFGTV20} provides the existence of a closed random conformal measure satisfying \eqref{C1}. 
The covering condition may be replaced with a stronger random contracting weight condition to achieve similar results (Main Theorem \cite{AFGTV-IVC}).

Lemma 2.5.10 \cite{AFGTV-TFPF} lists a series of hypotheses for random piecewise expanding interval maps $\{T_\omega\}_{\omega\in\Omega}$.
Each $T_\omega$ should be surjective, piecewise monotone, and nonsingular with respect to Lebesgue.
Furthermore, there are uniform-in-$\omega$ upper bounds on $|T'_\omega|$ and the number of inverse branches of $T_\omega$. 
There are uniform-in-$\omega$ bounds on the random weight functions $\{g_\omega\}_{\omega\in\Omega}$ and on the number of connected components of the random holes $H_{\omega}$.
The random map cocycles are uniformly covering, and the random weight $g_{\omega}$ is uniformly contracting.
We refer the reader to Lemma 2.5.10 \cite{AFGTV-TFPF} for details.
For such systems, if the random holes are sufficiently small, Lemma 2.5.10 \cite{AFGTV-TFPF} provides the existence of open random conformal measure $\nu_\infty$ \eqref{G1} and Theorem 2.5.12 \cite{AFGTV-TFPF} guarantees the existence of a $\mathcal{B}$-valued random variable $\phi$ that is uniformly bounded above.
The essential lower bound of $\nu_c(\phi)$ proved below gives the temperedness required in \eqref{G2}.
The exponential decay in \eqref{G3} is also provided by Theorem 2.5.12 \cite{AFGTV-TFPF}, now with the factor $D$ being a constant.
Thus, we also satisfy \eqref{LL1}.
For \eqref{LL3}, take $|u|=\|u\|_{L^1(\mathrm{Leb})}$.
By \cite[(2.5.19)--(2.5.20)]{AFGTV-TFPF}
and the subsequent $L^1$ estimate,
there exist constants $A,B\ge1$ and $q\in(0,1)$
such that
\[
\|(\lambda_{\omega,c}^{n})^{-1}\mathcal L_\omega^n u\|_{\mathrm{BV}}
\le Aq^n\|u\|_{\mathrm{BV}}
+B^n\|u\|_{L^1(\mathrm{Leb})},
\]
where $\lambda_{\omega,c}^{n}:=\lambda_{c,\sg^{n-1}\om}\cdots\lambda_{c,\om}$.
Moreover, the uniform perturbation estimates in
\cite[Lemma~2.5.10]{AFGTV-TFPF} give
$
R:=\operatorname*{ess\,sup}_\omega
\frac{\lambda_{\omega,c}}{\lambda_\omega}
\longrightarrow1
$
as the holes shrink uniformly.
Taking the holes sufficiently small that $Rq<1$,
we obtain
\[
\|\widetilde{\mathcal L}_\omega^{n_0}u\|_{\mathrm{BV}}
\le A(Rq)^{n_0}\|u\|_{\mathrm{BV}}
+(RB)^{n_0}\|u\|_{L^1(\mathrm{Leb})},
\]
where $n_0$ is chosen so that $A(Rq)^{n_0}<1$.
The compact embedding
$\mathrm{BV}([0,1])\hookrightarrow L^1(\mathrm{Leb})$
then verifies \eqref{LL3}.

Finally, we verify \eqref{LL4}.
Lemma~2.5.9, conditions (C5') and (C7'), and Appendix~2.A of \cite{AFGTV-TFPF} give 
\begin{equation}\label{eq:20260928}
\esssup_\omega\|\phi_{\omega,c}\|_{\mathrm{BV}}<\infty \qquad\text{and}\qquad
 \essinf_\omega\inf_{x\in[0,1]}\phi_{\omega,c}(x)>0,
 \end{equation}
and that there are $C\ge1$ and $\kappa_c\in(0,1)$ such that
\[
 \bigl\|\tilde{\mathcal L}_{\omega,c}^k u
       -\nu_{\omega,c}(u)\phi_{\sigma^k\omega,c}\bigr\|_{\mathrm{BV}}
 \le C\kappa_c^k\|u\|_{\mathrm{BV}},\qquad k\ge1.
\]
Here $\|\cdot\|_{\mathrm{BV}}$ denotes the $\mathrm{BV}_1$ norm of \cite[Section~2.5]{AFGTV-TFPF}, defined using Lebesgue measure; see \cite[Proposition~2.5.2 and (2.5.13)]{AFGTV-TFPF} for the norm equivalence used in this spectral-gap estimate.
The weak norm used below is $\|\cdot\|_{L^1(\mathrm{Leb})}$.
We emphasize that no absolute continuity of $\nu_{\omega,c}$ with respect to Lebesgue measure is assumed.
By (E1)--(E3) of \cite[Section~2.5]{AFGTV-TFPF},
we have
$\esssup_\omega\|g_\omega T'_\omega\|_\infty<\infty$
and $\essinf_\omega\inf g_\omega>0$.
Since $T_\omega$ is surjective and
$\nu_{\sigma\omega,c}$ is a probability measure,
we also have $\lambda_{\omega,c}\ge\inf g_\omega$.
Consequently,
\begin{align*}
 \|\tilde{\mathcal L}_{\omega,c}u\|_{L^1(\mathrm{Leb})}
 \le\lambda_{\omega,c}^{-1}
       \int_0^1g_\omega|T'_\omega|\,|u|\,d\mathrm{Leb}
 \le\lambda_{\omega,c}^{-1}\|g_\omega T'_\omega\|_\infty
                  \|u\|_{L^1(\mathrm{Leb})}
 \le C\|u\|_{L^1(\mathrm{Leb})}.
\end{align*}
Therefore, applying the above spectral gap to $|u|$, and using $\||u|\|_{\mathrm{BV}}\le\|u\|_{\mathrm{BV}}$, yields
\begin{align*}
 \nu_{\omega,c}(|u|)\|\phi_{\sigma^k\omega,c}\|_{L^1(\mathrm{Leb})}
 \le\|\tilde{\mathcal L}_{\omega,c}^k|u|\|_{L^1(\mathrm{Leb})}
            +C\kappa_c^k\|u\|_{\mathrm{BV}}
 \le C^k\|u\|_{L^1(\mathrm{Leb})}
            +C\kappa_c^k\|u\|_{\mathrm{BV}}.
\end{align*}
The term $\|\phi_{\omega,c}\|_{L^1(\mathrm{Leb})}$ has a positive lower bound independent of $\omega,k$ and of the holes, by \eqref{eq:20260928}.
Hence, increasing $C$ to absorb its reciprocal, we obtain, for every $k\ge1$,
\begin{equation}\label{eq:closed-weak-mass-comparison}
 \nu_{\omega,c}(|u|)
 \le C^{k+1}\|u\|_{L^1(\mathrm{Leb})}
       +C\kappa_c^k\|u\|_{\mathrm{BV}}.
\end{equation}

The weak projection convergence in
\cite[Theorem~4.8]{crimmins_stability_2019} is precisely
\[
 \esssup_\omega\sup_{\|u\|_{\mathrm{BV}}\le1}
 \bigl\|\nu_{\omega,\infty}(u)\phi_\omega
       -\nu_{\omega,c}(u)\phi_{\omega,c}\bigr\|_{L^1(\mathrm{Leb})}
 \longrightarrow0
\]
as the holes shrink uniformly within the stated family.
Evaluating at $u=\widehat\Xi$ gives
\begin{equation*}
 \esssup_\omega\|\phi_\omega-\phi_{\omega,c}\|_{L^1(\mathrm{Leb})}
 \longrightarrow0.
\end{equation*}
On the other hand, by the uniform strong-norm bound on the projections in the same theorem,
\begin{equation*}
 \esssup_\omega\|\phi_\omega\|_{\mathrm{BV}}
                    \le C.
\end{equation*}
Substitute $u=\phi_\omega-\phi_{\omega,c}$ into
\eqref{eq:closed-weak-mass-comparison}. After increasing $C$ once more, together with \eqref{eq:20260928},
\begin{align*}
 \esssup_\omega|\nu_{\omega,c}(\phi_\omega)-1|
 &\le C^{k+1}\esssup_\omega
                  \|\phi_\omega-\phi_{\omega,c}\|_{L^1(\mathrm{Leb})}
       +C\kappa_c^k.
\end{align*}
First fix $k$ so that $C\kappa_c^k<1/4$, and then make the holes small
enough that the first term is less than $1/4$. Thus
\begin{equation*}
 \essinf_\omega\nu_{\omega,c}(\phi_\omega)
 \ge1-\esssup_\omega|\nu_{\omega,c}(\phi_\omega)-1|
 \ge\frac12.
\end{equation*}
Together with the uniform bound on the number of hole components, this yields \eqref{LL4}.

Theorems \ref{main thm qstat and qerg measures}--\ref{be} therefore apply to this class.
More specific examples may be found in Section 2.7 \cite{AFGTV-TFPF}. 

\begin{remark}
As an immediate application of Theorem~\ref{thm:e1}(3), one can define
trajectory-wise conditional Lyapunov exponents.
Consider the deterministic case in the one-dimensional
setting of this subsection.
Assume in addition that $\log|T'|\in\mathrm{BV}([0,1])$ and that $T'$ is defined and nonzero along the survivor trajectories under consideration. Applying Theorem~\ref{thm:e1}(3) to
$
  f(x)=\log |T'(x)|
$
yields that the conditional Lyapunov exponent
\[
  \lim_{n\to\infty}\frac1n \log |(T^n)'(x_n)|
\]
exists and coincides with $\eta_\infty(\log |T'|)$ for $(\prod_{n=0}^\infty \eta_{n})$-a.e.~$(x_n)_{n\ge 0} \in \prod_{n=0}^\infty X _{n}$.
A similar but distinct notion of conditional Lyapunov exponent was studied in \cite{ELR2019} for Brownian-driven stochastic differential equations, using conditional expectations given survival rather than almost-sure convergence.
\end{remark}

\subsection{Multidimensional expanding maps with small random holes}
\label{s:multidimensional-expanding}
 All of our main theorems can be applied to multidimensional expanding maps. 
 Below we describe how to satisfy our hypotheses for a fixed multidimensional expanding map with small random holes. 
 We emphasise that the random holes may move about in the phase space arbitrarily from one $\omega$-fibre to the next.
 First we note that \eqref{M1} and \eqref{M2} are trivial.
 We will start by defining a class of expanding maps and a corresponding Banach space. 
Let $X$ be a compact, connected $d$-dimensional submanifold of $\mathbb R^d$ with boundary, and let $T:X\to X$ be a surjective map.
After a uniform rescaling of coordinates, we assume that $\mathrm{Leb}(X)=1$.
We assume that $\mathcal U=\{U_{i,1}\}_i$ is a finite partition of $X$ modulo Lebesgue-null sets into open sets, and that each restriction $T|_{U_{i,1}}$ has a $C^1$ diffeomorphic extension $T_i$ to an open neighbourhood $W_{i,1}$ of $\overline{U_{i,1}}$.
Write $T_i^{-1}:T_i(W_{i,1})\to W_{i,1}$ for its inverse.
Throughout this subsection, $\mathcal L_c$ is the transfer operator with weight $g=|\det DT|^{-1}$.

We further assume
\begin{itemize}
    \item[\mylabel{P1}{P1}] the Jacobian of $T$ is uniformly $\alpha$-H\"older for some $0<\alpha\le 1$, namely there are $c_1<\infty$ and $\eps_0>0$ depending on the map $T$ (these quantities will enter also the definition of the norm below), we have
\begin{equation*}
|\det D_xT_i^{-1}-\det D_yT_i^{-1}|\le c_1|\det D_zT_i^{-1}|\eps^{\alpha},
    \end{equation*}
for each $i$, every $z\in T_i(W_{i,1})$, every $x,y\in B_\varepsilon(z)\cap T_i(W_{i,1})$, and every $0<\varepsilon\le\varepsilon_0$.
\item[\mylabel{P2}{P2}]  $s(T):=\max_i\sup_{x\in T_i(W_{i,1})}\|D_xT_i^{-1}\|<1.$
\item[\mylabel{P3}{P3}]
There are compact $C^1$ codimension-one hypersurface pieces $\Gamma_1,\ldots,\Gamma_S$ such that $\bigcup_i\partial U_{i,1}\subset\bigcup_{\ell=1}^S\Gamma_\ell$.
\end{itemize}

The choices of $n_0$, $\varepsilon_0$, and the maximum hole radius below may depend on $\alpha$ in \eqref{P1}.
We now  introduce  the  quasi-H\"older spaces, 
whose importance for expanding dynamical systems was stressed in the seminal works by Keller \cite{KLG}  and Saussol \cite{saussol00}.
For $h\in L^1(X,\mathrm{Leb})$, we use its zero extension to $\mathbb R^d$, still denoted by $h$.
All oscillations in this subsection are computed for this extension.
For $0<\alpha\le1$ and $\varepsilon_0>0$, define
\[
|h|_\alpha
:=
\sup_{0<\varepsilon\le\varepsilon_0}
\varepsilon^{-\alpha}
\int_{\mathbb R^d}
\operatorname{osc}(h,B_\varepsilon(x))\,d\mathrm{Leb}(x),
\]
where $\text{osc}(h,A):=\esssup_{x,y\in A}|h(x)-h(y)|$ for any measurable set $A\subset\mathbb R^d$.
We set
\[
V_\alpha(X)
:=
\{h\in L^1(X,\mathrm{Leb}):|h|_\alpha<\infty\},
\qquad
\|h\|_\alpha:=|h|_\alpha+\|h\|_1,
\]
where $\|h\|_1=\int_X|h|\,d\mathrm{Leb}$.
Although the seminorm depends on $\varepsilon_0$, the space $V_\alpha(X)$ does not.
We take $\mathcal B=(V_\alpha(X),\|\cdot\|_\alpha)$. 
One can show  that $\mathcal{B}$ can be continuously injected into $L^{\infty}$, and in particular $||h||_{\infty}\le C_s ||h||_{\a}$, where the constant $C_s$ is  explicitly determined and depends in particular on  $\eps_0,$ see \cite{saussol00}, equation (6).
This demonstrates the second part of \eqref{B2}.
For the first part of \eqref{B2}, see Proposition 3.4(iii)\footnote{Note that the `inf' should be a `sup' in the expression in Proposition 3.4(iii) \cite{saussol00}.} in \cite{saussol00}.  
The single operator $\mathcal{L}_c$ is bounded in $\mathcal{B}$ so we obtain \eqref{C0}.
The closed conformal measure is $\nu_{\omega,c}=\mathrm{Leb}|_X$, with $\lambda_{\omega,c}=1$, so \eqref{C1} holds.
By \eqref{P3}, restrictions to $X$ of smooth functions defined on a neighbourhood of $X$ belong to $V_\alpha(X)$.
These restrictions are dense in both $C(X)$ and $L^1(X,\mathrm{Leb})$, which proves \eqref{B1}, see Proposition 1.8 in \cite{KLG}.

Following Saussol \cite{saussol00}, we define the {\em multiplicity entropy} $H_{\text{mult}}(\mathcal U,T)$ as
\begin{equation}\label{mpo}
H_{\text{mult}}(\mathcal{U},T)=\lim_{n\rightarrow \infty}\frac1n \log \sup_x \#\{U_{i,n}\in \mathcal{U}^{(n)}; x\in \overline{U_{i,n}}\},
\end{equation}
where $\mathcal{U}^{(n)}:=\mathcal{U}\wedge T^{-1}\mathcal{U}\dots \wedge T^{-(n-1)}\mathcal{U}$, 
and the {\em dilation coefficient} $\delta(T)$ as  
\begin{equation}\label{mqo}
\delta(T)=\lim_{n\rightarrow \infty}\frac1n \log s(T^n)
\end{equation}
(recall \eqref{P2} for the definition of $s(T^n)$).
The two limits exist by Fekete's lemma, since the corresponding multiplicities and the quantities \(s(T^n)\) are submultiplicative.
We now state our main assumption on the dynamics $T$:
\begin{itemize}
\item[\mylabel{P4}{P4}]
\begin{equation*}
%\label{mult}
H_{\text{mult}}(\mathcal U,T)+\delta(T)<0.
\end{equation*}
\end{itemize}
This controls the growth of the singularities of the map relative to its expanding properties.
Conditions \eqref{P1}--\eqref{P4} give quasi-compactness of $\mathcal L_c$ on $V_\alpha(X)$ with spectral radius one, see \cite[Lemma~2.2 and Theorem~5.1]{saussol00}.

We further assume
\begin{itemize}
\item[\mylabel{P5}{P5}]
The spectrum of $\mathcal L_c:V_\alpha(X)\to V_\alpha(X)$ on the unit circle consists only of the simple eigenvalue $1$.
\end{itemize}
Under \eqref{P1}--\eqref{P4}, a sufficient condition of \eqref{P5} is total transitivity (i.e.~for every $k\ge1$ and all nonempty open sets $U,V\subset X$, there exists $n\ge1$ such that $\Leb(U\cap T^{-kn}V)>0$).
Indeed, the finite cyclic decomposition in \cite[Lemma~2.2 and Theorem~5.1]{saussol00}, together with the local positivity of the component densities in \cite[Lemma~3.1]{saussol00}, then reduces to a single component of period one.

\medskip

We now introduce random holes $H_\omega\subset X$, $\omega\in\Omega$,
satisfying the following conditions.
\begin{itemize}
\item[\mylabel{P6}{P6}]
For $m$-a.e.\ $\omega$, there are $c_\omega\in X$ and
$r_\omega>0$ such that $H_\omega\subset B_{r_\omega}(c_\omega)$ and $r:=\esssup_{\oio}r_\omega <\infty$.
Moreover, there are constants $K_1\in\mathbb N$, $\rho>0$, and $K_2<\infty$, independent of $\omega$, such that $r_\omega^{-1}(\partial H_\omega-c_\omega)$ can be covered by at most $K_1$ pieces of Lipschitz graphs.
After a suitable translation and rotation, each piece lies in the graph of a $K_2$-Lipschitz function that is defined on a Euclidean ball of radius $\rho$ in $\mathbb R^{d-1}$ and can be extended to the concentric ball of radius $2\rho$ with the same Lipschitz constant.
\item[\mylabel{P7}{P7}]
The space $\Omega$ is Polish, $\sg:\Om\to\Om$ is a homeomorphism,
and $\omega\mapsto H_\omega$ has countable range, with
$\{\omega:H_\omega=H_j\}$ measurable for each $H_j$ in that range.
\end{itemize}

The Lipschitz graph condition in \eqref{P6} holds uniformly for ball holes $H_\omega=B_{r_\omega}(c_\omega)\subset X$, choosing $r_\omega$ to be the actual radius.
It also holds for any deterministic hole whose boundary admits a finite Lipschitz graph cover.

For condition \eqref{B3}, we verify $|\hat{X}_{\om,n}|_\alpha<\infty$ for $m$-a.e.\ $\oio$ and each $n\ge 0$.
Since $T$ has finitely many branches, this follows from Proposition~3.2(ii) in \cite{saussol00}, applied with $f=1$, using the finite Lipschitz boundary cover in \eqref{P6} and the $C^1$ inverse branches.

It remains to verify \eqref{LL}.
Following \cite[Lemma~2.5.10 and Appendix~2.B]{AFGTV-TFPF}, we apply the $m$-continuous version of
\cite[Theorem~4.8]{crimmins_stability_2019}.
The strong and weak norms are $\|\cdot\|_\alpha$ and $\|\cdot\|_1$, respectively.
We verify the following three conditions:
\begin{enumerate}[label=(\roman*),leftmargin=*,itemsep=2pt]
\item The closed transfer operator cocycle has a hyperbolic splitting with a one-dimensional leading subspace and an exponentially contracting complementary subspace.
\item The open cocycles satisfy a uniform Lasota--Yorke inequality:
there exist $r_0>0$, $0\le A,B<\infty$, and $0<\gamma<1$ such that, for every family of holes satisfying \eqref{P6} with $0<r\le r_0$, every $f\in V_\alpha$ and $m$-a.e.\ $\omega$,
\begin{equation}\label{saussolLY}
\|\mathcal L^k_\omega f\|_\alpha \le A\gamma^k\|f\|_\alpha+B\|f\|_1, \qquad k\ge0.
\end{equation}
\item The open operators converge to the closed operator in the strong-to-weak operator norm:
\[
\lim_{r\to0}\operatorname*{ess\,sup}_{\omega\in\Omega} \|\mathcal L_c-\mathcal L_\omega\|_{V_\alpha\to L^1}=0.
\]
\end{enumerate}
For (i), \cite[Theorem~5.1]{saussol00} gives quasi-compactness of $\mathcal L_c$ under \eqref{P1}--\eqref{P4}, and \eqref{P5} ensures that $1$ is its only peripheral eigenvalue and is simple.
Thus its leading eigenspace is spanned by the closed invariant density, and its complementary spectral subspace is exponentially contracting. Since $T$ is fixed, this gives the required splitting of the constant closed cocycle.

For (iii), let $\gamma_d$ be the volume of the unit ball in $\mathbb R^d$.
Since $\mathcal L_\omega f=\mathcal L_c(f\widehat{X_{\omega,0}})$, the $L^1$-contraction of $\mathcal L_c$, \eqref{P6}, and $\|f\|_\infty\le C_s\|f\|_\alpha$ give
\[
\begin{aligned}
\|(\mathcal L_c-\mathcal L_\omega)f\|_1
=\|\mathcal L_c(f\widehat{H_\omega})\|_1
\le\|f\|_\infty\,\mathrm{Leb}(H_\omega)
\le C_s\gamma_d r^d\|f\|_\alpha.
\end{aligned}
\]
For fixed $\varepsilon_0$, the constant $C_s$ is independent of
$\omega$ and $r$, so (iii) follows by taking the operator norm
and then letting $r\to0$.

It remains to prove (ii).
Recall that $\mathcal L_\omega^n f=\mathcal L_c^n(f\widehat{X_{\omega,n-1}})$.
For each $U_{i,n}\in\mathcal U^{(n)}$, let $T_i^n$ be the $C^1$ diffeomorphic extension of $T^n|_{U_{i,n}}$ obtained by composing the corresponding extensions of $T$, on an open neighbourhood $W_{i,n}$ of $\overline{U_{i,n}}$.
Write $T_i^{-n}$ for its inverse on $T_i^n(W_{i,n})$.
On this branch put $g_n=|\det DT_i^n|^{-1}$, and set
$V_{i,n,\omega}:=U_{i,n}\cap X_{\omega,n-1}$.
Then
\begin{equation}\label{roi}
\begin{aligned}
\operatorname{osc}(\mathcal L^n_\omega f,B_\varepsilon(x))
&\le\sum_i\operatorname{osc}\Bigl(
\bigl[(fg_n\widehat{X_{\omega,n-1}})\circ T_i^{-n}\bigr]
\widehat{T^nU_{i,n}},B_\varepsilon(x)\Bigr)\\
&=\sum_i\operatorname{osc}\Bigl(
\bigl[(fg_n)\circ T_i^{-n}\bigr]
\widehat{T^nV_{i,n,\omega}},B_\varepsilon(x)\Bigr).
\end{aligned}
\end{equation}
The boundary of $V_{i,n,\omega}$ is contained in $\partial U_{i,n}$ and the branchwise pullbacks of $\partial H_{\sigma^j\omega}$, $0\le j<n$.
All these boundaries are Lebesgue null; therefore replacing $V_{i,n,\omega}$ by its interior does not change the essential oscillations in \eqref{roi}.

To estimate \eqref{roi}, we use the quotient from \cite[condition (PE5)]{saussol00}:
\[
\begin{aligned}
G_{n,\omega}(\varepsilon,\varepsilon_0)
:=\sup_{x\in\mathbb R^d}\sum_i
\frac{\mathrm{Leb}\bigl(
(T_i^n)^{-1}B_\varepsilon(\partial(T_i^nV_{i,n,\omega}))
\cap B_{(1-s(T^n))\varepsilon_0}(x)\bigr)}
{\gamma_d((1-s(T^n))\varepsilon_0)^d},
\end{aligned}
\]
where $B_t(A):=\{x:\operatorname{dist}(x,A)<t\}$ for a set $A$.
The inverse branches are understood on their smooth extensions.
Let $c_{n_0}$ be the H\"older distortion constant in \eqref{P1} for $T^{n_0}$, and set
\[
\eta:=(1+c_{n_0}s(T^{n_0})^\alpha\varepsilon_0^\alpha)
\Bigl[s(T^{n_0})^\alpha+2\operatorname*{ess\,sup}_{\omega}
\sup_{0<\varepsilon\le\varepsilon_0}
\Bigl(\frac{\varepsilon_0}{\varepsilon}\Bigr)^\alpha
G_{n_0,\omega}(\varepsilon,\varepsilon_0)\Bigr].
\]
Notice that $\eta$ depends on $n_0,\varepsilon_0$ and $r$ (recall the definition of $V_{i,n,\omega}$ and \eqref{P6}).
When $\eta<\infty$, applying the calculation in
\cite[proof of Lemma~4.1]{saussol00} to \eqref{roi} gives\footnote{Combining
(9) and (11) in that proof gives, in Saussol's notation,
\[
\begin{aligned}
\int\operatorname{osc}(Pf,B_\varepsilon(x))\,dx
&\le (1+c(s\varepsilon)^\alpha)
\bigl[s^\alpha\varepsilon^\alpha
+2G(\varepsilon,\varepsilon_0)\varepsilon_0^\alpha\bigr]|f|_\alpha
+C\varepsilon^\alpha\|f\|_1.
\end{aligned}
\]
Here $C$ may depend on the fixed parameters.
Use $P=\mathcal L_\omega^{n_0}$, $T=T^{n_0}$, $U_i=V_{i,n_0,\omega}$, $s=s(T^{n_0})$, $c=c_{n_0}$, and $G=G_{n_0,\omega}$.
Dividing by $\varepsilon^\alpha$ and taking the supremum over $0<\varepsilon\le\varepsilon_0$ gives the stated inequality.}
\[
|\mathcal L_\omega^{n_0}f|_\alpha
\le\eta|f|_\alpha+D\|f\|_1
\]
with $D$ independent of $\omega$ and $r$.
\begin{lemma}
For the fixed $\alpha$ and the constants in \eqref{P6}, we can choose $n_0\ge1$, $\varepsilon_0>0$, and $r_0>0$ so that $\eta<1$ for every family of holes satisfying \eqref{P6} with $0<r\le r_0$.
\end{lemma}
\begin{proof}
Fix $n\ge1$. For $0\le j<n$, set
$\Gamma^H_{i,j,n,\omega}:=(T^j|_{\overline{U_{i,n}}})^{-1}
(\partial H_{\sigma^j\omega})$, where the restriction is understood
through its continuous branch extension. By the definition of
$V_{i,n,\omega}$,
\[
\partial V_{i,n,\omega}\subset\partial U_{i,n}
\cup\bigcup_{j=0}^{n-1}\Gamma^H_{i,j,n,\omega}.
\]
For sufficiently small $\varepsilon_0$, \eqref{P2} therefore gives,
with $0<\varepsilon\le\varepsilon_0$,
\begin{equation}\label{eq:boundary-pullback}
\begin{aligned}
(T_i^n)^{-1}
B_\varepsilon\bigl(\partial(T_i^n V_{i,n,\omega})\bigr)
\subset
B_{s(T^n)\varepsilon}(\partial V_{i,n,\omega})
\subset
B_{s(T^n)\varepsilon}(\partial U_{i,n})
\cup
\bigcup_{j=0}^{n-1}
B_{s(T^n)\varepsilon}(\Gamma^H_{i,j,n,\omega}).
\end{aligned}
\end{equation}
According to this inclusion, we have 
\begin{align}\label{eq:20261108}
G_{n,\omega}(\varepsilon,\varepsilon_0)\le G_{n,\omega}^{(1)}(\varepsilon,\varepsilon_0)+G_{n,\omega}^{(2)}(\varepsilon,\varepsilon_0),
\end{align}
where
\begin{align*}
&G_{n,\omega}^{(1)}(\varepsilon,\varepsilon_0):=\sup_{x\in \mathbb R^d}\sum_i\frac{\mathrm{Leb}\bigl(B_{s(T^n)\varepsilon}(\partial U_{i,n})
\cap B_{(1-s(T^n))\varepsilon_0}(x)\bigr)}
{\gamma_d((1-s(T^n))\varepsilon_0)^d},\\
&G_{n,\omega}^{(2)}(\varepsilon,\varepsilon_0):=\sup_{x\in \mathbb R^d}\sum_i\sum_{j=0}^{n-1}
\frac{\mathrm{Leb}\bigl(B_{s(T^n)\varepsilon}(\Gamma^H_{i,j,n,\omega})
\cap B_{(1-s(T^n))\varepsilon_0}(x)\bigr)}
{\gamma_d((1-s(T^n))\varepsilon_0)^d}.
\end{align*}

The estimate for $G_{n,\omega}^{(1)}$, which does not involve the holes, is obtained directly from \cite{saussol00} as follows.
Recall that $S$ is the number of hypersurface pieces in \eqref{P3},
and write $M_n:=\sup_x\#\{i:x\in\overline{U_{i,n}}\}$, as in \eqref{mpo}.
For each fixed $n$, the estimates in the proofs of
\cite[Lemmas~2.1 and~2.2]{saussol00}, applied to $T^n$, give,
for sufficiently small $\varepsilon_0$,\footnote{The proof of Lemma~2.1 uses
$T_i^{-1}B_\varepsilon(\partial T U_i)
\subset B_{s\varepsilon}(\partial U_i)$
and estimates the latter neighbourhoods. Thus its final bound also
applies to $G_{n,\omega}^{(1)}$ for $T^n$.
In Saussol's notation, when
$s\varepsilon<(1-s)\varepsilon_0$, this gives
\[
G(\varepsilon,\varepsilon_0)
\le \frac{2\gamma_{N-1}}{\gamma_N}Y
\frac{s\varepsilon}{(1-s)\varepsilon_0}(1+o(1)).
\]
The proof of Lemma~2.2 gives $Y=Y(T^n)\le SnM_n$.
Take $N=d$ and $s=s(T^n)$, choose $\varepsilon_0$ so that
$1+o(1)\le2$, and use $(\varepsilon/\varepsilon_0)^{1-\alpha}\le1$.
For fixed $n$, the error is uniform in $x$, the finitely many
compact $C^1$ pieces, and $0<\varepsilon\le\varepsilon_0$.
When $s\varepsilon\ge(1-s)\varepsilon_0$, the same bound follows
because $\Leb (B_{s\varepsilon}(\Gamma)\cap B_{(1-s)\varepsilon_0}(x))
\le \Leb (B_{(1-s)\varepsilon_0}(x))
=\gamma_N\bigl((1-s)\varepsilon_0\bigr)^N
\le 2\gamma_{N-1}s\varepsilon \bigl((1-s)\varepsilon_0\bigr)^{N-1}$
%with $\Gamma =\partial U_{i,n}$.
for each compact $C^1$ hypersurface piece $\Gamma$
covering $\partial U_{i,n}$.
For $d=1$, each such piece consists of finitely many points.
For fixed $n$, we may choose $\varepsilon_0$ so small that
every $B_{\varepsilon_0}(x)$ meets at most one point of
each piece. Since
$s\varepsilon+(1-s)\varepsilon_0\le\varepsilon_0$, we have
$
\Leb(B_{s\varepsilon}(\Gamma)
\cap B_{(1-s)\varepsilon_0}(x))
\le 2s\varepsilon.
$}
\begin{equation}\label{eq:partition-boundary-quotient}
\begin{aligned}
\sup_{0<\varepsilon\le\varepsilon_0}
\left(\frac{\varepsilon_0}{\varepsilon}\right)^\alpha
G_{n,\omega}^{(1)}(\varepsilon,\varepsilon_0)
\le
\frac{4\gamma_{d-1}}{\gamma_d}\,SnM_n
\frac{s(T^n)}{1-s(T^n)}.
\end{aligned}
\end{equation}

We estimate $G_{n,\omega}^{(2)}(\varepsilon,\varepsilon_0)$.
First, \eqref{P6} gives
\begin{equation}\label{eq:hole-tube-volume}
\mathrm{Leb}(B_t(\partial H_\omega))
\le C\,t(r_\omega+t)^{d-1}\qquad(t>0).
\end{equation}
Indeed, when $0<t\le\rho r_\omega$, integration in the vertical direction
in the extended graph coordinates of \eqref{P6} gives
$C t r_\omega^{d-1}$. When $t>\rho r_\omega$, the inclusion
$B_t(\partial H_\omega)\subset B_{r_\omega+t}(c_\omega)$ gives
\eqref{eq:hole-tube-volume} directly.
For fixed $n$, \eqref{eq:hole-tube-volume} and change of variables
in the finitely many $C^1$ branch extensions give
$\mathrm{Leb}(B_t(\Gamma^H_{i,j,n,\omega}))\le C_n t(r+t)^{d-1}$
for sufficiently small $r,t$. Also, $\Gamma^H_{i,j,n,\omega}$ is
contained in a ball of radius $C_n r$. Hence
\begin{equation}\label{eq:hole-pullback-volume}
\mathrm{Leb}(B_t(\Gamma^H_{i,j,n,\omega}))
\le \min\bigl\{C_n t(r+t)^{d-1},\,
\gamma_d(C_n r+t)^d\bigr\}.
\end{equation}
Constants denoted by $C$ and $C_n$ may increase between estimates;
$C$ is independent of $n,\omega,r,\varepsilon,\varepsilon_0$, whereas
$C_n$ may depend on $n$ but not on $\omega$ or these radii.

When $d\ge2$ and $t=s(T^n)\varepsilon\le C_n r$,
the first bound in \eqref{eq:hole-pullback-volume} gives
\[
\begin{aligned}
\left(\frac{\varepsilon_0}{\varepsilon}\right)^\alpha
\frac{\mathrm{Leb}(B_t(\Gamma^H_{i,j,n,\omega}))}
{\gamma_d((1-s(T^n))\varepsilon_0)^d}
\le C_n\frac{s(T^n)\varepsilon^{1-\alpha}r^{d-1}}
{(1-s(T^n))^d\varepsilon_0^{d-\alpha}}
\le C_n\left(\frac r{\varepsilon_0}\right)^{d-\alpha}.
\end{aligned}
\]
The last inequality uses $\varepsilon\le C_n r/s(T^n)$ and
$1-\alpha\ge0$, with the factors depending on $n$ absorbed into $C_n$.
When $d\ge2$ and $t>C_n r$, the second bound in
\eqref{eq:hole-pullback-volume} gives
$\mathrm{Leb}(B_t(\Gamma^H_{i,j,n,\omega}))\le\gamma_d(2t)^d$,
and therefore
\[
\begin{aligned}
\left(\frac{\varepsilon_0}{\varepsilon}\right)^\alpha
\frac{\mathrm{Leb}(B_t(\Gamma^H_{i,j,n,\omega}))}
{\gamma_d((1-s(T^n))\varepsilon_0)^d}
\le 2^d\left(\frac{s(T^n)}{1-s(T^n)}\right)^d
\left(\frac{\varepsilon}{\varepsilon_0}\right)^{d-\alpha}
\le 2^d\left(\frac{s(T^n)}{1-s(T^n)}\right)^d.
\end{aligned}
\]
When $d=1$, $\#\partial H_{\sigma^j\omega}\le K_1$ by \eqref{P6},
so injectivity of the branch implies
$\mathrm{Leb}(B_t(\Gamma^H_{i,j,n,\omega}))\le2K_1t$.
Since $\gamma_1=2$, this gives
\[
\left(\frac{\varepsilon_0}{\varepsilon}\right)^\alpha
\frac{\mathrm{Leb}(B_t(\Gamma^H_{i,j,n,\omega}))}
{2(1-s(T^n))\varepsilon_0}
\le K_1\frac{s(T^n)}{1-s(T^n)}
\left(\frac{\varepsilon}{\varepsilon_0}\right)^{1-\alpha}
\le K_1\frac{s(T^n)}{1-s(T^n)}.
\]
For this fixed $n$, decrease $\varepsilon_0$ so that 
\begin{equation}\label{eq:20261108c}
\sup_{x\in\mathbb R^d}
\#\{i:B_{\varepsilon_0}(x)\cap\overline{U_{i,n}}
\ne\varnothing\}\le M_n.
\end{equation}
This is possible by finiteness and compactness: otherwise,
$M_n+1$ fixed closures would contain a common limit point
of centres of balls with radii tending to zero.
Since $\Gamma^H_{i,j,n,\omega}\subset\overline{U_{i,n}}$ and
$t+(1-s(T^n))\varepsilon_0\le\varepsilon_0$, a nonempty
intersection
$
B_t(\Gamma^H_{i,j,n,\omega})
\cap B_{(1-s(T^n))\varepsilon_0}(x)
$
implies that $\overline{U_{i,n}}$ meets
$B_{\varepsilon_0}(x)$.
Thus at most $M_n$ indices $i$ occur.
For each $i$ there are $n$ indices $j$, so summing the
preceding estimates gives
\begin{equation}\label{eq:hole-sum-bound}
\begin{aligned}
\sup_{0<\varepsilon\le\varepsilon_0}
\left(\frac{\varepsilon_0}{\varepsilon}\right)^\alpha
G_{n,\omega}^{(2)}(\varepsilon,\varepsilon_0)
\le nM_n
\begin{cases}
2^d\bigl(s(T^n)/(1-s(T^n))\bigr)^d
+C_n(r/\varepsilon_0)^{d-\alpha},&d\ge2,\\
K_1s(T^n)/(1-s(T^n)),&d=1.
\end{cases}
\end{aligned}
\end{equation}

Combining \eqref{eq:20261108}, \eqref{eq:partition-boundary-quotient} and \eqref{eq:hole-sum-bound}, absorbing the fixed $S$ into $C$,
and using
$(s(T^n)/(1-s(T^n)))^d\le s(T^n)/(1-s(T^n))$ when
$s(T^n)\le1/2$, we obtain
\begin{equation}\label{eq:direct-boundary-quotient}
\begin{aligned}
2\sup_{0<\varepsilon\le\varepsilon_0}
\left(\frac{\varepsilon_0}{\varepsilon}\right)^\alpha
G_{n,\omega}(\varepsilon,\varepsilon_0)
\le C nM_n\frac{s(T^n)}{1-s(T^n)}
+\begin{cases}
 C_n nM_n(r/\varepsilon_0)^{d-\alpha},&d\ge2,\\
 0,&d=1.
\end{cases}
\end{aligned}
\end{equation}
By \eqref{mpo}--\eqref{mqo}, the exponential growth rate of
$nM_ns(T^n)$ is $H_{\mathrm{mult}}(\mathcal U,T)+\delta(T)$,
which is negative by \eqref{P4}; also $s(T^n)\to0$ by \eqref{P2}.
Thus
\[
s(T^n)^\alpha+nM_ns(T^n)\longrightarrow0.
\]
Since $(1-s(T^n))^{-1}\le2$ for all sufficiently large $n$,
we can first choose $n_0$ so that $s(T^{n_0})\le1/2$ and the
first term on the right-hand side of
\eqref{eq:direct-boundary-quotient}, together with
$s(T^{n_0})^\alpha$, is less than $1/4$.
Then choose $\varepsilon_0>0$ sufficiently small that 
\eqref{eq:boundary-pullback},
\eqref{eq:partition-boundary-quotient} and \eqref{eq:20261108c} hold for both $n=1$ and $n=n_0$.
We also require that, for $n=1,n_0$,
\eqref{eq:hole-pullback-volume} holds for all
$0<t\le s(T^n)\varepsilon_0$ whenever $r$ is sufficiently
small, and that
\[
1+c_{n_0}s(T^{n_0})^\alpha\varepsilon_0^\alpha<2.
\]
Finally, choose $r_0>0$ so that \eqref{eq:hole-pullback-volume}
holds for $0<r\le r_0$ at both $n=1$ and $n=n_0$ and, when $d\ge2$,
\[
C_{n_0}n_0M_{n_0}(r_0/\varepsilon_0)^{d-\alpha}<1/4.
\]
This last choice is possible because $d-\alpha>0$ for $d\ge2$.
Consequently, $\eta<1$, with an upper bound strictly below $1$
independent of the family of holes with $0<r\le r_0$.
\end{proof}

Applying the same argument with $n=1$ gives
$\operatorname*{ess\,sup}_\omega
\|\mathcal L_\omega\|_{V_\alpha\to V_\alpha}<\infty$,
without requiring contraction. Since the open operators contract
$L^1(\mathrm{Leb})$,
iteration in blocks of $n_0$ steps, with at most $n_0-1$ steps left
over, proves \eqref{saussolLY} and completes (ii).

Now (i)--(iii) allow us to apply \cite[Theorem~3.9]{crimmins_stability_2019}.
For sufficiently small $r$, this gives measurable equivariant rank-one projections $P_\omega$ with uniformly contracting complements.
The projections are uniformly bounded and converge in the strong-to-weak operator norm to the closed-system projection $P_cu=\mathrm{Leb}(u)\phi_c$ as $r\to0$, where $\phi_c$ is the closed invariant density normalized by $\mathrm{Leb}(\phi_c)=1$.
In particular, $P_\omega\widehat\Xi\neq0$.
Set $\phi_\omega:=P_\omega\widehat\Xi$.
Since $P_\omega$ has rank one, there exists a unique bounded linear functional $\nu_{\omega,\infty}\in V_\alpha^*$ such that
$P_\omega u=\nu_{\omega,\infty}(u)\phi_\omega$ for $u\in V_\alpha$.
Since $P_\omega^2=P_\omega$, we have
$\nu_{\omega,\infty}(\widehat\Xi) =\nu_{\omega,\infty}(\phi_\omega)=1$.
The functionals $\nu_{\omega,\infty}$ are uniformly bounded, since the projections are uniformly bounded and $\|\phi_\omega\|_\alpha$ is uniformly bounded away from zero.
Equivariance and the reality of the splitting give real multipliers $\lambda_\omega$ such that
$\mathcal L_\omega\phi_\omega=\lambda_\omega\phi_{\sigma\omega}$ and $\nu_{\sigma\omega,\infty}(\mathcal L_\omega u)=\lambda_\omega\nu_{\omega,\infty}(u).$

Applying the projection convergence to $u=\widehat\Xi$ gives
\begin{equation}\label{eq:20261108b}
\operatorname*{ess\,sup}_\omega\|\phi_\omega-\phi_c\|_{L^1(\mathrm{Leb})}\longrightarrow0 \qquad(r\to0).
\end{equation}
Integrating $\mathcal L_\omega\phi_\omega=\lambda_\omega\phi_{\sigma\omega}$ gives
\[
\lambda_\omega=\frac{\mathrm{Leb}(\phi_\omega)-\int_{H_\omega}\phi_\omega\,d\mathrm{Leb}}{\mathrm{Leb}(\phi_{\sigma\omega})}.
\]
Since
\[
\left|\int_{H_\omega}\phi_\omega\,d\mathrm{Leb}\right|\le\|\phi_\omega-\phi_c\|_1+\|\phi_c\|_\infty\gamma_d r^d,
\]
we obtain
$\operatorname*{ess\,sup}_\omega|\lambda_\omega-1|\to0$ as $r\to0$.
Let $\mu\in(0,1)$ be a uniform contraction rate for the complementary bundle.
For sufficiently small $r$, we have
\[
\lambda_-:=\operatorname*{ess\,inf}_\omega\lambda_\omega>\max\{\mu,\gamma\},
\]
and $\log\lambda\in L^\infty(m)\subset L^1(m)$.
The uniform bounds on the complementary cocycle and the projections therefore give
\[
\left\|\widetilde{\mathcal L}_\omega^n u -\nu_{\omega,\infty}(u)\phi_{\sigma^n\omega} \right\|_\alpha \le D(\mu/\lambda_-)^n\|u\|_\alpha
\]
for some $D<\infty$ independent of $\omega$.
This proves \eqref{G3}.

To verify positivity, note that $\lambda_\omega>0$ and that the open operators preserve nonnegative functions.
Applying the preceding estimate along the backward orbit to $u=\widehat\Xi$ shows that $\phi_\omega\ge0$ almost everywhere.
Applying it forward to any nonnegative $u\in V_\alpha$ and integrating gives $\nu_{\omega,\infty}(u)\ge0$, since $\mathrm{Leb}(\phi_\omega)$ is uniformly bounded away from zero by \eqref{eq:20261108b}.
Thus \eqref{G1} holds, while \eqref{G2} follows from \eqref{eq:20261108b} and $\nu_{\omega,\infty}(\phi_\omega)=1$.
Moreover, the $k=1$ case of \eqref{saussolLY} gives a uniform bound in \eqref{G0}.
Together with the above estimates, this proves \eqref{LL1}.

Condition \eqref{LL2} follows directly from \eqref{P7}, since $T$ is fixed.

Since $\lambda_->\gamma$, \eqref{saussolLY} gives
\[
\|\widetilde{\mathcal L}_\omega^k f\|_\alpha
\le A(\gamma/\lambda_-)^k\|f\|_\alpha
+B\lambda_-^{-k}\|f\|_1.
\]
Choosing $k$ sufficiently large so that
$A(\gamma/\lambda_-)^k<1$ and using the compact
embedding \cite[Proposition~3.3]{saussol00}
verifies \eqref{LL3}.

Finally, it follows from \eqref{eq:20261108b} that, for all sufficiently small $r$,
\[
\operatorname*{ess\,inf}_\omega
\nu_{\omega,c}(\phi_\omega)
=\operatorname*{ess\,inf}_\omega
\mathrm{Leb}(\phi_\omega)\ge\frac12,
\]
while \eqref{eq:hole-tube-volume} yields
\[
\operatorname*{ess\,sup}_\omega
\|\hat X_{\omega,0}\|_\alpha
\le \mathrm{Leb}(X)+C\epsilon_0^{1-\alpha}
<\infty.
\]
Thus \eqref{LL4} holds.

 \subsection{Deterministic subshifts of finite type with small holes}\label{s:6.4}

Consider a finite set $E$ with at least two elements, called the alphabet.
Let $A:E\times E\to\{0,1\}$ be a matrix, which we call the incidence matrix. For each $n\geq 0$ denote the set of finite words of length $n$ by 
\begin{align*}
    E_A^n:=\{\gm\in E^n: A_{\gm_j\gm_{j+1}}=1 \quad  \forall \, 0\leq j< n-1\}
\end{align*}
and the set of (one-sided) infinite-length words
\begin{align*}
    E_A^\infty:=\{\gm\in E^\NN: A_{\gm_j\gm_{j+1}}=1 \quad \forall \, j\geq 0\}.
\end{align*}
Naturally, we let $E_A^*:=\cup_{n=0}^\infty E_A^n$ denote the set of all words of finite length.
Given $\gm\in E_A^*$, we let $[\gm]$ denote the cylinder set generated by $\gm$. 
Fix $\al\in(0,1)$, and denote the usual metric $d_\al$ on $E_A^\infty$ given by $d_\al(\gm,\rho):=\al^{|\gm\wedge\rho|}$, where $|\gm|$ denotes the length of $\gm$ and $\gm\wedge\rho$ denotes the initial segment common to $\gm$ and $\rho$. 

We consider the dynamical system $\vta:E_A^\infty\to E_A^\infty$, where $\vta$ denotes the usual left-shift map, i.e. $\vta(\gm)_n:=\gm_{n+1}$ for each $n\geq 0$.  
We assume that every symbol of $E$ occurs in $E_A^\infty$ and that $\vartheta$ is topologically mixing. Since $E$ is finite, $A$ is then primitive (i.e.~there exists an integer $n \ge 1$ such that every entry of \(A^{n}\) is strictly positive).
Let $\varphi:E_A^\infty\to\mathbb R$ be a H\"older continuous potential.
Increasing $\alpha$ if necessary, we ensure that $\varphi$ is Lipschitz continuous with respect to $d_\alpha$. The potential $\varphi$ is summable, meaning that
$\sum_{e\in E}\exp(\sup_{[e]}\varphi)<\infty$, because $E$ is finite and $\varphi$ is bounded on the compact space $E_A^\infty$.
We define the transfer operator (associated to $\vp$) acting (continuously) on the space $H_\al^b(E_A^\infty)$ of bounded, complex-valued Lipschitz continuous (with respect to the metric $d_\al$) functions by 
\begin{align*}
    \cL_\vp(g)(\gm):=\sum_{e\in E: \, A_{e\gm_0}=1}g(e\gm)\exp(\vp(e\gm)).
\end{align*}
Theorem 2.1.1 of \cite{pollicott_open_2017} gives that there is a unique shift-invariant Gibbs/equilibrium measure $\mu_\vp$. 
To use $\mu_\varphi$ as the reference measure, we normalize
$\varphi$ by adding a constant and a coboundary, without changing
its equilibrium measure, so that
\[
\mathcal L_\varphi\widehat{E_A^\infty}=\widehat{E_A^\infty},
\qquad
\mathcal L_\varphi^*\mu_\varphi=\mu_\varphi.
\]
We take $\Xi=E_A^\infty$, $\nu_c=\mu_\varphi$,
$\mathcal L_c=\mathcal L_\varphi$, and $\lambda_c=1$.
Given a function $g\in L^1(\mu_\vp)$ and an integer $m\geq0$, we define the oscillation function $\osc_m(g):E_A^\infty\to[0,\infty]$ by 
\begin{align*}
    \osc_m(g)(\gm):=\esssup\{|g(\rho)-g(\vrho)|:\rho,\vrho\in[\gm\rvert_m]\}.
\end{align*}
Here $\gamma|_m$ consists of the first $m$ symbols, $[\gamma|_0]=E_A^\infty$, and the essential supremum is taken with respect to $\mu_\varphi\otimes\mu_\varphi$. 
We further define the semi-norm 
\begin{align*}
    |g|_\al:=\sup_{m\geq 0}\{\al^{-m}\|\osc_m(g)\|_{L^1(\mu_\vp)}\}.
\end{align*}
We will consider the space 
\begin{align*}
    \cB_\al:=\{g\in L^1(\mu_\vp): |g|_\al<\infty\}, 
\end{align*}
which can be considered a Banach space together with the norm 
\begin{align*}
    \|g\|_\al:=\|g\|_{L^1(\mu_\vp)}+|g|_{\al}.
\end{align*}
For each $n\geq 0$ we let $H_n\sub E_A^\infty$ be a measurable set such that $(H_{n})_{n=0}^\infty$ is a nested sequence of open sets (holes), i.e. $H_{n}\sub H_{n-1}$. We assume the following:
\begin{enumerate}
    \item[\mylabel{U0}{U0}] $H_{0}=E_A^\infty$, 
    \item[\mylabel{U1}{U1}] For each $n\geq 1$, the open set $H_{n}$ is a disjoint union of cylinders each of which is of length $n$, 
    \item[\mylabel{U2}{U2}] There exists $\vrho\in(0,1)$ such that $\mu_\vp(H_{n})\leq \vrho^n$ for each $n\geq 0$.
\end{enumerate}
Increasing $\alpha$ if necessary, we assume from now on that $\varrho<\alpha<1$.
The perturbation estimates in \cite{pollicott_open_2017} use the weak norm $\|\cdot\|_*$ defined by
\begin{align}\label{star norm}
    \|g\|_{*}:=\|g\|_{L^1(\mu_\vp)}+|g|_{*},
\end{align}
where the semi-norm $|\cdot|_{*}$ is given by
\begin{align*}
    |g|_{*}:=\sup_{j\geq 0}\sup_{m\geq 0}\lt\{\al^{-m}\int_{\vta^{-j}(H_{m})}|g|\, d\mu_\vp\rt\}.
\end{align*}
We take $\cB=(\cB_\al,\|\cdot\|_\al)$,
and use $\|\cdot\|_*$ as the weak norm
in the perturbation estimates. 
For all $g\in\cB_\al$, the norm defined in \eqref{star norm} satisfies the following:
\begin{align*}
    \|g\|_{*}\leq 2\|g\|_\infty\leq 2\|g\|_\al.
\end{align*}
For each $n\in\NN$, the perturbed transfer operator $\cL_{n}:\cB\to\cB$ is naturally defined by 
\begin{align*}
    \cL_{n}(g):=\cL_\vp(g\cdot\ind_{H_{n}^c}).
\end{align*}

We are finally ready to check that the assumptions required here hold in the setting of \cite{pollicott_open_2017} for each sufficiently large fixed
$n\geq 0$.
We begin by noting that in the setting of deterministic dynamics with deterministic holes, where we take $\Omega$ to be a singleton and $\sigma$ to be the identity, our assumptions follow from the work
\cite{pollicott_open_2017}, with \eqref{LL3} omitted as permitted by Remark~\ref{rmk:5.6}.
In fact, \eqref{M1}, \eqref{M2}, \eqref{LL1} (provided \eqref{G} holds), and \eqref{LL2} all hold trivially in the deterministic setting, and \eqref{B1} follows from the definition of $\mathcal B$.
Hence, we are only left to check our assumptions \eqref{C0}--\eqref{C1}, \eqref{B2}--\eqref{B3}, \eqref{G0}--\eqref{G3}, and \eqref{LL4}.
Conditions \eqref{C0}--\eqref{C1} follow from Proposition 2.2.3 of \cite{pollicott_open_2017}, using the normalization above and the non-atomicity of $\mu_\varphi$ for \eqref{C1}.
Condition \eqref{B2} follows from Lemma 2.3.1 of \cite{pollicott_open_2017} together with the estimate
\[
\|fg\|_\alpha\leq 2\|f\|_\alpha\|g\|_\alpha.
\]
Conditions \eqref{B3} and
\eqref{G1}--\eqref{G3} follow from
Lemmas~2.3.3--2.3.4, Proposition~2.4.2,
and Corollary~2.4.4 of
\cite{pollicott_open_2017}.
Once \eqref{G1}--\eqref{G3} are established,
condition \eqref{G0} follows in this
deterministic setting by taking $n=1$
in \eqref{eq:stronggap}.
Finally, \eqref{LL4} follows from \eqref{B3} and \eqref{G}, as explained in Remark~\ref{rmk:5.6}.
Thus, we obtain all of our conditional limit theorems in the deterministic finite-alphabet setting considered by Pollicott and Urba\'nski \cite{pollicott_open_2017}.

\begin{remark}
We used Pollicott and Urba\'nski \cite{pollicott_open_2017} to verify the
hypotheses in this section. For related spectral and perturbative results on
finite-state systems with holes, see Tanaka \cite{tanaka_ruelle_2019,tanaka_finiteholes_2020}.

The work of Pollicott and Urba\'nski also allows countably infinite alphabets.
Countable-state systems with holes have also been studied by Demers et al.\
\cite{demers_countableholes_2017} and Tanaka \cite{tanaka_countableholes_2024}.
To avoid additional technical complications, we restrict $\Xi$ to be compact
and therefore consider only finite-alphabet subshifts of finite type in this
section. An appropriate extension of our arguments to a Polish phase space
$\Xi$ might allow us to treat the countable-alphabet setting of \cite{pollicott_open_2017}.

Our arguments may also extend to deterministic expanding interval maps with countably many branches and small holes.
Furthermore, the inducing schemes developed in Chapter 5 of \cite{pollicott_open_2017} suggest possible extensions to suitable topological Collet--Eckmann interval maps, rational maps on the Riemann sphere, and transcendental meromorphic maps on the complex plane.
\end{remark}

\section{Proof of Theorems   \ref{main thm qstat and qerg measures}, \ref{main thm yaglom_limit}, \ref{thm:c}}
\label{sec: proof of A,B,C}

\subsection{Preliminaries}
We will use the following lemma repeatedly.
\begin{lemma}\label{eq: iterate pert tr op2}
For every $\zeta\in\mathcal P_{\mathcal B,+}(\nu_c)$ with random density $\psi$, every $f,v\in\mathcal B$, every integers $N,n$ with $0\le n\le N$ and $m$-a.e.~$\oio$, we have
\begin{align*}
\zeta_{\om,N}(f\circ T^n_\om \cdot v) =
\frac{ \nu_{\sigma ^n\om,c}(f \cdot \mathcal L^n_\om(v\cdot \psi_\om)\cdot \hat X_{\sigma ^n\om,N-n})}{\nu_{\sigma ^n\om,c}(\mathcal L^n_\om\psi_\om \cdot \hat X_{\sigma ^n\om,N-n})}. 
\end{align*}
\end{lemma}
\begin{proof}
The case $n=0$ follows directly from \eqref{eq:con}.
For $n\ge 1$, by the fact that $\hat X_{\om,N}=\hat X_{\om,n-1}\cdot \hat X_{\sg^n\om,N-n}\circ T^n_\om$ \mas ~and the duality \eqref{eq:duality},
\begin{align*}
\zeta _\om (f\circ T^n_\om \cdot v\cdot \hat X_{\om,N})&=
\nu_{\om,c} \big((f \hat X_{\sg^n\om,N-n})\circ T^n_\om \cdot (v \psi_\om \hat X_{\om,n-1})\big)\\
&= \lambda_{\om,c}^{-n} \nu_{\sg^n\om,c} \big( f \hat X_{\sg^n\om,N-n} \mathcal L_{\om,c}^n( v \psi_\om \hat X_{\om,n-1})\big)
\end{align*}
for $\mas ~\oio$, where $\lambda_{\om,c}^{-n}:=(\lambda_{\om,c}^{n})^{-1}$, $\lambda_{\om,c}^{n}:=\lambda_{\sg^{n-1}\om,c}\lambda_{\sg^{n-2}\om,c}\cdots\lambda_{\om,c}$.
Thus, it follows from \eqref{eq: iterate pert tr op} that
\[
\zeta _\om (f\circ T^n_\om \cdot v\cdot \hat X_{\om,N})=\lambda_{\om,c}^{-n} \nu_{\sg^n\om,c} \big( f  \mathcal L_{\om}^n( v \psi_\om )\hat X_{\sg^n\om,N-n}\big)\quad \mas
\]
In a similar manner, we get
\[
\zeta _\om(X_{\om,N})=\lambda_{\omega,c}^{-n} \nu_{\sigma ^n\om,c}(\mathcal L^n_\om\psi_\om \hat X_{\sigma ^n\om,N-n}) \quad \mas
\]
The claim immediately follows from these equations.
\end{proof}

Recall that $\eta _{\om,0}:=(\phi\cdot \nu_c)_{\om,0}$ $m$-a.s.
Lemma \ref{eq: iterate pert tr op2} with $N=n$,
$\psi_\om=\phi_\om\hat X_{\om,0}/
\nu_{\om,c}(\phi_\om\hat X_{\om,0})$,
$v=\hat X_\om$, and
$\tilde{\mathcal L}_\om\phi_\om=\phi_{\sg\om}$
in \eqref{G} gives, for every $h\in\cB$
and every $n\in\NN$,
\begin{align}\label{eq: cond inv}
    \eta_{\om,n}( h\circ T^n_\om ) = \eta_{\sg^n\om,0}(h)
\end{align}
for $m$-a.e.~$\omega$.
In particular, $\eta_0$ is conditionally invariant.

\subsection{Auxiliary version of Theorem \ref{thm:c}}
\label{s:7.2}

The exponential decaying estimates in Theorems \ref{main thm qstat and qerg measures}, \ref{main thm yaglom_limit}, and \ref{thm:c} rely on a bootstrap argument.
The estimate in Lemma~\ref{lem:cmix} below is almost the same as that in Theorem \ref{thm:c}.
The essential difference is the use of $\zeta_{\omega,\infty}(f)$ instead of $||f||_\infty$.
This replacement is justified by the bound $|\zeta_{\omega,\infty}(f)|\le ||f||_\infty$ for every $f\in\cB$, and the bound can be established once it is shown that $\zeta_{\om,n}(h)\to\zt_{\om,\infty}(h)$ for every $h\in\cB$.
Since Lemma~\ref{lem:cmix} is used exactly for this convergence, we first prove Lemma~\ref{lem:cmix} without relying on any a priori estimates.
\begin{lemma}\label{lem:cmix}
Assume that \eqref{G} holds.  Then, for every such $\zeta$, every $f,v\in\cB$, all integers $N,n$ with $n_{\omega,\zeta}\le n\le N$, and $m$-a.e.~$\omega$,
\begin{align*}
&\left|
\zeta_{\om,N}\left(f\circ T_\om^n\cdot v\right)
-
\eta_{\sg^n\om,N-n}(f)\zeta_{\omega,\infty}(v)
\right| \notag\\
&\qquad\le
\frac{2C_0\left(1+\Vert f\Vert _\infty |\zeta_{\omega,\infty}(v)|\right)
\max\{\Vert v\psi_\omega\Vert \Vert f\Vert _\infty,\Vert \psi_\omega\Vert\}}
{\nu_{\omega,\infty}(\psi_\omega)}
\left(D_{\sg^n\om}C_{\sg^n\om,\phi}\right)
\kappa^n .
\end{align*}
\end{lemma}

Recalling $\psi$ is the random density of $\zeta$, and taking Lemma \ref{eq: iterate pert tr op2} and \eqref{G} into account, for each $\oio$, we introduce 
\begin{align*}
& L_1:= 
  \nu_{\om,\infty}(v\psi_\om) \nu_{\sg^n\om,c}(f \phi_{\sg^n\om} \hat X_{\sg^n\om,N-n}), \\
 &R_1:= \nu_{\sg^n\om,c}(f  (\tilde{\mathcal L}^n_\om (v\psi_\om) -  \nu_{\om,\infty} (v\psi_\om)  \phi _{\sg^n\om}) \hat X_{\sg^n\om,N-n}),\\
&L_2:= \nu_{\om,\infty}(\psi_\om) \nu_{\sg^n\om,c}(\phi_{\sg^n\om} \hat X_{\sg^n\om,N-n})  , \\
&R_2:=  \nu_{\sg^n\om,c}( (\tilde{\mathcal L}^n_\om \psi_\om -\nu_{\om,\infty}(\psi_\om) \phi _{\sg^n\om}) \hat X_{\sg^n\om,N-n}).
\end{align*}
Then Lemma~\ref{eq: iterate pert tr op2} gives
\[
\zeta_{\om,N}(f\circ T^n_\om \cdot v)=\frac{L_1+R_1}{L_2+R_2},
\]
and by the definition of $\zeta_{\omega,\infty}$,
\[
\eta_{\sg^n\om,N-n}(f)\zeta_{\omega,\infty}(v)=\frac{L_1}{L_2}.
\]
Thus $\frac{\vert L_1\vert }{\vert L _2\vert }  \le \Vert f \Vert _\infty |\zt_{\om,\infty}(v)|$ and
\begin{align}\label{eq:20260215a}
\notag\left|
\zeta_{\om,N}(f\circ T^n_\om\cdot v)
-
\eta_{\sg^n\om,N-n}(f)\zeta_{\omega,\infty}(v)
\right| &= \left\vert\frac{L _1 + R_1}{L _2+R_2} - \frac{L _1 }{L _2}\right\vert \\\notag
&\le  \frac{\vert R_1\vert }{\vert L _2+ R_2\vert } +\left\vert \frac{L _1 }{L _2+R_2} - \frac{L _1 }{L _2}\right\vert\notag \\
&\le  \frac{\vert R_1\vert }{\vert L_2 + R_2\vert } + \frac{\vert R_2\vert }{\vert L_2+ R_2\vert } \cdot \frac{\vert L_1\vert }{\vert L_2\vert} \notag\\
&
\le\left(1+\|f\|_\infty |\zeta_{\omega,\infty}(v)|\right)
\frac{\max\{|R_1|,|R_2|\}}{|L_2+R_2|}.
\end{align}

On the one hand, it immediately follows from \eqref{B2} and \eqref{eq:stronggap} 
with $f=v\psi_\om$ %\gf{Where do we assume that \eqref{hyp:G} holds with these choices?} 
that 
\begin{align}\label{eq:20260215b}
\vert R_1\vert &\le \Vert f\Vert _\infty \Vert \tilde{\mathcal L}^n_\om(v \psi_\om) - \phi_{\sg^n\om}  \nu_{\om,\infty} (v\psi_{\omega})\Vert _\infty \nu_{\sg^n\om,c}(\hat X_{\sg^n\om,N-n})\\
&\le C_0  \Vert f\Vert _\infty\Vert \tilde{\mathcal L}^n_\om(v \psi_\om) - \phi _{\sg^n\om}\nu_{\om,\infty} (v\psi_\om) \Vert   \nu_{\sg^n\om,c}(\hat X_{\sg^n\om,N-n})\notag \\
&\le  C_0D_{\sg^n\om}\Vert v\psi_\om \Vert   \Vert f\Vert_\infty   \kappa ^n  \nu_{\sg^n\om,c}(\hat X_{\sg^n\om,N-n})\notag. 
\end{align}
In a similar manner, we get
\begin{align}\label{eq:20260215c}
\vert R_2\vert \le C_0D_{\sg^n\om}\Vert \psi _\om\Vert    \kappa ^n  \nu_{\sg^n\om,c}(\hat X_{\sg^n\om,N-n}).
\end{align}

We need to make the denominator $L_2$ comparable to $R_1$ and $R_2$, except the exponential decay factor $\kappa^n$, uniformly with respect to $N - n$, which constitutes a key difference from conventional statistical quantities (see also Section~\ref{s:9}). For that purpose, the following is elementary but important.
Recall \eqref{eq:0820} for $C_\phi$.
\begin{lemma}\label{eq:20240618a}
For every $n\in\NN$ and $\mae~\oio$, 
\begin{align}\label{eq:0821b}
\nu_{\om,c}(\hat X_{\om,n}) \le C_{\om,\phi} \nu_{\om,c}(\phi_\om \hat X_{\om,n}).
 \end{align}
 Furthermore, the random variable $C_{\phi}$ is tempered.
 \end{lemma}
\begin{proof}
By \eqref{C1} and \eqref{G},
\begin{align*}
\nu_{\om,c}(\phi_\om \hat X_{\om,n})
&=\lambda_{\om,c}^{-{(n+1)}}\nu_{\sg^{n+1}\om,c}(\mathcal L_{\om,c}^{n+1}(\phi_\om \hat X_{\om,n}))\\
&= \lambda_{\om,c}^{-{(n+1)}} \lambda^{n+1}_\om \nu_{\sg^{n+1}\om,c}(\tilde{\mathcal L}^{n+1}_\om\phi_\om)
=  \lambda_{\om,c}^{-{(n+1)}}\lambda^{n+1}_\om \nu_{\sg^{n+1}\om,c}(\phi_{\sg^{n+1}\om})
\end{align*}
(recall $\lambda_{\om,c}^{n}=\lambda_{\sg^{n-1}\om,c}\lambda_{\sg^{n-2}\om,c}\cdots\lambda_{\om,c}$), and
\begin{align*}
\nu_{\om,c}(\hat X_{\om,n})&=\lambda_{\om,c}^{-{(n+1)}}\nu_{\sg^{n+1}\om,c}(\mathcal L_{\om,c}^{n+1}\hat X_{\om,n})
=\lambda_{\om,c}^{-{(n+1)}}\lambda^{n+1}_\om  \nu_{\sg^{n+1}\om,c}(\tilde{\mathcal L}^{n+1}_\om\hat X_{\om,{0}}),
\end{align*}
so we have
\[
\nu_{\om,c}(\hat X_{\om,n})
=\left(1+\frac{\nu_{\sg^{n+1}\om,c}(\tilde{\mathcal L}^{n+1}_\om(\hat X_{\om,{0}})-\phi_{\sg^{n+1}\om})}{\nu_{\sg^{n+1}\om,c}(\phi_{\sg^{n+1}\om})}\right)\nu_{\om,c}(\phi_\om \hat X_{\om,n}).
\]
Thus,
since $\nu_{\om,\infty}(\hat X_{\om,{0}})=1$ \mas, it follows from \eqref{B2} and \eqref{G} that 
\begin{align*}
\nu_{\om,c}(\hat X_{\om,n})
\le \left(1+  \frac{C_0\Vert \hat X_{\om,{0}}  \Vert D_{\sg^{n+1}\om}  }{\nu_{\sg^{n+1}\om,c}(\phi_{\sg^{n+1}\om})}  \kappa ^{n+1}\right)\nu_{\om,c}(\phi _\om\hat X_{\om,n}).
\end{align*}
 \eqref{eq:0821b} immediately follows.

We next prove that $C_\phi$ is tempered.
Set
\[
a_\omega:=\|\hat X_{\omega,0}\|,
\qquad
b_\omega:=
\frac{D_\omega}{\nu_{\omega,c}(\phi_\omega)}.
\]
By \eqref{B3}, \eqref{G2}, and \eqref{G3},
both $a$ and $b$ are tempered.
Fix $\varepsilon>0$ such that
$e^\varepsilon\kappa<1$.
There are finite positive random variables
$A_\varepsilon$ and $B_\varepsilon$ such that
\[
a_{\sigma^j\omega}
\le A_\varepsilon(\omega)e^{\varepsilon j},
\qquad
b_{\sigma^j\omega}
\le B_\varepsilon(\omega)e^{\varepsilon j}
\]
for all $j\ge0$ and $m$-a.e.~$\omega$.
Consequently, by \eqref{eq:0820},
\begin{align*}
C_{\sigma^j\omega,\phi}
&\le
1+C_0a_{\sigma^j\omega}
\sup_{k\ge0}
b_{\sigma^{j+k}\omega}\kappa^k\\
&\le
1+C_0A_\varepsilon(\omega)
B_\varepsilon(\omega)e^{2\varepsilon j}
\sup_{k\ge0}(e^\varepsilon\kappa)^k\\
&\le
\bigl(1+C_0A_\varepsilon(\omega)
B_\varepsilon(\omega)\bigr)
e^{2\varepsilon j}.
\end{align*}
Since $C_{\omega,\phi}\ge1$ and
$\varepsilon>0$ can be chosen arbitrarily
small, it follows that
\[
\lim_{j\to\infty}
\frac1j\log C_{\sigma^j\omega,\phi}=0
\]
for $m$-a.e.~$\omega$.
Thus $C_\phi$ is tempered.
\end{proof}

By Lemma \ref{eq:20240618a} we have
\[
\vert L_2\vert  = \nu_{\om,\infty}(\psi_\om)\nu_{\sg^n\om,c}(\phi_{\sg^n\om}\hat X_{\sg^n\om,N-n})  \ge  \frac{ \nu_{\om,\infty}(\psi_\om)  \nu_{\sg^n\om,c}( \hat X_{\sg^n\om,N-n})}{C_{\sg^n\om,\phi}}.
\]
Thus, we have
\[
\vert L_2 +R_2\vert  \ge \Big(\frac{ \nu_{\om,\infty}(\psi_\om)}{C_{\sg^n\om,\phi}} -  C_0D_{\sg^n\om} \Vert \psi _\om\Vert    \kappa ^n\Big)  \nu_{\sg^n\om,c}( \hat X_{\sg^n\om,N-n})  \ge \frac{1}{2} \frac{ \nu_{\om,\infty}(\psi_\om)  \nu_{\sg^n\om,c}( \hat X_{\sg^n\om,N-n})}{C_{\sg^n\om,\phi}}
\]
\forevery $n\ge n_{\om,\zt}$ (recall \eqref{eq:nomzt} for $n_{\om,\zt}$).
Therefore,
the desired bound of Lemma \ref{lem:cmix}
immediately follows from \eqref{eq:20260215a}, \eqref{eq:20260215b} and \eqref{eq:20260215c}.

\subsection{Proof of Theorem  \ref{main thm qstat and qerg measures}}
Since
$D$ and $C_\phi$ are tempered by Lemma \ref{eq:20240618a},
it follows from Lemma \ref{lem:cmix} 
that
\begin{equation}\label{eq:20260705b}
\zeta_{\omega,n}(h)\to \zeta_{\omega,\infty}(h),\qquad h\in\cB.
\end{equation}
For $h\in\cB$ satisfying $h\ge 0$, this convergence gives
\begin{align}\label{eq:20260528g}
  \zeta_{\omega,\infty}(h)=\lim_{n\to\infty}\zeta_{\omega,n}(h)\ge0,\qquad
   \zeta_{\omega,\infty}(\hat \Xi)=\lim_{n\to\infty}\zeta_{\omega,n}(\hat \Xi)=1
\end{align}
and
\begin{align}\label{eq:20260528}
   |\zeta_{\omega,\infty}(h)|
   =\lim_{n\to\infty}|\zeta_{\omega,n}(h)|
   \le \Vert h\Vert_\infty,
\end{align}
because each $\zeta_{\omega,n}$ is a probability measure.
By applying \eqref{eq:20260528} to the bound in Lemma \ref{lem:cmix} with $N=n$, $f=\hat X$, and $v=h$, together with \eqref{B2},  
we get item (1).

By \eqref{eq:20260528},
 the restriction of
$\zeta_{\omega,\infty}$ to $\cB\cap C(\Xi)$ is continuous for the supremum norm.
By the density
assumption \eqref{B1}, it extends uniquely to a bounded linear functional
on $C(\Xi)$, 
which is positive and normalized by \eqref{eq:20260528g}.
Therefore, the Riesz--Markov representation theorem gives 
a Borel probability measure 
$\widehat\zeta_{\omega,\infty}$ representing the functional induced by $\zeta_{\om,\infty}$.
In particular,
\begin{align}\label{eq:20260705a}
\zt_{\om,\infty}(f)=\widehat\zeta_{\om,\infty}(f),\qquad f\in\cB\cap C(\Xi)
\end{align}
for $m$-a.e.~$\oio$.
This proves item~(2).

For $\zeta=\nu_c$, the definition of $\zeta_\infty$ gives $\zeta_\infty=\nu_\infty$.
For $\zeta=\eta_0$, we have
\[
\frac{d\eta_{\omega,0}}{d\nu_{\omega,c}}
=
\frac{\phi_\omega\hat X_{\omega,0}}
{\nu_{\omega,c}(\phi_\omega\hat X_{\omega,0})}.
\]
Since $\nu_{\omega,\infty}(h\hat X_{\omega,0}) =\nu_{\omega,\infty}(h)$ for $h\in\cB$ and $\nu_{\omega,\infty}(\phi_\omega)=1$, we obtain $\zeta_\infty=\eta_\infty$.
This proves item~(3).

If $g\in C(\Xi)$ and
$\epsilon>0$, choose $g_0\in\cB\cap C(\Xi)$ with $\Vert g-g_0\Vert_\infty<\epsilon$.  Then, it follows from \eqref{eq:20260528} and \eqref{eq:20260705a} that
\[
\begin{aligned}
&|\zeta_{\omega,n}(g)-\widehat\zeta_{\omega,\infty}(g)|\\
&\le
|\zeta_{\omega,n}(g-g_0)|
+|\zeta_{\omega,n}(g_0)-\widehat\zeta_{\omega,\infty}(g_0)| 
+|\widehat\zeta_{\omega,\infty}(g_0)-\widehat\zeta_{\omega,\infty}(g)|\\
&\le
2\epsilon
+|\zeta_{\omega,n}(g_0)-\widehat\zeta_{\omega,\infty}(g_0)|
=2\epsilon
+|\zeta_{\omega,n}(g_0)-\zeta_{\omega,\infty}(g_0)|,
\end{aligned}
\]
and the last term tends to zero by \eqref{eq:20260705b}.  Hence $\zeta_{\omega,n}(g)\to
\widehat\zeta_{\omega,\infty}(g)$ for every $g\in C(\Xi)$.

For each $k\ge0$, the measures $\zeta_{\omega,n}$
are supported on $X_{\omega,k}$ whenever $n\ge k$.
Thus, by the Portmanteau theorem,
$
\widehat\zeta_{\omega,\infty}
(\overline{X_{\omega,k}})=1.
$
Consequently,
\[
\widehat\zeta_{\omega,\infty}
\Big(\bigcap_{k\ge0}\overline{X_{\omega,k}}\Big)=1.
\]

Moreover, for every Borel set $A\subset\Xi$,
\[
\zeta_{\omega,n}(A)
\le
\frac{\|\psi_\omega\|_\infty}
{\nu_{\omega,n}(\psi_\omega)}
\nu_{\omega,n}(A).
\]
By item~(1) with $\zeta=\nu_c$,
$\nu_{\omega,n}(\psi_\omega)$ converges to
$\nu_{\omega,\infty}(\psi_\omega)>0$.
Hence the coefficient is bounded for all sufficiently large $n$.
Passing to the weak limits in the corresponding inequality
for nonnegative continuous functions yields
$
\widehat\zeta_{\omega,\infty}
\ll\widehat\nu_{\omega,\infty}.
$
This proves item~(4).

\subsection{Proof of Theorem \ref{main thm yaglom_limit}}
For every $h\in\cB$,
\[
T_*^n(\zeta_n)_\omega(h)=\zeta_{\sg^{-n}\om,n}(h\circ T^n_{\sg^{-n}\om}).
\]
If $n\ge n_{\sigma^{-n}\omega,\zeta}$, Lemma~\ref{lem:cmix} with base point $\sigma^{-n}\omega$, $N=n$, $f=h$ and $v=\widehat X_{\sigma^{-n}\omega,0}$, together with \eqref{eq: cond inv} and \eqref{eq:20260528}, yields
\[
\left|T_*^n(\zeta_n)_\omega(h)-\eta_{\omega,0}(h)\right|
\le
\frac{2C_0(1+\|h\|_\infty)^2\|\psi_{\sigma^{-n}\omega}\|}{\nu_{\sigma^{-n}\omega,\infty}(\psi_{\sigma^{-n}\omega})}D_\omega C_{\omega,\phi}\kappa^n,
\]
which is bounded by the asserted right-hand side.
If $n<n_{\sigma^{-n}\omega,\zeta}$, the definition \eqref{eq:nomzt} implies that for some $k\ge n$,
\[
\frac{2C_0\|\psi_{\sigma^{-n}\omega}\|}{\nu_{\sigma^{-n}\omega,\infty}(\psi_{\sigma^{-n}\omega})}D_{\sigma^{k-n}\omega}C_{\sigma^{k-n}\omega,\phi}\kappa^k\ge1.
\]
Writing $k=n+j$, we have $j\ge0$ and $\sigma^{k-n}\omega=\sigma^j\omega$, hence
\[
\frac{2C_0\|\psi_{\sigma^{-n}\omega}\|}{\nu_{\sigma^{-n}\omega,\infty}(\psi_{\sigma^{-n}\omega})}\widetilde C_{\omega,\phi}\kappa^n\ge1.
\]
Both measures on the left-hand side of item~(1) are probability measures, so their difference on $h$ is bounded by $2\|h\|_\infty\le(1+\|h\|_\infty)^2$.
This proves the asserted bound in the second case as well, and therefore proves item~(1).

If $\|\psi\|$ and $\nu_\infty(\psi)^{-1}$ are tempered with respect to $(\sigma^{-1},m)$, then the coefficient on the right-hand side has subexponential growth in $n$, while $\widetilde C_{\omega,\phi}\kappa^n\to 0$ for $m$-a.e.~$\omega$ because $\widetilde C_{\omega,\phi}$ is finite.
Therefore item~(1) implies
\[
    T_*^n(\zeta_n)_\omega(h)\to\eta_{\omega,0}(h), \qquad h\in\mathcal B
\]
for \(m\)-a.e.~\(\omega\). 
Let \(g\in C(\Xi)\) and let \(\varepsilon>0\).  By \eqref{B1}, choose
\(h\in\mathcal B\cap C(\Xi)\) such that
$\|g-h\|_\infty<\varepsilon.$
Since \(T_*^n(\zeta_n)_\omega\) and \(\eta_{\omega,0}\) are probability
measures, we have
\[
\begin{aligned}
&\left| T_*^n(\zeta_n)_\omega(g)-\eta_{\omega,0}(g) \right|\\
&\le
\left| T_*^n(\zeta_n)_\omega(g-h) \right| 
+ \left| T_*^n(\zeta_n)_\omega(h)-\eta_{\omega,0}(h) \right| 
+ \left| \eta_{\omega,0}(h-g) \right|\\
&\le
2\varepsilon
+ \left| T_*^n(\zeta_n)_\omega(h)-\eta_{\omega,0}(h) \right|.
\end{aligned}
\]
Since the last term tends to zero as \(n\to\infty\) and \(\varepsilon>0\) is arbitrary, we obtain 
\[
    T_*^n(\zeta_n)_\omega(g)\to\eta_{\omega,0}(g),
    \qquad g\in C(\Xi)
\]
for \(m\)-a.e.~\(\omega\).
This proves item~(2).

Finally, let 
$\rho\in \mathcal P_{\mathcal B,+}(\nu_c)$ be a conditionally invariant random probability measure, and let $\psi$ be its random density.
Then $T_*^n(\rho_n)=\rho$ for every $n\ge 0$.
Replacing $\omega$ by $\sigma^n\omega$ in item~(1), we obtain
\[
\left|\rho_{\sg^n\om}(h)-\eta_{\sg^n\om,0}(h)\right|
\le
\frac{2C_0(1+\|h\|_\infty)^2\|\psi_\om\|}{\nu_{\om,\infty}(\psi_\om)}
\widetilde C_{\sg^n\omega,\phi}\kappa^n
\]
for every $h\in\cB$, every $n\ge n_{\om,\rho}$, and $m$-a.e.~$\omega$.
Since $\widetilde C_{\phi}$ is tempered, the right-hand side tends to $0$ as $n\to\infty$.
Thus, for every $h\in\cB$,
\[
\left|\rho_{\sg^n\om}(h)-\eta_{\sg^n\om,0}(h)\right|\to 0
\]
for $m$-a.e.~$\omega$.
If we write
$
U_h(\omega):=\left|\rho_\omega(h)-\eta_{\omega,0}(h)\right|,
$
then $U_h\circ\sigma^n\to 0$ $m$-a.s.
Since $m$ is $\sigma$-invariant, for every $\varepsilon>0$,
\[
m(U_h>\varepsilon)=m(U_h\circ\sigma^n>\varepsilon)\to 0,
\]
so $U_h=0$ $m$-a.s.
Choose a countable uniformly dense subset
$\mathcal D\subset\cB\cap C(\Xi)$.
Intersecting the full-measure sets for $h\in\mathcal D$,
we obtain
$\rho_\omega(h)=\eta_{\omega,0}(h)$
for every $h\in\mathcal D$ and $m$-a.e.~$\omega$.
Since both are probability measures,
$\rho_\omega=\eta_{\omega,0}$ for $m$-a.e.~$\omega$.
This proves item~(3).

\subsection{Proof of Theorem \ref{thm:c}}

Theorem~\ref{thm:c} is an immediate consequence of Lemma \ref{lem:cmix} and \eqref{eq:20260528}.

\section{Dynamical variance}\label{ss:variance}
In this section, we record the second moment estimate used in the proof of Theorem~\ref{thm:e1}(2).
The quantity $\Sigma^2(f)$ introduced below will reappear in Section~\ref{s:taylor} as the quadratic coefficient in the Taylor expansion of the leading multiplier.

For $n\geq 1$, write
\[
S_{\omega,n}f:=\sum_{k=0}^{n-1} f_{\sigma^k\omega}\circ T_\omega^k.
\]
Recall that the hypothesis for items (3) and (4) were used in Corollary \ref{cor:CLTLDPXinfty}.
\begin{proposition}\label{prop:0625a}
Assume that the limit law condition \eqref{LL} holds.
Let \(f\in L^\infty(\Omega,\mathcal B)\) be real-valued and satisfy $\eta_{\omega,\infty}(f_\omega)=0$ for $m$-a.e.~$\oio$.
Then the following hold.
\begin{enumerate}
\item
The following limits exist and coincide for $m$-a.e.~$\oio$:
\[
\lim_{n\to\infty}\eta_{\omega,n}\left[\frac{1}{n}\left(S_{\omega,n}f\right)^2\right]=\Sigma^2(f),
\]
where
$\Sigma^2(f)$ is given in \eqref{eq:20260705c}.
\item
For every $\zeta\in\mathcal P_{\mathcal B,+}(\nu_c)$,
\[
\lim_{n\to\infty}\frac{1}{n^2}\zeta_{\omega,n}\bigl((S_{\omega,n}f)^2\bigr)=0
\]
for $m$-a.e.~$\oio$.
\end{enumerate}
Assume in addition that for $m$-a.e.~$\oio$, every function in $\cB$ is $\nu_{\om,\infty}$-a.e.~continuous. Then the following hold.
\begin{enumerate}
\item[(3)]
For every $k\ge 1$ and $m$-a.e.~$\oio$.
\[
\nu_{\sg^k\omega,\infty}\bigl(f_{\sg^k\omega}\, \tilde{\mathcal L}_\omega^k(\phi_\om f_{\omega})\bigr)
=\widehat \eta_{\omega,\infty}\bigl(f_\omega\cdot f_{\sigma^k\omega}\circ T_\omega^k\bigr).
\]
\item[(4)] 
$\Sigma^2(f)=0$ if and only if there exists $\tilde f\in L^2(\widehat \eta_{\infty})$ such that
\[
f_\omega(x)=\tilde f_\omega(x)-\tilde f_{\sigma\omega}\circ T_\omega(x)
\]
for $m$-a.e.~$\oio$ and
$\widehat\eta_{\omega,\infty}$-a.e.~$x\in X_{\omega,\infty}$.
\end{enumerate}
\end{proposition}

\begin{proof}
For item~(1), we expand the square and split the expression into diagonal and off-diagonal terms:
\begin{align*}
\eta_{\omega,n}\left[\frac{1}{n}\left(S_{\omega,n}f\right)^2\right]
&=\frac{1}{n}\sum_{k=0}^{n-1}\eta_{\omega,n}\bigl((f_{\sigma^k\omega}\circ T_\omega^k)^2\bigr) \\
&\quad +\frac{2}{n}\sum_{i=0}^{n-2}\sum_{j=i+1}^{n-1}
\eta_{\omega,n}\bigl(f_{\sigma^i\omega}\circ T_\omega^i\cdot f_{\sigma^j\omega}\circ T_\omega^j\bigr)
\end{align*}
for $m$-a.e.~$\oio$.

We first consider the diagonal term.
By the conditional invariance of $\eta$ in Theorem~\ref{main thm yaglom_limit},
\[
\eta_{\omega,n}\bigl((f_{\sigma^k\omega}\circ T_\omega^k)^2\bigr)
=
\eta_{\sigma^k\omega,n-k}(f_{\sigma^k\omega}^2)
\]
for $m$-a.e.~$\omega$ and every $0\leq k\leq n-1$.
Note also that \(\|\phi\|\) is essentially bounded. Indeed,
\eqref{eq:stronggap}, \(\nu_{\omega,\infty}(\hat X_{\omega,0})=1\),
and \eqref{LL} give
\[
\|\phi_{\sigma\omega}\|
\le
D\|\hat X_{\omega,0}\|\kappa
+\lambda_\omega^{-1}K_\omega\|\hat X_{\omega,0}\|,
\]
and the right-hand side is essentially bounded.

Lemma \ref{eq: iterate pert tr op2}, applied with \(N=l\), gives
\[
\eta_{\omega,l}(f_\omega^2)
=
\frac{
\nu_{\sigma^l\omega,c}\left(
\hat X_{\sigma^l\omega,0}\,
\widetilde{\mathcal L}_\omega^l(\phi_\omega f_\omega^2)
\right)
}{
\nu_{\sigma^l\omega,c}\left(
\hat X_{\sigma^l\omega,0}\phi_{\sigma^l\omega}
\right)
}.
\]
Since
\(\nu_{\omega,\infty}(\phi_\omega f_\omega^2)
=\eta_{\omega,\infty}(f_\omega^2)\), \eqref{eq:stronggap},
\eqref{B2} and Lemma \ref{eq:20240618a} yield
\[
\begin{aligned}
\left|
\eta_{\omega,l}(f_\omega^2)
-\eta_{\omega,\infty}(f_\omega^2)
\right|
\le
C_0D\|\phi_\omega f_\omega^2\|\kappa^l
\frac{
\nu_{\sigma^l\omega,c}(\hat X_{\sigma^l\omega,0})
}{
\nu_{\sigma^l\omega,c}
(\phi_{\sigma^l\omega}\hat X_{\sigma^l\omega,0})
}
\le
C_0DC_{\sg^l\om,\phi}
\|\phi_\omega f_\omega^2\|\kappa^l.
\end{aligned}
\]
On the other hand, by \eqref{eq:0820}, \eqref{LL1}, and \eqref{LL4},
\begin{equation}\label{eq:LL4-uniform-Cphi}
\esssup_\om C_{\om,\phi}
\le 1+
\frac{C_0D\,\esssup_\om\|\hat X_{\om,0}\|}
     {\essinf_\om\nu_{\om,c}(\phi_\om)}<\infty.
\end{equation}
Therefore, the essential boundedness of \(\|\phi\|, \|f\|\) and \eqref{B2} give a constant
\(C<\infty\) such that
\begin{equation*}
\left|
\eta_{\omega,l}(f_\omega^2)-\eta_{\omega,\infty}(f_\omega^2)
\right|
\le C\kappa^l
\end{equation*}
for every \(l\ge1\) and $m$-a.e.~$\omega$. Consequently,
\[
\left|
\frac{1}{n}\sum_{k=0}^{n-1}\eta_{\sigma^k\omega,n-k}(f_{\sigma^k\omega}^2)
-
\frac{1}{n}\sum_{k=0}^{n-1}\eta_{\sigma^k\omega,\infty}(f_{\sigma^k\omega}^2)
\right|\le
\frac{C}{n}\sum_{r=1}^{n}\kappa^r\longrightarrow0.
\]
Then, by Birkhoff's ergodic theorem,
\[
\frac{1}{n}\sum_{k=0}^{n-1}\eta_{\sigma^k\omega,\infty}(f_{\sigma^k\omega}^2)
\to
\int_\Omega \eta_{\omega,\infty}(f_\omega^2)\,dm(\omega)
\]
for $m$-a.e.~$\omega$.

We next consider the off-diagonal term.
Write $k=j-i$.
Using again the conditional invariance of $\eta$ and Lemma \ref{eq: iterate pert tr op2}, we have
\begin{align*}
\eta_{\omega,n}\bigl(f_{\sigma^i\omega}\circ T_\omega^i \cdot f_{\sigma^{i+k}\omega}\circ T_\omega^{i+k}\bigr)
=
\eta_{\sigma^i\omega,n-i}\bigl(f_{\sigma^i\omega}\cdot f_{\sigma^{i+k}\omega}\circ T_{\sigma^i\omega}^k\bigr)
\end{align*}
for $m$-a.e.~$\omega$.
Thus,
\begin{equation}\label{eq:20260604}
\begin{split} 
&\frac{2}{n}\sum_{i=0}^{n-2}\sum_{j=i+1}^{n-1}
\eta_{\omega,n}\bigl(f_{\sigma^i\omega}\circ T_\omega^i \cdot f_{\sigma^j\omega}\circ T_\omega^j\bigr)\\
&=
2\sum_{k=1}^{n-1}\frac{n-k}{n}
\left(
\frac{1}{n-k}
\sum_{i=0}^{n-k-1}
\eta_{\sigma^i\omega,n-i}\bigl(f_{\sigma^i\omega}\cdot f_{\sigma^{i+k}\omega}\circ T_{\sigma^i\omega}^k\bigr)
\right).
\end{split} 
\end{equation}

We claim that for every fixed $k\ge1$ and for $m$-a.e.~$\omega$,
\begin{align}
\label{eq:fixed-gap-functional-limit}
&\eta_{\omega,l}(f_\omega\cdot f_{\sigma^k\omega}\circ T_\omega^k) \to
a_k(\omega):=\nu_{\sigma^k\omega,\infty}
(
f_{\sigma^k\omega}\,
\widetilde{\mathcal L}_\omega^k(\phi_\omega f_\omega)
)
\qquad (l\to\infty).
\end{align}
Indeed, it follows from Lemma \eqref{eq: iterate pert tr op2} that
\begin{align*}
\eta_{\omega,l}(f_\omega\cdot f_{\sigma^k\omega}\circ T_\omega^k)
=
\frac{
\nu_{\sigma^k\omega,c}(f_{\sigma^k\omega}\,\widetilde{\mathcal L}_\omega^k(\phi_\omega f_\omega)\,\hat X_{\sigma^k\omega,l-k})
}{
\nu_{\sigma^k\omega,c}(\phi_{\sigma^k\omega}\,\hat X_{\sigma^k\omega,l-k})
}
\end{align*}
for all $l\ge k$.
By \eqref{eq:duality} and \eqref{G2}, the numerator can be written as
\[
\lambda_{\sigma^k\omega,c}^{-l+k}\lambda_{\sigma^k\omega}^{l-k}
\nu_{\sigma^l\omega,c}(\hat X_{\sigma^l\omega,0}\,\widetilde{\mathcal L}_{\sigma^k\omega}^{l-k} (f_{\sigma^k\omega}\,\widetilde{\mathcal L}_\omega^k(\phi_\omega f_\omega))),
\]
while the denominator is
\[
\lambda_{\sigma^k\omega,c}^{-l+k}\lambda_{\sigma^k\omega}^{l-k}
\nu_{\sigma^l\omega,c}(\hat X_{\sigma^l\omega,0}\,\phi_{\sigma^l\omega}).
\]
Hence, with $\tilde f:=f_{\sigma^k\omega}\,\widetilde{\mathcal L}_\omega^k(\phi_\omega f_\omega)\in\cB$ for short notation,
\begin{align*}
\eta_{\omega,l}
(
f_\omega\cdot f_{\sigma^k\omega}\circ T_\omega^k
) 
&=
\frac{
\nu_{\sigma^l\omega,c}
\left(
\hat X_{\sigma^l\omega,0}\,
\widetilde{\mathcal L}_{\sigma^k\omega}^{l-k}(\tilde f)
\right)
}{
\nu_{\sigma^l\omega,c}
\left(
\hat X_{\sigma^l\omega,0}\,
\phi_{\sigma^l\omega}
\right)}\\
&=
a_k(\omega)
+
\frac{
\nu_{\sigma^l\omega,c}
(
\hat X_{\sigma^l\omega,0}[
\widetilde{\mathcal L}_{\sigma^k\omega}^{l-k}
(\tilde f) - \phi_{\sg^l\om}\nu_{\sigma^k\omega,\infty}
(\tilde f)) ]
)
}{
\nu_{\sigma^l\omega,c}
(
\hat X_{\sigma^l\omega,0}\,
\phi_{\sigma^l\omega}
)
}.
\end{align*}
By \eqref{G3} and Lemma~\ref{eq:20240618a}, the numerator of 
the second term on the right-hand side is bounded by 
$
C_0DC_{\sigma^l\omega,\phi}
\|\tilde f\|\kappa^{l-k} \,
\nu_{\sigma^l\omega,c}(\hat X_{\sigma^l\omega,0} \,\phi_{\sigma^l\omega})$.
Hence, it follows from \eqref{eq:LL4-uniform-Cphi}, the essential boundedness of $\|f\|$, \eqref{B2} and \eqref{LL} that
\begin{align}\label{eq:20261018b}
|\eta_{\omega,l}(f_\omega\cdot f_{\sigma^k\omega}\circ T_\omega^k)  - a_k(\omega)|
\le C\kappa^{l-k}
\end{align}
with a constant $C>0$.
Therefore, we obtain \eqref{eq:fixed-gap-functional-limit}.

For each fixed $k\ge1$, it follows from \eqref{eq:20261018b} that
\[
 \frac{1}{n-k}\sum_{i=0}^{n-k-1}
 \left|\eta_{\sigma^i\omega,n-i}
 \bigl(f_{\sigma^i\omega}f_{\sigma^{i+k}\omega}\circ T_{\sigma^i\omega}^k\bigr)
       -a_k(\sigma^i\omega)\right|
 \le \frac{C}{n-k}\sum_{l=1}^{n-k}\kappa^l\longrightarrow0.
\]
Consequently,
\[
\frac{1}{n-k}
\sum_{i=0}^{n-k-1}
\eta_{\sigma^i\omega,n-i}
\bigl(
f_{\sigma^i\omega}\cdot
f_{\sigma^{i+k}\omega}\circ T_{\sigma^i\omega}^k
\bigr)
-
\frac{1}{n-k}
\sum_{i=0}^{n-k-1}
a_k(\sigma^i\omega)
\to 0
\]
for \(m\)-a.e. \(\omega\).  By Birkhoff's ergodic theorem,
\[
\frac{1}{n-k}
\sum_{i=0}^{n-k-1}
a_k(\sigma^i\omega)
\to
\int_\Omega a_k(\omega)\,dm(\omega).
\]
Therefore, for every fixed \(K\ge1\),
\begin{align}
\label{eq:fixed-K-offdiag-limit}
&2\sum_{k=1}^{K}\frac{n-k}{n}
\left(
\frac{1}{n-k}
\sum_{i=0}^{n-k-1}
\eta_{\sigma^i\omega,n-i}
\bigl(
f_{\sigma^i\omega}\cdot
f_{\sigma^{i+k}\omega}\circ T_{\sigma^i\omega}^k
\bigr)
\right) \to
2\sum_{k=1}^{K}
\int_\Omega a_k(\omega)\,dm(\omega)
\end{align}
for \(m\)-a.e. \(\omega\).

It remains to control the tail $2\sum_{k=K+1}^{n-1}(\cdots)$ of \eqref{eq:20260604}.
%\eqref{eq:fixed-K-offdiag-limit} in \(k\). 
For that purpose, we prepare some estimates.
Notice that \eqref{LL} 
implies
\[
\esssup_{\omega} C_{\omega,\phi}<\infty.
\]
Moreover, since
$
\nu_{\omega,\infty}(\hat X_{\omega,0})=1,
$
applying \eqref{LL} gives
\begin{align*}
\esssup_{\omega} \|\phi_{\sigma\omega}\|
&\le
\esssup_{\omega} \Big(\|
\widetilde{\mathcal L}_\omega\hat X_{\omega,0}
-
\phi_{\sigma\omega}
\| +\|\widetilde{\mathcal L}_\omega\hat X_{\omega,0}\|\Big)\\
&\le
\esssup_{\omega} \Big(D\|\hat X_{\omega,0}\|\kappa +
\lambda_\omega^{-1}K_\omega\|\hat X_{\omega,0}\|\Big) < \infty .
\end{align*}
Since \(\nu_{\omega,\infty}(\phi_\omega\hat X_{\omega,0}) =\nu_{\omega,\infty}(\phi_\omega)=1\)
by \eqref{G1}, we have
\[
\psi_\omega
=
\frac{\phi_\omega\hat X_{\omega,0}}
     {\nu_{\omega,c}(\phi_\omega\hat X_{\omega,0})},
     \qquad
\nu_{\omega,\infty}(\psi_\omega)
=
\frac{1}{\nu_{\omega,c}(\phi_\omega\hat X_{\omega,0})}.
\]
Hence
\[
\frac{\|\psi_\omega\|}{\nu_{\omega,\infty}(\psi_\omega)}
=
\|\phi_\omega\hat X_{\omega,0}\|,
\qquad
\frac{\|v\psi_\omega\|}{\nu_{\omega,\infty}(\psi_\omega)}
=
\|v\phi_\omega\hat X_{\omega,0}\|.
\]
Thus, since \(f\in L^\infty(\Omega,\mathcal B)\), the constants
appearing in Theorem~\ref{thm:c}, with \(\zeta=\eta_0\),
\(v=f_\omega\), and \(f=f_{\sigma^k\omega}\), are bounded uniformly in
\(\omega\) and \(k\).  Also, the thresholds \(n_{\omega,\eta_0}\) are
bounded uniformly in \(\omega\).  Hence there exist constants
\(C_f<\infty\) and \(k_f\ge1\) such that, for every \(k\ge k_f\), every
\(l\ge k\), and \(m\)-a.e. \(\omega\),
\begin{equation}
\label{eq:uniform-finite-gap-bound}
\left|
\eta_{\omega,l}
\left(
f_\omega\cdot f_{\sigma^k\omega}\circ T_\omega^k
\right)
\right|
\le
C_f\kappa^k .
\end{equation}
Indeed, this is exactly Theorem~\ref{thm:c} applied with
\(\zeta=\eta_0\).  The product term in Theorem~\ref{thm:c} vanishes
because
$
\eta_{\omega,\infty}(f_\omega)=0$ for \(m\)-a.e.~$\omega$.
Passing to the limit \(l\to\infty\) in
\eqref{eq:uniform-finite-gap-bound} and using
\eqref{eq:fixed-gap-functional-limit}, we also obtain
\begin{equation}
\label{eq:uniform-ak-bound}
|a_k(\omega)|
\le
C_f\kappa^k
\qquad
\text{for every }k\ge k_f
\text{ and for \(m\)-a.e. }\omega .
\end{equation}

Let \(K\ge k_f\).  By
\eqref{eq:uniform-finite-gap-bound},
\begin{align*}
&2\sum_{k=K+1}^{n-1}\frac{n-k}{n}
\left|
\frac{1}{n-k}
\sum_{i=0}^{n-k-1}
\eta_{\sigma^i\omega,n-i}
\bigl(
f_{\sigma^i\omega}\cdot
f_{\sigma^{i+k}\omega}\circ T_{\sigma^i\omega}^k
\bigr)
\right|\\
&\le
2\sum_{k=K+1}^{n-1}\frac{n-k}{n}
\left(
\frac{1}{n-k}
\sum_{i=0}^{n-k-1}
C_f\kappa^k
\right)\le
2C_f\sum_{k=K+1}^{\infty}\kappa^k .
\end{align*}
Hence
\begin{equation*}
\lim_{K\to\infty}
\limsup_{n\to\infty}
2\sum_{k=K+1}^{n-1}\frac{n-k}{n}
\left|
\frac{1}{n-k}
\sum_{i=0}^{n-k-1}
\eta_{\sigma^i\omega,n-i}
\bigl(
f_{\sigma^i\omega}\cdot
f_{\sigma^{i+k}\omega}\circ T_{\sigma^i\omega}^k
\bigr)
\right|
=0.
\end{equation*}
Similarly, by \eqref{eq:uniform-ak-bound},
\[
\sum_{k=K+1}^{\infty}
\left|
\int_\Omega a_k(\omega)\,dm(\omega)
\right|
\le
C_f\sum_{k=K+1}^{\infty}\kappa^k
\to 0
\qquad (K\to\infty).
\]
It follows from
\eqref{eq:fixed-K-offdiag-limit} and the above tail estimates that
\begin{align*}
&\lim_{n\to\infty}
\frac{2}{n}\sum_{i=0}^{n-2}\sum_{j=i+1}^{n-1}
\eta_{\omega,n}
\bigl(
f_{\sigma^i\omega}\circ T_\omega^i
\cdot
f_{\sigma^j\omega}\circ T_\omega^j
\bigr)=
2\sum_{k=1}^{\infty}
\int_\Omega a_k(\omega)\,dm(\omega)
\end{align*}
for \(m\)-a.e. \(\omega\).
Combining the diagonal and off-diagonal parts proves item~(1).

For item~(2), fix \(\zeta\in\mathcal P_{\mathcal B,+}(\nu_c)\), and let \(\psi\) be its random density. 
For $\max\{1,n_{\omega,\zeta}\}\le i<j\le n-1$, repeated application of
Lemma~\ref{eq: iterate pert tr op2} gives
\begin{align*}
\zeta_{\omega,n}
 (f_{\sigma^i\omega}\circ T_\omega^i\,
  f_{\sigma^j\omega}\circ T_\omega^j)
  &=\nu_{\omega,\infty}(\psi_\omega)
\frac{\nu_{\sigma^j\omega,c}\!\left(
 \hat X_{\sigma^j\omega,n-j}f_{\sigma^j\omega}
 \widetilde{\mathcal L}_{\sigma^i\omega}^{j-i}
   (f_{\sigma^i\omega}\phi_{\sigma^i\omega})\right)}
 {\nu_{\sigma^j\omega,c}
   (\hat X_{\sigma^j\omega,n-j}\widetilde{\mathcal L}_\omega^j\psi_\omega)}\\
&\qquad+
\frac{\nu_{\sigma^j\omega,c}\!\left(
 \hat X_{\sigma^j\omega,n-j}f_{\sigma^j\omega}
 \widetilde{\mathcal L}_{\sigma^i\omega}^{j-i}\!\left(
  f_{\sigma^i\omega}\left[
   \widetilde{\mathcal L}_\omega^i\psi_\omega
    -\nu_{\omega,\infty}(\psi_\omega)\phi_{\sigma^i\omega}
  \right]\right)\right)}
 {\nu_{\sigma^j\omega,c}
   (\hat X_{\sigma^j\omega,n-j}\widetilde{\mathcal L}_\omega^j\psi_\omega)}.
\end{align*}
The denominator estimate after Lemma~\ref{eq:20240618a} is
\[
 \nu_{\sigma^j\omega,c}
 (\hat X_{\sigma^j\omega,n-j}\widetilde{\mathcal L}_\omega^j\psi_\omega)
 \ge
 \frac{\nu_{\omega,\infty}(\psi_\omega)}{2C_{\sigma^j\omega,\phi}}
 \nu_{\sigma^j\omega,c}(\hat X_{\sigma^j\omega,n-j}).
\]
For the first term, the centering condition gives
$\nu_{\sigma^i\omega,\infty}
 (f_{\sigma^i\omega}\phi_{\sigma^i\omega})=0$,
thus \eqref{LL1} yields
\[
 \left\|\widetilde{\mathcal L}_{\sigma^i\omega}^{j-i}
       (f_{\sigma^i\omega}\phi_{\sigma^i\omega})\right\|
 \le D\kappa^{j-i}\|f_{\sigma^i\omega}\phi_{\sigma^i\omega}\|.
\]
For the second term, \eqref{LL1} gives
\[
 \left\|\widetilde{\mathcal L}_\omega^i\psi_\omega
   -\nu_{\omega,\infty}(\psi_\omega)\phi_{\sigma^i\omega}\right\|
 \le D\kappa^i\|\psi_\omega\|.
\]
Therefore, one can find a constant $C>0$, depending only on $C_0$, $D$ and the essential suprema of $\|\phi_\om\|$, $\|f_\om\|$ and $C_{\om,\phi}$, such that
\begin{equation*}
 \left|\zeta_{\omega,n}
 (f_{\sigma^i\omega}\circ T_\omega^i\,
  f_{\sigma^j\omega}\circ T_\omega^j)\right|
 \le
 C\Big(
       \kappa^{j-i}+\frac{\|\psi_\omega\|}{\nu_{\omega,\infty}(\psi_\omega)}\kappa^i\Big).
\end{equation*}
(A similar argument to the one used to obtain the constant $C>0$ in \eqref{eq:20261018b} 
applies here.)
In the square of $S_{\omega,n}f$, the diagonal terms and the pairs with
$i<\max\{1,n_{\omega,\zeta}\}$ contribute at most
$(2n_{\omega,\zeta}+3)nC_0^2\esssup_{\oio}\|f_\om\|^2$.
 For the other pairs,
 \begin{align*}
2C\sum_{0\le i<j<n}\Big(
       \kappa^{j-i}+\frac{\|\psi_\omega\|}{\nu_{\omega,\infty}(\psi_\omega)}\kappa^i\Big)
\le \frac{2Cn}{1-\kappa}\Big(1+\frac{\|\psi_\omega\|}{\nu_{\omega,\infty}(\psi_\omega)}\Big).
\end{align*}
Thus
\begin{align*}
\frac{\zeta_{\omega,n}((S_{\omega,n}f)^2)}{n^2}
\le\frac{(2n_{\omega,\zeta}+3)C_0^2
                 \esssup_{\oio}\|f_\om\|^2}{n}
+  \frac{2C}{(1-\kappa)n}\Big(1+\frac{\|\psi_\omega\|}{\nu_{\omega,\infty}(\psi_\omega)}\Big)
\longrightarrow0.
\end{align*}
This proves item~(2).

For item~(3), apply 
Lemma~\ref{lem:borel-formulas-on-survivor} with
$n=k$, $u=\phi_\omega f_\omega$, and
$v=f_{\sigma^k\omega}$.  Together with
\eqref{eq:representing-measure-densities}  and the argument below it, this gives
\begin{align*}
&\widehat\eta_{\omega,\infty}
\bigl(f_\omega f_{\sigma^k\omega}\circ T_\omega^k\bigr)
=
\widehat\nu_{\omega,\infty}
\bigl(\phi_\omega f_\omega
      f_{\sigma^k\omega}\circ T_\omega^k\bigr)
=
\widehat\nu_{\sigma^k\omega,\infty}
\left(
 f_{\sigma^k\omega}\,
 \widetilde{\mathcal L}_\omega^k(\phi_\omega f_\omega)
\right).
\end{align*}
By \eqref{B2}, the function inside the last integral belongs to
$\mathcal B$.  Equation~\eqref{eq:zeta-functional-measure-identification}, applied
with $\zeta=\nu_c$, therefore identifies that integral with
\[
\nu_{\sigma^k\omega,\infty}
\left(
 f_{\sigma^k\omega}\,
 \widetilde{\mathcal L}_\omega^k(\phi_\omega f_\omega)
\right),
\]
which proves item~(3).

Finally, item~(4) is a statement about the closed random dynamical system $T:X_\infty\to X_\infty$ with its random invariant probability measure $\widehat\eta_\infty$.
This is exactly the setting used in the proof of Proposition~3 of \cite{DFGTV18a}.
Therefore item~(4) follows from the proof of \cite[Proposition~3]{DFGTV18a}.
\end{proof}

\section{Proof of Theorem  \ref{thm:e1}}\label{s:thmD}
\label{sec: proof of D}

\subsection{Proof of Theorem \ref{thm:e1}(1)}

We first show that $\widehat\eta_\infty$ is the unique conditionally ergodic limit over $\mathcal P_{\mathcal B,+}(\nu_c)$.
Fix $\zeta\in\mathcal P_{\mathcal B,+}(\nu_c)$ and $h\in\mathcal B\cap C(\Xi)$.
For every $N\ge1$, by the definition of pushforward,
\begin{align*}
\frac1N\sum_{j=0}^{N-1} (T_*^j(\zeta_{N-1}))_\omega(h)=\frac1N\sum_{j=0}^{N-1}
\zeta_{\sigma^{-j}\omega,N-1}
\bigl(h\circ T_{\sigma^{-j}\omega}^j\bigr).
\end{align*}
For each $j$ satisfying $j\ge n_{\sigma^{-j}\omega,\zeta}$, Theorem \ref{thm:c}, applied with base point $\sigma^{-j}\omega$, conditioning time $N-1$, and $v=\hat X_{\sigma^{-j}\omega}$, gives
\begin{align}\label{eq:20261110b}
\left| \zeta_{\sigma^{-j}\omega,N-1} \bigl(h\circ T_{\sigma^{-j}\omega}^j\bigr) -\eta_{\omega,N-1-j}(h) \right|
\le
C_{\sigma^{-j}\omega,\zeta,h,\hat X_{\sigma^{-j}\omega}} D_\omega C_{\omega,\phi}\kappa^j.
\end{align}
Here we used $\zeta_{\sigma^{-j}\omega,\infty}(\hat X_{\sigma^{-j}\omega})=1$.
For the remaining $j$, the left-hand side is bounded by $2\Vert h\Vert_\infty$.

For $M\ge1$, define
\[
E_M:=\left\{\oio : \max\{ n_{\omega,\zeta}, C_{\omega,\zeta,h,\hat X_\omega}\}>M \right\}.
\]
For $j\ge M$ with $\sigma^{-j}\omega\notin E_M$, the right-hand side of \eqref{eq:20261110b} is at most $M D_\omega C_{\omega,\phi}\kappa^j$.
Consequently, for $m$-a.e.~$\omega$,
\begin{align*}
&\limsup_{N\to\infty}\frac1N\sum_{j=0}^{N-1} \left| \zeta_{\sigma^{-j}\omega,N-1} \bigl(h\circ T_{\sigma^{-j}\omega}^j\bigr) -\eta_{\omega,N-1-j}(h) \right|\\
&\quad\le 2\Vert h\Vert_\infty \lim_{N\to\infty}\frac1N \sum_{j=0}^{N-1}\widehat{E_M}(\sigma^{-j}\omega) =2\Vert h\Vert_\infty m(E_M)\longrightarrow 0,
\end{align*}
where the equality follows from Birkhoff's ergodic theorem for $(\sigma^{-1},m)$
and the convergence holds because both $n_{\zeta}$ and $C_{\zeta,h,\hat X_\omega}$ are finite $m$-a.s.

On the other hand, Theorem~\ref{main thm qstat and qerg measures}(4) and the Ces\`aro theorem yield
\[
\frac1N\sum_{j=0}^{N-1}\eta_{\omega,N-1-j}(h) =\frac1N\sum_{k=0}^{N-1}\eta_{\omega,k}(h) \longrightarrow \widehat\eta_{\omega,\infty}(h).
\]
We have therefore proved the desired convergence against every $h\in\mathcal B\cap C(\Xi)$.
To pass to $C(\Xi)$, choose a countable subset
$\{h_k:k\ge1\}\subset\mathcal B\cap C(\Xi)$ that is dense in
$C(\Xi)$ in the uniform norm, using \eqref{B1} and the separability of
$C(\Xi)$.
For the fixed $\zeta$, the preceding convergence holds for every
$h_k$ on a common full $m$-measure set.
For any $h\in C(\Xi)$ and any $k\ge1$, since all measures
involved are probabilities,
\[
\begin{aligned}
\left|\frac1N\sum_{j=0}^{N-1}
(T_*^j\zeta_{N-1})_\omega(h)
-\widehat\eta_{\omega,\infty}(h)\right|\le 2\|h-h_k\|_\infty
+\left|\frac1N\sum_{j=0}^{N-1}
(T_*^j\zeta_{N-1})_\omega(h_k)
-\widehat\eta_{\omega,\infty}(h_k)\right|.
\end{aligned}
\]
Letting first $N\to\infty$ and then approximating $h$
uniformly by $h_k$ proves the weak convergence for
all $h\in C(\Xi)$ on the same full-measure set.

We now prove the weak convergence
\[
\frac1N\sum_{j=0}^{N-1} \bigl(T_*^j(\zeta_{N-1})\bigr)_{\sigma^j\omega}\longrightarrow \overline \eta_{\infty}
\]
for every $\zt\in\mathcal P_{\mathcal B,+}(\nu_c)$ and $m$-a.e.~$\oio$.
Fix $h\in\mathcal B\cap C(\Xi)$.
Theorem~\ref{thm:c} with $f=h$ and $v=\hat X_{\omega,0}$ yields
\[
\left|\zeta_{\omega,N-1}(h\circ T_\omega^j) -\eta_{\sigma^j\omega,N-1-j}(h)\right| \le C_{\omega,\zeta,h,\hat X_{\omega,0}} D_{\sigma^j\omega}C_{\sigma^j\omega,\phi}\kappa^j
\]
for $n_{\omega,\zeta}\le j\le N-1$.
For the remaining indices the difference is at most $2\|h\|_\infty$.
Consequently, since $DC_\phi$ is tempered,
we have
\begin{equation}\label{eq:20261110}
\frac1N\sum_{j=0}^{N-1} \left|\zeta_{\omega,N-1}(h\circ T_\omega^j) -\eta_{\sigma^j\omega,N-1-j}(h)\right| \longrightarrow0.
\end{equation}

For each integer $L\ge1$, put
\[
a_L(\omega):=\sup_{k\ge L} \left|\eta_{\omega,k}(h)-\widehat\eta_{\omega,\infty}(h)\right|.
\]
By Theorem~\ref{main thm qstat and qerg measures}(4), $a_L\downarrow0$ for $m$-a.e.~$\omega$, while $0\le a_L\le2\|h\|_\infty$.
Thus $\int a_L\,dm\to0$ by dominated convergence.
For $N>L$,
\begin{align*}
\frac1N\sum_{j=0}^{N-1} \left|\eta_{\sigma^j\omega,N-1-j}(h) -\widehat\eta_{\sigma^j\omega,\infty}(h)\right|
\le\frac1N\sum_{j=0}^{N-1}a_L(\sigma^j\omega) +\frac{2L\|h\|_\infty}{N}.
\end{align*}
By Birkhoff's ergodic theorem for $(\sigma,m)$, the limsup of the left-hand side is at most $\int a_L\,dm$ for $m$-a.e.~$\omega$.
Letting $L\to\infty$, we obtain
\[
\frac1N\sum_{j=0}^{N-1} \left|\eta_{\sigma^j\omega,N-1-j}(h) -\widehat\eta_{\sigma^j\omega,\infty}(h)\right| \longrightarrow0.
\]
Therefore, by \eqref{eq:20261110} and Birkhoff's ergodic theorem, we have
\[
\frac1N\sum_{j=0}^{N-1} \zeta_{\omega,N-1}(h\circ T_\omega^j) \longrightarrow \int_\Omega\widehat\eta_{\omega,\infty}(h)\,dm(\omega) = \overline\eta_\infty(h).
\]
Thus, the claimed weak convergence follows by taking a countable uniformly dense subset of $\mathcal B\cap C(\Xi)$, as in the preceding argument.

\subsection{Proof of Theorem \ref{thm:e1}(3)}
Let $\psi$ be the random density of $\zeta\in\mathcal P_{\mathcal B,+}(\nu_c)$.
Fix $\omega$ in a full-measure set on which Theorems \ref{main thm qstat and qerg measures}, \ref{main thm yaglom_limit} and \ref{thm:c} hold, on which the Birkhoff theorem can be applied to the indicator functions used below, and on which the tempered estimates for $D$, $C_\phi$ and $\nu_c(\phi)^{-1}$ hold.
For $N\ge 1$ and $x\in X_{\omega,N}$, write
\[
A_{\omega,N}(x):=\frac1N S_{\omega,N}f(x)
=\frac1N\sum_{j=0}^{N-1}f_{\sigma^j\omega}\circ T_\omega^j(x),
\]
and let
\[
Z_\omega:=\bigotimes_{N\ge 1}\zeta_{\omega,N}.
\]

We first reduce to the $\eta_{\omega,\infty}$-centered case.
Set
\[
\widetilde f_\omega:=f_\omega-\eta_{\omega,\infty}(f_\omega).
\]
Since $\eta_{\omega,\infty}$ is a probability measure, \eqref{B2} and $f\in L^\infty(\Omega,\cB)$ imply $\widetilde f\in L^\infty(\Omega,\cB)$.
Moreover, for every $N\ge 1$ and every $(x_N)_{N\ge 1}\in\prod_{N\ge 1}X_{\omega,N}$,
\[
\frac1N\sum_{j=0}^{N-1}\widetilde f_{\sigma^j\omega}\circ T_\omega^j(x_N)
=
A_{\omega,N}(x_N)-\frac1N\sum_{j=0}^{N-1}\eta_{\sigma^j\omega,\infty}(f_{\sigma^j\omega}).
\]
By Birkhoff's theorem,
\[
\frac1N\sum_{j=0}^{N-1}\eta_{\sigma^j\omega,\infty}(f_{\sigma^j\omega})
\to\int_\Omega\eta_{\omega,\infty}(f_\omega)\,dm(\omega).
\]
Thus it suffices to prove the theorem with $f$ replaced by $\widetilde f$.
From now on we rename $\widetilde f$ as $f$, and assume that
\[
\eta_{\omega,\infty}(f_\omega)=0
\qquad\text{for $m$-a.e.~$\oio$.}
\]

\medskip

A guiding idea is the following.
For a fixed $N$, under $\zeta_{\omega,N}$, the quantity $f_{\sigma^j\omega}\circ T_\omega^j$ is observed under the condition that the orbit survives the remaining $N-j$ iterates.
For that reason, the natural reference value at time $j$ is
$
\eta_{\sigma^j\omega,N-j}(f_{\sigma^j\omega}),
$
not $\eta_{\sigma^j\omega,\infty}(f_{\sigma^j\omega})$.
This leads to the random triangular array
\[
Y_{\omega,N,j}
:=
f_{\sigma^j\omega}\circ T_\omega^j-
\eta_{\sigma^j\omega,N-j}(f_{\sigma^j\omega}).
\]
We will estimate the fourth moment of $\sum_jY_{\omega,N,j}$ after first restricting to times for which the shifted coefficients are bounded.
Once the fourth power is expanded, the terms to estimate are exactly
$
\zeta_{\omega,N}(Y_{\omega,N,a_1}Y_{\omega,N,a_2}Y_{\omega,N,a_3}Y_{\omega,N,a_4})
$
with $0\le a_1\le a_2\le a_3\le a_4\le N-1$.
The advantage of the above centering is that, in these four-point terms, the leading eigendirection disappears exactly at the last step.
That exact cancellation forces decay from the last gap $d_3:=a_4-a_3$.

There is one additional point.
The coefficients in \eqref{eq:stronggap} are tempered, but they need not be bounded along the whole forward orbit of $\omega$.
We therefore work first in the window $0\le j\le (1-\rho)N$ and only at times for which the relevant shifted coefficients are at most $b$.
The remaining times are removed by Birkhoff's theorem applied to indicators of level sets, and the terminal window is removed at the end by letting $\rho\downarrow0$.

\medskip
For $m\ge 0$, define
\[
\bar f_{\omega,m}:=f_\omega-\eta_{\omega,m}(f_\omega).
\]
Then $Y_{\omega,N,j}=\bar f_{\sigma^j\omega,N-j}\circ T_\omega^j$ for each $0\le j\le N-1$, and
\begin{equation}\label{eq:qslln-centered-sum}
A_{\omega,N}(x)
=
\frac1N\sum_{j=0}^{N-1}Y_{\omega,N,j}(x)
+
\frac1N\sum_{j=0}^{N-1}\eta_{\sigma^j\omega,N-j}(f_{\sigma^j\omega})
\qquad (x\in X_{\omega,N}).
\end{equation}
Since $\eta_{\omega,m}$ is a probability measure and \eqref{B2} holds, there is a constant $K_f$ depending only on $C_0$ and $\|f\|_{L^\infty(\Omega,\cB)}$ such that
\[
\|\bar f_{\omega,m}\|\le K_f,
\qquad
\|Y_{\omega,N,j}\|_\infty\le K_f
\]
for all $\omega,m,N,j$.

Fix $0<\rho<1$ and $b>1$.
Put $N_\rho:=\lfloor(1-\rho)N\rfloor$ and
\[
I_N:=
\left\{n_{\omega,\zeta}\le j\le N_\rho:
D_{\sigma^j\omega}\le b,
\ C_{\sigma^j\omega,\phi}\le b,
\ \|\phi_{\sigma^j\omega}\|\le b
\right\}.
\]
The dependence of $I_N$ on $b$, $\rho$ and $\omega$ is suppressed.

The key estimate is the following four-point bound.

\begin{lemma}\label{lem:qslln-4pt}
For the fixed $b$ and $\rho$ above, there exist a random variable $K_\omega>0$ and a number $\beta=\beta(\rho)>0$ such that, for every $N\ge1$ and every ordered quadruple
$
0\le a_1\le a_2\le a_3\le a_4\le N_\rho
$
with $a_r\in I_N$ for $1\le r\le4$, we have
\begin{equation}\label{eq:qslln-4pt-bound}
\left|
\zeta_{\omega,N}
\bigl(
Y_{\omega,N,a_1}Y_{\omega,N,a_2}Y_{\omega,N,a_3}Y_{\omega,N,a_4}
\bigr)
\right|
\le
K_\omega\bigl(\kappa^{a_1}+\kappa^{d_1}+e^{-\beta N}\bigr)\kappa^{d_3},
\end{equation}
where $d_r:=a_{r+1}-a_r$ for $1\le r\le3$.
Here $K_\omega$ may depend on $b$, $\rho$, $\zeta$ and $f$, but not on $N$ or the ordered quadruple.
\end{lemma}

\begin{proof}
Let
\[
m_r:=N-a_r,
\qquad
\bar f_r:=\bar f_{\sigma^{a_r}\omega,m_r}
\qquad (1\le r\le4).
\]
Thus $Y_{\omega,N,a_r}=\bar f_r\circ T_\omega^{a_r}$.
We shall use \eqref{eq:stronggap} also for zero iterates; when the time is in $I_N$, the case $n=0$ follows from
$h=\nu_{\tau,\infty}(h)\phi_\tau+(h-\nu_{\tau,\infty}(h)\phi_\tau)$, from \eqref{B2}, and from $\|\phi_\tau\|\le b$.

First, for $a\in I_N$ and $h\in\cB$, the denominator estimate used in the proof of Theorem \ref{thm:c}, namely \eqref{eq:20260215c}, Lemma~\ref{eq:20240618a} and \eqref{eq:nomzt}, gives
\[
\left|
\frac{
\nu_{\sigma^a\omega,c}\left(h\,\hat X_{\sigma^a\omega,N-a}\right)
}{
\nu_{\sigma^a\omega,c}\left(\tilde{\mathcal L}_\omega^a\psi_\omega\,\hat X_{\sigma^a\omega,N-a}\right)
}
\right|
\le
\frac{2C_0C_{\sigma^a\omega,\phi}}
{\nu_{\omega,\infty}(\psi_\omega)}\|h\|
\le K_\omega\|h\|.
\]

We also need an estimate for the centering constants.
Let $a\in I_N$ and $m=N-a$.
By the definition of $\eta_{\sigma^a\omega,m}$ and by
$\tilde{\mathcal L}_{\sigma^a\omega}^m\phi_{\sigma^a\omega}=\phi_{\sigma^{a+m}\omega}$,
\[
\eta_{\sigma^a\omega,m}(f_{\sigma^a\omega})-
\eta_{\sigma^a\omega,\infty}(f_{\sigma^a\omega})
=
\frac{
\nu_{\sigma^{a+m}\omega,c}
\left(
\tilde{\mathcal L}_{\sigma^a\omega}^m
\left[
\left(f_{\sigma^a\omega}-\eta_{\sigma^a\omega,\infty}(f_{\sigma^a\omega})\right)
\phi_{\sigma^a\omega}
\right]
\right)
}
{\nu_{\sigma^{a+m}\omega,c}(\phi_{\sigma^{a+m}\omega})}.
\]
The function inside the square brackets has zero $\nu_{\sigma^a\omega,\infty}$-average.
Hence \eqref{eq:stronggap} and \eqref{B2} give
\[
\left|
\eta_{\sigma^a\omega,m}(f_{\sigma^a\omega})-
\eta_{\sigma^a\omega,\infty}(f_{\sigma^a\omega})
\right|
\le
K_fD_{\sigma^a\omega}\|\phi_{\sigma^a\omega}\|
\nu_{\sigma^{a+m}\omega,c}(\phi_{\sigma^{a+m}\omega})^{-1}\kappa^m.
\]
Since $a\in I_N$ and $a\le N_\rho$, while $\nu_c(\phi)^{-1}$ is tempered and $a+m=N$, there are $K_\omega>0$ and $\beta=\beta(\rho)>0$ such that
\[
\left|
\eta_{\sigma^a\omega,N-a}(f_{\sigma^a\omega})-
\eta_{\sigma^a\omega,\infty}(f_{\sigma^a\omega})
\right|
\le K_\omega e^{-\beta N}.
\]

Repeated use of Lemma \ref{eq: iterate pert tr op2} gives the exact identity
\begin{align*}
&\zeta_{\omega,N}
\bigl(
Y_{\omega,N,a_1}Y_{\omega,N,a_2}Y_{\omega,N,a_3}Y_{\omega,N,a_4}
\bigr) \\
&\quad=
\frac{
\nu_{\sigma^{a_4}\omega,c}
\left(
\bar f_4\,
\tilde{\mathcal L}_{\sigma^{a_3}\omega}^{d_3}
\left[
\bar f_3\,
\tilde{\mathcal L}_{\sigma^{a_2}\omega}^{d_2}
\left[
\bar f_2\,
\tilde{\mathcal L}_{\sigma^{a_1}\omega}^{d_1}
\left[
\bar f_1\,\tilde{\mathcal L}_\omega^{a_1}\psi_\omega
\right]
\right]
\right]
\hat X_{\sigma^{a_4}\omega,m_4}
\right)
}{
\nu_{\sigma^{a_4}\omega,c}
\left(
\tilde{\mathcal L}_\omega^{a_4}\psi_\omega\,
\hat X_{\sigma^{a_4}\omega,m_4}
\right)
}.
\end{align*}
We now estimate the expression inside the last application of
$\tilde{\mathcal L}_{\sigma^{a_3}\omega}^{d_3}$.
By \eqref{eq:stronggap},
\[
\tilde{\mathcal L}_\omega^{a_1}\psi_\omega
=
\nu_{\omega,\infty}(\psi_\omega)\phi_{\sigma^{a_1}\omega}+r_0,
\qquad
\|r_0\|\le K_\omega\kappa^{a_1}.
\]
Using the centered assumption and the finite-time estimate above,
\begin{align*}
&\left|
\nu_{\sigma^{a_1}\omega,\infty}
\left(\bar f_1\,\tilde{\mathcal L}_\omega^{a_1}\psi_\omega\right)
\right| \\
&\quad\le
|\nu_{\omega,\infty}(\psi_\omega)|
\left|
\eta_{\sigma^{a_1}\omega,m_1}(f_{\sigma^{a_1}\omega})-
\eta_{\sigma^{a_1}\omega,\infty}(f_{\sigma^{a_1}\omega})
\right|
+
\left|
\nu_{\sigma^{a_1}\omega,\infty}(\bar f_1r_0)
\right| \\
&\quad\le K_\omega\left(e^{-\beta N}+\kappa^{a_1}\right).
\end{align*}
Also $\|\bar f_1\,\tilde{\mathcal L}_\omega^{a_1}\psi_\omega\|\le K_\omega$.
Applying \eqref{eq:stronggap} from $a_1$ to $a_2$, multiplying by $\bar f_2$, and using $a_1,a_2\in I_N$, we obtain
\[
\left\|
\bar f_2\,
\tilde{\mathcal L}_{\sigma^{a_1}\omega}^{d_1}
\left[
\bar f_1\,\tilde{\mathcal L}_\omega^{a_1}\psi_\omega
\right]
\right\|
\le
K_\omega\left(\kappa^{a_1}+\kappa^{d_1}+e^{-\beta N}\right).
\]
Applying \eqref{eq:stronggap} from $a_2$ to $a_3$ and multiplying by $\bar f_3$ gives
\[
\left\|
\bar f_3\,
\tilde{\mathcal L}_{\sigma^{a_2}\omega}^{d_2}
\left[
\bar f_2\,
\tilde{\mathcal L}_{\sigma^{a_1}\omega}^{d_1}
\left[
\bar f_1\,\tilde{\mathcal L}_\omega^{a_1}\psi_\omega
\right]
\right]
\right\|
\le
K_\omega\left(\kappa^{a_1}+\kappa^{d_1}+e^{-\beta N}\right).
\]
One last application of \eqref{eq:stronggap}, now from $a_3$ to $a_4$, shows that the difference between the last pushed term and its leading eigendirection has norm at most
\[
K_\omega\left(\kappa^{a_1}+\kappa^{d_1}+e^{-\beta N}\right)\kappa^{d_3},
\]
since $D_{\sigma^{a_4}\omega}\le b$ by $a_4\in I_N$.
The leading eigendirection gives no contribution to the numerator.  Indeed, by the definition of $\eta_{\sigma^{a_4}\omega,m_4}$,
\begin{align*}
&\nu_{\sigma^{a_4}\omega,c}
\left(
\bar f_4\phi_{\sigma^{a_4}\omega}
\hat X_{\sigma^{a_4}\omega,m_4}
\right) \\
&\quad=
\nu_{\sigma^{a_4}\omega,c}
\left(
 f_{\sigma^{a_4}\omega}\phi_{\sigma^{a_4}\omega}
\hat X_{\sigma^{a_4}\omega,m_4}
\right)
-
\eta_{\sigma^{a_4}\omega,m_4}(f_{\sigma^{a_4}\omega})
\nu_{\sigma^{a_4}\omega,c}
\left(
\phi_{\sigma^{a_4}\omega}
\hat X_{\sigma^{a_4}\omega,m_4}
\right)=0.
\end{align*}
Using the denominator estimate with $a=a_4$ and then \eqref{B2}, we obtain \eqref{eq:qslln-4pt-bound}.
\end{proof}

\begin{lemma}\label{prop:qslln-4th}
For the fixed $b$ and $\rho$ above, there exists a random variable $K_\omega>0$ such that, for every $N\ge1$,
\begin{equation}\label{eq:qslln-4th-moment}
\zeta_{\omega,N}
\left(
\left(\sum_{j\in I_N}Y_{\omega,N,j}\right)^4
\right)
\le K_\omega N^2.
\end{equation}
\end{lemma}

\begin{proof}
Expand the fourth power and order the indices.
At most $4!$ tuples correspond to the same ordered quadruple.
Thus Lemma \ref{lem:qslln-4pt} gives
\begin{align*}
\zeta_{\omega,N}
\left(
\left(\sum_{j\in I_N}Y_{\omega,N,j}\right)^4
\right)
&\le
24K_\omega
\sum_{0\le a_1\le a_2\le a_3\le a_4\le N_\rho}
\left(\kappa^{a_1}+\kappa^{d_1}+e^{-\beta N}\right)\kappa^{d_3}.
\end{align*}
The first two sums are bounded by a constant multiple of $N^2$.  Indeed,
\[
\sum_{0\le a_1\le a_2\le a_3\le a_4\le N_\rho}\kappa^{a_1+d_3}
\le
\sum_{0\le a_2\le a_3\le N_\rho}
\left(\sum_{a_1=0}^{a_2}\kappa^{a_1}\right)
\left(\sum_{a_4=a_3}^{N_\rho}\kappa^{a_4-a_3}\right)
\le K N^2,
\]
and
\[
\sum_{0\le a_1\le a_2\le a_3\le a_4\le N_\rho}\kappa^{d_1+d_3}
\le
\sum_{0\le a_2\le a_3\le N_\rho}
\left(\sum_{a_1=0}^{a_2}\kappa^{a_2-a_1}\right)
\left(\sum_{a_4=a_3}^{N_\rho}\kappa^{a_4-a_3}\right)
\le K N^2.
\]
The remaining sum is at most $K e^{-\beta N}N^4$, hence also at most $K N^2$ after increasing $K$.
This proves \eqref{eq:qslln-4th-moment}.
\end{proof}

We now complete the proof of Theorem \ref{thm:e1}(3).

\begin{proof}
By Theorem~\ref{main thm qstat and qerg measures}(1), for $m$-a.e.~$\omega$,
\[
\eta_{\omega,\ell}(f_\omega)\to \eta_{\omega,\infty}(f_\omega)=0
\qquad (\ell\to\infty).
\]
Therefore, by the same diagonal Ces\`aro argument as in the proof of Proposition~\ref{prop:0625a}(1),
\begin{equation}\label{eq:qslln-deterministic-term}
\frac1N\sum_{j=0}^{N-1}\eta_{\sigma^j\omega,N-j}(f_{\sigma^j\omega})\to 0
\qquad\text{for $m$-a.e.~$\oio$.}
\end{equation}

Fix rational $b>1$ and rational $0<\rho<1$.
By Lemma~\ref{prop:qslln-4th} and Markov's inequality, for every $\delta>0$,
\[
Z_\omega
\left(
\left|
\frac1N\sum_{j\in I_N}Y_{\omega,N,j}(x_N)
\right|>\delta
\right)
\le \frac{K_\omega}{\delta^4N^2}.
\]
The right-hand side is summable in $N$, so the Borel--Cantelli lemma gives
\[
\frac1N\sum_{j\in I_N}Y_{\omega,N,j}(x_N)\to0
\]
for $Z_\omega$-a.e.~$(x_N)_{N\ge1}$.
Taking the intersection over rational $b>1$ and rational $0<\rho<1$, we may use this convergence for all such pairs on the same full $Z_\omega$-measure set.

It remains to remove the discarded times and then the terminal window.
For $b>1$, let
\[
B_b:=\{D>b\}\cup\{C_\phi>b\}\cup\{\|\phi\|>b\}.
\]
Since $D$, $C_\phi$ and $\|\phi\|$ are finite $m$-a.s., $m(B_b)\to0$ as $b\to\infty$.
For $n_{\omega,\zeta}\le j\le N_\rho$, the definition of $I_N$ implies
$j\notin I_N$ if and only if $\sigma^j\omega\in B_b$.
Splitting the sum accordingly,
\[
\frac1N\sum_{j=n_{\omega,\zeta}}^{N_\rho}Y_{\omega,N,j}(x_N)
=
\frac1N\sum_{j\in I_N}Y_{\omega,N,j}(x_N)
+
\frac1N
\sum_{\substack{n_{\omega,\zeta}\le j\le N_\rho\\ \sigma^j\omega\in B_b}}
Y_{\omega,N,j}(x_N).
\]
The first sum tends to $0$ by the Borel--Cantelli argument above.
For the second, $\|Y_{\omega,N,j}\|_\infty\le K_f$ gives
the upper bound by
$K_f\cdot\frac1N\sum_{j=0}^{N_\rho}1_{B_b}(\sigma^j\omega)$.
Since $\frac{N_\rho}{N}\to1-\rho$, Birkhoff's theorem applied to $1_{B_b}$ yields
$
\frac1N\sum_{j=0}^{N_\rho}1_{B_b}(\sigma^j\omega)
 \to (1-\rho)\,m(B_b)\le m(B_b).
$
Combining these,
\[
\limsup_{N\to\infty}
\left|
\frac1N\sum_{j=0}^{N_\rho}Y_{\omega,N,j}(x_N)
\right|
\le K_f\,m(B_b),
\]
where the finitely many indices $j<n_{\omega,\zeta}$ contribute $O(1/N)$ after division by $N$.
Letting $b\to\infty$ through rational values gives
\[
\frac1N\sum_{j=0}^{N_\rho}Y_{\omega,N,j}(x_N)\to0
\]
for every rational $0<\rho<1$.
Finally,
\[
\left|
\frac1N\sum_{j=N_\rho+1}^{N-1}Y_{\omega,N,j}(x_N)
\right|
\le K_f\frac{N-N_\rho}{N}
\le K_f\rho+o(1).
\]
Letting $\rho\downarrow0$ through rational values yields
\begin{equation}\label{eq:qslln-fluctuation-limit}
\frac1N\sum_{j=0}^{N-1}Y_{\omega,N,j}(x_N)\to0
\qquad\text{for $Z_\omega$-a.e.~$(x_N)_{N\ge1}$.}
\end{equation}
Combining \eqref{eq:qslln-centered-sum}, \eqref{eq:qslln-deterministic-term}, and \eqref{eq:qslln-fluctuation-limit}, we get
\[
A_{\omega,N}(x_N)\to0
\qquad\text{for $Z_\omega$-a.e.~$(x_N)_{N\ge1}$}
\]
in the centered case.
Undoing the $\eta_{\omega,\infty}$-centering completes the proof of Theorem \ref{thm:e1}(3).
\end{proof}

\subsection{Proof of Theorem \ref{thm:e1}(2)}
Fix $\zeta\in\mathcal P_{\mathcal B,+}(\nu_c)$, a real-valued $f\in L^\infty(\Omega,\mathcal B)$, and $\delta>0$.
Use $A_{\omega,N}$ and $Z_\omega$ from the proof of item~(3).
For $m$-a.e.~$\omega$, item~(3) gives $A_{\omega,N}(x_N)\to \int_\Omega\eta_{\omega,\infty}(f_\omega)\,dm(\omega)$
for $Z_\omega$-a.e.~$(x_N)_{N\ge1}$.
Hence, by dominated convergence,
\[
\begin{aligned}
\zeta_{\omega,N}\left(\left|\frac1N S_{\omega,N}f-\int_\Omega\eta_{\omega,\infty}(f_\omega)\,dm(\omega)\right|>\delta\right)
=\int\mathbf1_{\left\{\left|A_{\omega,N}(x_N)-\int_\Omega\eta_{\omega,\infty}(f_\omega)\,dm(\omega)\right|>\delta\right\}}\,dZ_\omega\longrightarrow0.
\end{aligned}
\]
This proves item~(2).

\begin{remark}
Under \eqref{LL}, item~(2) also follows directly from the second-moment estimate in Proposition~\ref{prop:0625a}, without using item~(3).
Indeed, set $\tilde f_\omega:=f_\omega-\eta_{\omega,\infty}(f_\omega)$.
By \eqref{B2}, $\tilde f\in L^\infty(\Omega,\mathcal B)$ and $\eta_{\omega,\infty}(\tilde f_\omega)=0$.
Birkhoff's ergodic theorem gives
\[
\frac1n\sum_{j=0}^{n-1}\eta_{\sigma^j\omega,\infty}(f_{\sigma^j\omega})
\longrightarrow\int_\Omega\eta_{\omega,\infty}(f_\omega)\,dm(\omega)
\]
for $m$-a.e.~$\omega$. 
Hence, for every $\delta>0$ and all sufficiently large $n$, Chebyshev's inequality and Proposition~\ref{prop:0625a}(2) yield
\[
\begin{aligned}
\zeta_{\omega,n}\left(
\left|\frac1n S_{\omega,n}f-\int_\Omega\eta_{\omega,\infty}(f_\omega)\,dm(\omega)\right|>\delta
\right)\le\frac{4}{\delta^2n^2}
\zeta_{\omega,n}\bigl((S_{\omega,n}\tilde f)^2\bigr)
\longrightarrow0.
\end{aligned}
\]
\end{remark}

\section{Taylor expansion of the leading multiplier}
\label{s:taylor}

\subsection{Twisted open transfer operators cocycle}
\label{ss:twisted}
For each $\theta\in\mathbb C$, we define a twisted open transfer operator cocycle $\tilde\cL^\theta$ by
\[
\tilde\cL_\om^\theta (\varphi) := \tilde{\mathcal L}_\omega (e^{\theta f_\omega}\cdot \varphi)=\lambda_{ \omega   }^{-1}\mathcal L_{\omega ,c}(e^{\theta f_\omega} \hat X_{\omega,0} \cdot \varphi), \qquad \varphi\in \mathcal B .
\]
Let $\zt\in\cP_{\cB,+}(\nu_c)$, with density $\psi=\frac{d\zeta}{d\nu_c}$.
Then, for $N\ge n\ge 1$, by the duality
\[
\nu _{\omega ,c} ( v \circ T_\omega  \cdot u )  =\lambda_{\om,c}^{-1}\nu _{\sigma \omega ,c}( v \cdot  \mathcal L_{\omega ,c} u), 
\]
one has 
\begin{align*}
 \zeta_{\om,N}(e^{\theta (S_{n}f)_\om})
  &=\frac{\nu _{\om,c}(e^{\theta (S_{n}f)_\om} \hat X_{\om,N} \cdot \psi _\om)}{\nu _{\om, c}(  \hat X_{\om,N} \cdot \psi_\om)} \\
 &= \frac{\nu _{\om,c}((e^{\theta (S_{n-1}f)_{\sg\om}}  \hat X_{\sg\om,N-1} ) \circ T_\om \cdot (e^{\theta f_\om }\hat X_{\om,0}  \psi_\om ))}{\nu _{\om,c}( \hat X_{\sg\om,N-1}  \circ T_\om  \cdot (\hat X_{\om,0} \psi_\om ))}\\
  &= \frac{\lambda_{\om,c}^{-1} \, \nu _{\sg\om,c}(e^{(\theta S_{n-1}f)_{\sg\om}}  \hat X_{\sg\om,N-1}  \cdot  \lambda_{\om}\tilde{\mathcal L}_{\om}^{\theta}  \psi _{\om}) }{\lambda_{\om,c}^{-1}\, \nu _{\sg\om,c}(\hat X_{\sg\om,N-1}  \cdot \lambda _{\om}\tilde{\mathcal L}_\om  \psi_\om ) }\\
      &= \frac{ \nu _{\sg\om,c}(e^{(\theta S_{n-1}f)_{\sg\om}}  \hat X_{\sg\om,N-1}  \cdot   \tilde{\mathcal L}_{\om}^{\theta}  \psi _{\om}) }{ \nu _{\sg\om,c}(\hat X_{\sg\om,N-1}  \cdot  \tilde{\mathcal L}_\om  \psi_\om ) },
\end{align*}
where $(S_{n}f)_\om:=\sum_{j=0}^{n-1}f_{\sg^j\om} \circ T^j_\om$.
Repeating this procedure, with the notation
  \[
  \tilde{\mathcal L}_{\omega}^{\theta ,n}=\tilde{\mathcal L}_{\sigma^{n-1}\omega}^\theta \tilde{\mathcal L}_{\sigma^{n-2}\omega}^\theta  \cdots \tilde{\mathcal L}_{\omega}^\theta , 
  \]
  we get that
\begin{equation}\label{eq:1001a-noomega}
\zeta_{\om,N}(e^{\theta (S_nf)_\om})=
 \frac{\nu_{\sg^n\om,c}( \hat{X}_{\sg^n\om,N-n}\tilde{\mathcal L}^{\theta,n}_\om  \psi _\om)}{\nu_{\sg^n\om,c}(  \hat{X}_{\sg^n\om,N-n} \tilde\cL^n_\om  \psi_\om ) }.
\end{equation}

\subsection{Taylor expansion of the leading random eigenvalue}

The goal of this section is to obtain information about the Taylor expansion around $\theta=0$ of $\lm_\om^\theta$, the leading Lyapunov multiplier of $\tilde\cL_\om^\theta$. We do this by adapting the spectral approach presented in \cite{DFGTV18b} to the current setting.

For a Banach space $B$, we denote by $B^*$ its dual Banach space (equipped with the operator norm) and by $L^\infty (\Omega ,B)$ 
the $L^\infty$-normed Bochner space consisting of $B$-valued functions on $(\Omega, \mathcal F,m)$.

\begin{theorem}\label{conj:1}
Assume \eqref{LL}, and let
$f\in L^\infty(\Omega,\cB)$ be real-valued with
$\eta_{\omega,\infty}(f_\omega)=0$
for $m$-a.e.~$\omega$.
Then, there are an open set $U \subset \mathbb C$ containing $0$,  $\mathscr C^1$ maps   $\phi: U \to L^\infty(\Omega, \cB)$, $\nu_\infty: U \to L^\infty(\Omega, \cB^*)$, and a  $\mathscr C^3$ map  $\lambda : U \to L^\infty (m)$   such that the following hold with the notation $\phi ^\theta _\omega  , \nu ^\theta_{\infty , \omega} , \lambda _{\theta ,\omega }$ for $\phi (\theta )( \omega) , \nu_\infty(\theta )( \omega) , \lambda (\theta )(\omega )$, respectively:
\begin{enumerate}
\item
 For every $\theta \in U $, $h\in\mathcal B$ and $m$-a.e.~$\oio$,
 \[
\tilde{\mathcal L} ^\theta _\om\phi ^\theta_\om   = \lambda^\theta _{\om} \phi ^\theta _{\sg\om},
\quad 
\nu_{\sg\om,\infty}^\theta (\tilde{ \mathcal L} ^\theta_{\om} h )= \lambda ^\theta _{\om} \nu_{\om,\infty} ^\theta (h)
\]
and
\[
 \nu_{\om,\infty} (\phi^\theta _\om)=  \nu^\theta _{\om,\infty}(\phi^\theta _\om) =1. 
\]
\item
For every $\theta \in U $,
there is an $\tilde{\mathcal L}_\omega ^\theta$-equivariant splitting $\mathcal B =\mathbb C\phi _{\omega}^\theta \oplus  \mathrm{Ker}( \nu^\theta _{\omega,\infty}) $, where $\mathrm{Ker}(\nu^\theta _{\omega,\infty})=\{u \in \mathcal B : 
\nu^\theta _{\omega,\infty}(u) = 0\}$, that is, $\tilde{\mathcal L}^\theta_{\omega} (\mathrm{Ker}(\nu^\theta _{\omega,\infty})) \subset \mathrm{Ker}(\nu^\theta_{\sigma\omega ,\infty})$ for $m$-a.e.~$\oio$.
Let  $\widehat{\Lambda} : U \to L^\infty (m)$  be given by 
\[
\widehat{\Lambda } (\theta)(\omega)  = \log \lambda ^\theta_\omega .
\]
Then, for every $\theta\in U\cap\mathbb R$ and $m$-a.e.~$\oio$,
 \[
\lim_{n\to\infty} \frac{1}{n}\log \Vert \tilde{\mathcal L} ^{\theta ,n}_\om \vert _{\mathbb C\phi ^\theta_\om}\Vert   = \int \widehat{\Lambda } (\theta ) \,dm. 
\]
Moreover, there is $\kappa _1\in (0,1)$ such that 
 for every $\theta \in U $ and $m$-a.e.~$\oio$,
\[
\lim_{n\to\infty} \frac{1}{n}\log \Vert \tilde{\mathcal L} ^{\theta ,n}_\om \vert _{\mathrm{Ker}(\nu^\theta_\om)}\Vert  
\le \log \kappa _1.
\]
Replacing $U$ by a smaller neighborhood of $0$, there are constants $C_2>0$ and $r\in (0,1)$ 
such that for every $\theta \in U$, $n\ge 1$ and $m$-a.e.~$\oio$,
\[
 \Vert \tilde\cL_\om^{\theta ,n} \vert_{\mathrm{Ker}(\nu^\theta_{\om,\infty})}\Vert
\le C_2 r^n.
\]
Consequently,
\begin{align*}
   \norm{\tilde\cL ^{\theta ,n}_\om h -\lambda^{\theta, n}_{\om} \nu_{\om,\infty}^\theta (h)\phi ^\theta _{\sg^n\om}}
   \le C_2\norm{h} r^n \quad \mas
\end{align*}
for every $\theta\in U$, $n\ge 1$, $h\in\mathcal B$ and $m$-a.e.~$\oio$, where $\lambda^{\theta, n}_\om=\lambda_{\sg^{n-1}\om}^\theta \cdots \lambda _{\om}^\theta$.
\item
It holds that
\[
\widehat{\Lambda }  (0)(\omega) =\widehat{\Lambda } ' (0)(\omega)=  0 
\]
for $m$-a.e.~$\oio$,
and 
\[ 
\int  \widehat{ \Lambda } ^{\prime \prime} (0)  \,dm =  \Sigma ^2(f) .
\] 
Moreover, we have the following $\omega$-uniform estimate on the remainder:
\[
\lim_{\theta \to 0}\left\Vert\frac{\widehat{ \Lambda }(\theta )}{\theta ^2} - \frac{\widehat{ \Lambda } ^{\prime \prime}(0)}{2} \right\Vert _{L^\infty} =0
\] 
\item
It holds that
\[
\phi ^0_\om  =\phi_\om , \quad \nu_{\om,\infty}^0  =\nu_{\om,\infty}   \quad  \text{and} \quad \lambda ^0_\om =1
\]
for $m$-a.e.~$\oio$.
Furthermore, there is a constant $C>0$ such that
for any $\theta \in U$, integer $k\ge 0$ and $\mathcal B$-valued random variable $\psi$ satisfying $\nu_\infty (\psi )>0$ $\mas$, 
\[
\left\vert \frac{\nu_{\om,c}(\hat X_{\om,k} \phi^\theta _\om) }{ \nu_{ \om,c}(\hat X_{\om,k} \phi _\om)}-1 \right\vert \le C |\theta|  \quad \text{and} \quad \left\vert \frac{\nu_{\om,\infty}^\theta (\psi_\om) }{ \nu_{\om,\infty}(\psi _\om)} -1 \right\vert \le C|\theta| \frac{||\psi _\om||}{\nu_{\om,\infty}(\psi _\om)} 
\]
for $m$-a.e.~$\oio$.
\end{enumerate}
\end{theorem}

We finally introduce the random variable given in Theorem \ref{be}:
\begin{equation}\label{eq:0805b}
\sg_{\om,n}^2(f)=
\begin{cases}
    \sum_{j=0}^{n-1}\hat{\Lambda}''(0)(\sigma^{j}\omega)\quad &(\text{if $\sum_{j=0}^{n-1}\hat{\Lambda}''(0)(\sigma^{j}\omega)> 0$})\\
    n\Sigma^2(f) \quad &(\text{otherwise})
    \end{cases}.
\end{equation}
Notice that it follows from Theorem \ref{conj:1}(3) and Birkhoff ergodic theorem that $\frac{\sg_n^2(f)}{n}$ converges to $\Sigma^2 (f)$ $\mas$

\subsection{Outline of the proof of Theorem \ref{conj:1}}
\label{ss:implicit}

The proof of Theorem \ref{conj:1} is an adaptation of the spectral perturbation argument of \cite{DFGTV18b}; see also \cite{CNW,AFGTV-TFPF}. 
Since the hypotheses in \cite{DFGTV18b} are formulated in terms of variation-type estimates, whereas here they are stated for an abstract Banach algebra, we first record the precise comparison needed below.

\begin{itemize}
    \item \textbf{Condition \cite[(C0)]{DFGTV18b}: measurability for the multiplicative ergodic theorem.}
    This condition is used to apply a multiplicative ergodic theorem.
    %In \cite{DFGTV18b},
    \cite[(C0)]{DFGTV18b} is the topological measurability condition used to invoke the theorem of Froyland--Lloyd--Quas \cite{FLQ2}.
    As shown in \cite{CNW}, the alternative measurable formulation based on the theorem of Gonz\'alez-Tokman--Quas \cite{GTQ14} is sufficient for the same purpose.
    The two alternatives are precisely those recorded in \eqref{LL2}.

\item \textbf{Condition \cite[(C1)]{DFGTV18b}: boundedness of the cocycle on \(\mathcal B\).}
This follows immediately from \eqref{LL1} and the definition of \(\tilde{\mathcal L}_\omega\).
    \item \textbf{Condition \cite[(C2)]{DFGTV18b}: Lasota--Yorke estimate.}
    This is supplied by the Lasota--Yorke part of \eqref{LL3}. 
\item \textbf{Condition \cite[(C3)]{DFGTV18b}: contraction on the complement of the top direction.}
   This immediately follows from \eqref{LL1}:
    \eqref{eq:stronggap} in \eqref{LL1} 
gives
    $       \bigl\|\tilde{\mathcal L}_\omega^n\varphi\bigr\|       \leq        D\kappa^n\|\varphi\|$
    for each $n\geq 0$ and $\varphi\in\mathcal B$ satisfying
    \(\nu_{\omega,\infty}(\varphi)=0\).
    \item \textbf{Condition \cite[(C4)]{DFGTV18b}: positivity.} 
We do not assume this condition.
In \cite{DFGTV18b}, (C4) is a uniform positivity condition from which the unperturbed density has a positive pointwise lower bound; see \cite[Lemma~2.11(2), equation~(18)]{DFGTV18b}.
In our notation, the corresponding pointwise condition would be $\essinf_\om\inf_{x\in X_\om}\phi_\om(x)>0$, which we do not require; compare \eqref{LL4}.
This lower bound is used to obtain nonnegativity of the perturbed densities $\phi_\omega^\theta$ and hence positivity of the corresponding multipliers $\lambda_\omega^\theta$ for real $\theta$ sufficiently close to $0$; see \cite[Proposition~3.6 and equation~(43)]{DFGTV18b}.
This positivity allows the use of real logarithms when identifying the limiting logarithmic moment generating function with $\int_\Omega\log\lambda_\omega^\theta\,dm(\om)$ for real $\theta$, before applying the G\"artner--Ellis theorem in the proof of the LDP; see \cite[Section~4.1]{DFGTV18b}.
Here, however, we will obtain $\lambda_\omega^\theta>0$ without this uniform pointwise lower bound, as briefly explained below.

For the spectral construction below, the denominators $\mathbf I\circ\mathbf L_\theta(\varphi)$ and $\mathbf J\circ\mathbf L_\theta^*(u)$ equal $1$ at $(\theta,\varphi)=(0,0)$ and $(\theta,u)=(0,0)$, respectively,
and remain uniformly away from zero in neighbourhoods of these points.
The quotient derivatives therefore involve only bounded powers of these scalar reciprocals. 
The inverses of the linearized equations follow from the uniform gap, as do the bounds on the complementary iterates.
None of these steps divides by $\phi_\omega(x)$.
For real $\theta$, uniqueness of the normalized implicit solutions gives real $\phi_\omega^\theta$ and $\lambda_\omega^\theta$.
The fact that $\lambda_\omega^0=1$ and the $L^\infty(m)$-continuity of $\theta\mapsto\lambda^\theta$ then give $\lambda_\omega^\theta>1/2$ on a common smaller neighbourhood of $0$, for $m$-a.e.\ $\omega$.
Thus $\log\lambda_\omega^\theta$ is real and essentially bounded, without using $\phi_\omega^\theta>0$.
The details of these arguments are given below.
For a preceding abstract formulation without (C4), see \cite{CNW}.
\end{itemize}

Besides these cocycle assumptions, one also needs the function space properties that are recorded in \cite{DFGTV18b} as $V1$--$V9$. For the present argument, their roles are as follows.

\begin{itemize}
    \item The conditions \(V1\) and \(V2\) are part of the Banach space setting. 
    \item The product estimates in \(V3\) and \(V8\) are replaced by the Banach algebra property \eqref{B2}.
    \item The density condition $V6$ is \eqref{B1}, and $V5$ corresponds to \eqref{B3}.
    \item The compactness condition in $V4$ is provided by the compact embedding in \eqref{LL3}.
    \item 
    The role of \(V9\) in \cite{DFGTV18b} is to control the regularity of the twist
\(\theta\mapsto \mathcal L_\omega(e^{\theta g_\omega}\cdot)\).
Here this is obtained instead from the Banach algebra property \eqref{B2}.
Indeed, since \(f\in L^\infty(\Omega,\mathcal B)\), the multiplication operators
$
h\mapsto e^{\theta f_\omega}h
$
are uniformly bounded on \(\mathcal B\) for \(\theta\) in a neighbourhood of \(0\), and depend analytically on \(\theta\) in operator norm.
Together with the uniform bound in \eqref{LL1}, this gives the required regularity of
$
\tilde{\mathcal L}_\omega(e^{\theta f_\omega}h).
$
Thus the variation-type axiom \(V9\) is not needed as a separate abstract assumption in the present setting. See also \cite{CNW}.
\end{itemize}

With these identifications, the perturbation argument proceeds as in \cite{DFGTV18b}. 
First, on the primal side, one solves by the implicit function theorem a nonlinear fixed-point equation associated with the twisted cocycle. This produces the perturbed top Lyapunov eigenfunction $\phi^\theta$ and the corresponding Lyapunov eigenvalue $\lambda^\theta$. Next one carries out the same construction on the dual side to obtain the dual eigenfunctionals $\nu_\infty^\theta$. These two objects give the equivariant splitting
\[
\mathcal B=\mathbb C\phi_\omega^\theta\oplus \mathrm{Ker}(\nu_{\omega,\infty}^\theta),
\]
and the spectral gap on the complementary bundle yields the exponential remainder estimate in Theorem~\ref{conj:1}(2). The extra smoothness of $\lambda^\theta$, compared with the $\mathscr C^1$ regularity of $\phi^\theta$ and $\nu_\infty^\theta$, comes from the fact that $\theta\mapsto \tilde{\mathcal L}^\theta$ is analytic. Finally, differentiating the fixed-point equations and the identity
\[
\lambda^\theta=\nu_\infty(\tilde{\mathcal L}^\theta\phi^\theta)
\]
at $\theta=0$ gives the formulas in Theorem~\ref{conj:1}(3)--(4); in particular, the centering condition $\eta_\infty(f)=0$ forces the first derivative to vanish, while the second derivative is identified with $\Sigma^2(f)$.

\medskip

For the interested reader, in what follows, we explain the construction of $\phi^\theta$ and show the formulas of the first and second derivatives of $\lambda^\theta$ in more detail.
Recall that given a map $\mathbf A$ from 
an open subset $U$ of 
a normed vector space $\mathscr X$ to a  normed vector space $\mathscr Y$, the Fr\'echet derivative $D\mathbf A (x )$ of $\mathbf A$ at $x \in U$ is a bounded operator from $\mathscr X$ to $\mathscr Y$ such that
\[
\lim _{\Vert y  \Vert  _{\mathscr X} \to 0} \frac{\Vert \mathbf A (x +y  ) -\mathbf A( x) -  D\mathbf A (x )(y)  \Vert _{\mathscr Y} }{\Vert y \Vert _{\mathscr X}} =0 .
\]
The Fr\'echet derivative  of $D\mathbf A : U\to L(\mathscr X, \mathscr Y) $ at $x \in U$ is denoted by $D^2\mathbf A (x)$, 
 where $L(\mathscr X,\mathscr Y)$ is the space of bounded linear operators from $\mathscr X$ to $\mathscr Y$. 
$\mathbf A$ 
 is said to be of class $\mathscr C^1$ (resp.~$\mathscr C^2$) if the map $D\mathbf A $ 
 (resp.~$D^2\mathbf A$)
  is continuous. 
We also use the notation $\mathbf A^\prime$, $\mathbf A^{\prime \prime}$ for $D\mathbf A $, $D^2\mathbf A$. Higher derivatives are defined analogously.
Below we 
use $U$ as a sufficiently small ball in $\mathbb C$ centered at $0$ on which the  maps given below are well-defined, 
even if it may change between occurrences. 
Furthermore, for simple description, we identify each $\varphi \in \mathscr X$ with a map $t\mapsto t\varphi$ from $\mathbb C$ to $\mathscr X$ for $\mathscr X= L^\infty (\Omega , \mathcal B)$ or $L^\infty (\Omega , \mathbb C)$,  when no confusion can arise.
Moreover, following standard practice, for a $\mathscr C^1$ map $\varphi:U\to\mathscr X$, we write $\varphi^\prime(\theta)(\omega)=\psi(\omega)$ with $\psi\in\mathscr X$ to mean $\varphi^\prime(\theta)(t)(\omega)=t\psi(\omega)$ for every $t\in U$.

  Define $\mathbf I : L^\infty(\Omega,\mathcal B)\to L^\infty (m)$ and $\mathbf L_\theta :  L^\infty(\Omega,\mathcal B) \to L^\infty(\Omega,\mathcal B)$ by
\[
\mathbf I(\varphi ) (\omega ) =  \nu_{\omega,\infty} \left( \varphi _\omega \right), \quad 
 \mathbf L _\theta  (\varphi ) (\omega  )= \tilde{\mathcal L}^\theta _{\sigma ^{-1} \omega }  ( \phi_{\sigma^{-1}\omega}  +\varphi _{\sigma ^{-1} \omega }  )
\]
for   $\varphi \in L^\infty(\Omega,\mathcal B)$ and $\omega \in \Omega$. 
It follows from \eqref{G} that 
$\theta \mapsto \mathbf I \circ \mathbf L_\theta (\phi)$ is a continuous map 
 from $\mathbb C$ to $L^\infty (m)$.
The map $(\theta,\varphi)\mapsto
\mathbf I\circ\mathbf L_\theta(\varphi)$ is holomorphic into $L^\infty(m)$,
and
$
 \mathbf I(\phi)=\mathbf I\circ\mathbf L_0(0)=1.
$
Consequently, for $|\theta|$ and $\|\varphi\|_{L^\infty(\Omega,\mathcal B)}$
sufficiently small,
\[
 \essinf_\omega\left|
 (\mathbf I\circ\mathbf L_\theta(\varphi))(\omega)\right|>\frac12.
\]
On this neighbourhood with $\varphi\in\mathrm{Ker}(\mathbf I)$,
\[
 \mathbf F(\theta ,  \varphi ) :=
 \frac{ \mathbf L_\theta (  \varphi )}{\mathbf I \circ \mathbf L_\theta ( \varphi ) }  -  (\phi + \varphi ) 
\]
is well-defined as an element of 
$L^\infty(\Omega,\mathcal B)$. 
Notice that it is in fact in $\mathrm{Ker} (\mathbf I )$.
It also holds that if $\mathbf F( \theta , h ^\theta ) =0$, then
\begin{equation}
\label{eq:1205c}
\tilde{\mathcal L}^\theta_\omega    \phi  ^\theta_\om  = \lambda ^\theta_{\om} \phi ^\theta_{\sg\om} \quad \text{and} \quad  \nu_{\om,\infty} (\phi^\theta _\om)=1 \quad \text{for $m$-a.e.~$\oio$},
\end{equation}
 where 
 \begin{equation}\label{eq:1101}
\phi  ^\theta _\om= \phi_\om    + h ^{\theta}_\om , \quad  \lambda ^\theta_\om  
=   \nu_{\sg\om,\infty} (\tilde{ \mathcal L} ^\theta _{\om}\phi^\theta  _{\om}) \quad \text{for $\oio$.} 
 \end{equation}

Furthermore, 
we  define $\mathbf J : L^\infty (\Omega , \mathcal B^*)\to L^\infty (m)$ and $\mathbf L_\theta^* :  L^\infty (\Omega , \mathcal B^*) \to L^\infty (\Omega , \mathcal B^*)$ by
\[
\mathbf J(u ) (\omega ) = u_\omega (\phi_ \omega ), \quad 
 \mathbf L _\theta ^* (u) (\omega )(\varphi ) 
 %=   [(\nu_{\infty} +u)(\tilde{\mathcal L}_{\theta }\varphi )]_{\sg\om} 
 =   (\nu_{\infty,\sigma \omega} +u _{\sigma  \omega })(\tilde{\mathcal L}_{\theta ,  \omega }\varphi)
\]
for   $u \in L^\infty (\Omega , \mathcal B^* )$, $\omega \in \Omega$ and $\varphi\in\mathcal B$. 
Since $\mathbf J\circ\mathbf L_0^*(0)(\om) =\nu_{\sg\om,\infty} (\tilde{\mathcal L} _{\om}\phi  _{\om}) =\nu_{\sg\om,\infty} (\phi  _{\sg\om})  =1$,
for $|\theta|$ and $\|u\|_{L^\infty(\Omega,\mathcal B^*)}$ sufficiently small,
\[
 \essinf_\omega\left|
 (\mathbf J\circ\mathbf L_\theta^*(u))(\omega)\right|>\frac12.
\]
Thus the map
\begin{equation*}
 \mathbf G(\theta ,  u ) =
 \frac{ \mathbf L_\theta ^* (  u)}{\mathbf J \circ \mathbf L_\theta ^*( u ) }  -   (\nu_\infty +u) 
\end{equation*}
is well-defined 
and in $\mathrm{Ker} (\mathbf J )$ on this neighbourhood with $u\in\mathrm{Ker}(\mathbf J)$.
Since $\mathbf G(\theta,w^\theta)=0$, we have
\[
(\nu_{\sigma\omega,\infty}+w_{\sigma\omega}^\theta)
  (\widetilde{\mathcal L}_\omega^\theta\varphi)
=
\hat\lambda_\omega^\theta
  (\nu_{\omega,\infty}+w_\omega^\theta)(\varphi)
\]
for every $\varphi\in\mathcal B$ and $m$-a.e.\ $\omega$.
Taking $\varphi=\phi_\omega^\theta$ and using \eqref{eq:1205c},
we obtain
\[
\hat\lambda_\omega^\theta
  (\nu_{\omega,\infty}+w_\omega^\theta)(\phi_\omega^\theta)
=
\lambda_\omega^\theta
  (\nu_{\sigma\omega,\infty}+w_{\sigma\omega}^\theta)
  (\phi_{\sigma\omega}^\theta).
\]
Consequently, by the definition of $\nu_{\omega,\infty}^\theta$,
\begin{align*}
\nu_{\sigma\omega,\infty}^\theta
  (\widetilde{\mathcal L}_\omega^\theta\varphi)
=
\frac{
  \hat\lambda_\omega^\theta
  (\nu_{\omega,\infty}+w_\omega^\theta)(\varphi)
}{
  (\nu_{\sigma\omega,\infty}+w_{\sigma\omega}^\theta)
  (\phi_{\sigma\omega}^\theta)
}=
\lambda_\omega^\theta
\frac{
  (\nu_{\omega,\infty}+w_\omega^\theta)(\varphi)
}{
  (\nu_{\omega,\infty}+w_\omega^\theta)(\phi_\omega^\theta)
}.
\end{align*}
Thus
\begin{equation}\label{eq:1205cb}
\nu_{\sigma\omega,\infty}^\theta
  (\widetilde{\mathcal L}_\omega^\theta\varphi)
=
\lambda_\omega^\theta\nu_{\omega,\infty}^\theta(\varphi),
\qquad
\nu_{\omega,\infty}^\theta(\phi_\omega^\theta)=1
\end{equation}
for every $\varphi\in\mathcal B$ and $m$-a.e.\ $\omega$.
Therefore, by applying the implicit function theorem to the maps $\mathbf F$ and $\mathbf G$, one can investigate the Taylor expansion of $\lambda ^\theta $.

By a computation analogous to that of \cite[Section~3.3]{DFGTV18b}, it can be seen that the derivative of $\phi$ at $\theta=0$ 
is given as
\[
\phi ^\prime (0 )    = - \left(\partial _2\mathbf F (0, 0) \right )^{-1}\left( \partial _1\mathbf F (0, 0)\right),
\]
where $\partial _1\mathbf F$ and $\partial _2\mathbf F$ are the partial derivatives of $\mathbf F$ with respect to $\theta$ and $\varphi$, respectively.
On the other hand, 
\begin{equation*}
 \nu_{\om,\infty  }(\tilde{\mathcal L} _{\sg^{-1}\om}( \varphi\phi _{\sg^{-1}\om}  )) 
= \nu_{\sg^{-1}\om,\infty}( \varphi \phi  _{\sg^{-1}\om} ) =   \eta_{\sg^{-1}\om,\infty  }(\varphi )
\end{equation*}
 \forevery $\varphi\in\mathcal B$ and $m$-a.e.~$\oio$.
So, by the assumption $\eta_{\om,\infty}(f_\om)=0$ $\mas$ and the fact that $\nu_{\om,\infty}(\tilde{\mathcal L}_{\sg^{-1}\om}\phi_{\sg^{-1}\om}) = \nu_{\sg^{-1}\om,\infty}(\phi_{\sg^{-1}\om}) =1 \;\; \mas$ 
 we get
\begin{align*}
 &\partial _1  \mathbf F ( 0 , 0 )(\omega) \\
 &= \frac{  \nu_{\om,\infty }(\tilde{\mathcal L}_{\sg^{-1}\om}  \phi_{\sg^{-1}\om} ) \, \tilde{\mathcal L}_{\sg^{-1}\om}   (f_{\sg^{-1}\om} \phi _{\sg^{-1}\om}) -     \nu_{\om,\infty} (\tilde{\mathcal L}_{\sg^{-1}\om}    (f_{\sg^{-1}\om}  \phi  _{\sg^{-1}\om} )) \, \tilde{\mathcal L}_{\sg^{-1}\om}  \phi _{\sg^{-1}\om}}{\{  \nu_{\om,\infty  }(\tilde{\mathcal L}_{\sg^{-1}\om}    \phi _{\sg^{-1}\om}  ) \}^2}\\
 &=  \tilde{\mathcal L} _{\sg^{-1}\om}   (f _{\sg^{-1}\om}  \phi  _{\sg^{-1}\om} ) 
\end{align*}
for $m$-a.e.~$\oio$.
Similarly, since $\nu_{\om,\infty  } (\tilde{\mathcal L}_{\sg^{-1}\om}    \psi _{\sg^{-1}\om} ) =0 \; \;\mas$ for any $\psi \in \mathrm{Ker}( \mathbf I)$, 
\begin{align*}
& \partial _2  \mathbf F ( 0 , 0 )(\psi )(\om) \\
&= \frac{  \nu_{\om,\infty  }(\tilde{\mathcal L}_{\sg^{-1}\om}  \phi _{\sg^{-1}\om} ) \, \tilde{\mathcal L} _{\sg^{-1}\om}  \psi _{\sg^{-1}\om}   - \nu_{\om,\infty  } (\tilde{\mathcal L}_{\sg^{-1}\om}     \psi _{\sg^{-1}\om}   ) \, \tilde{\mathcal L}_{\sg^{-1}\om}   \phi _{\sg^{-1}\om}   }{\{  \nu_{\om,\infty  }(\tilde{\mathcal L} _{\sg^{-1}\om}  \phi _{\sg^{-1}\om}  ) \}^2} -  \psi _\om\\
& = \tilde{\mathcal L}_{\sg^{-1}\om}    \psi _{\sg^{-1}\om}     -\psi   _\om
\end{align*}
for every $\psi \in \mathrm{Ker}( \mathbf I)$ and $m$-a.e.~$\oio$.
Furthermore, 
by a computation analogous to that of \cite[Section~3.3]{DFGTV18b}, 
one can see that $\partial _2  \mathbf F ( 0 , 0 ): \mathrm{Ker} (\mathbf I) \to   \mathrm{Ker} (\mathbf I) $ is bijective, and 
\[
(\partial _2  \mathbf F ( 0 , 0 ))^{-1}(\varphi )(\om)= -  \sum _{n=0}^\infty \tilde{\mathcal L} ^n_{\sg^{-n}\om} \varphi_{\sg^{-n}\om}    
\]
for $m$-a.e.~$\oio$
(the series converges due to \eqref{LL}; the invertibility of $\partial_2\mathbf G(0,0)$ follows similarly from the uniform gap in \eqref{LL1}).
Hence, it follows from the implicit function theorem that
\begin{align}\label{eq;0822}
\phi ^\prime (0 )(\om)
= 
\sum _{n=1}^\infty   \tilde{\mathcal L} ^n_{\sg^{-n}\om}  (f_{\sg^{-n}\om} \phi_{\sg^{-n}\om} ) 
\end{align}
for $m$-a.e.~$\oio$.
Since 
\[
\nu_{\om,\infty }(\tilde{\mathcal L} ^n_{\sg^{-n}\om}  (f_{\sg^{-n}\om} \phi_{\sg^{-n}\om} ) )
=
 \nu_{\sg^{-n}\om,\infty}(f_{\sg^{-n}\om} \phi_{\sg^{-n}\om} ) 
=  \eta_{\sg^{-n}\om,\infty}(f_{\sg^{-n}\om}  )  =0
\]
by \eqref{eq:1205cb} and the assumption on $f$,
we get 
\begin{equation}\label{eq:0207e3}
\nu_{\om,\infty} (\phi ^\prime (0 )(\om))  =0
\end{equation}
for $m$-a.e.~$\oio$.
By performing a similar calculation, one can also show 
\begin{equation}\label{eq:0207e3b}
\nu_{\om,\infty} (\phi ^{\prime \prime}(0)(\om))=0
\end{equation}
for $m$-a.e.~$\oio$
(cf.~\cite[Lemma 3.15]{DFGTV18b}).

On the other hand, by a computation analogous to that of \cite[Section~3.4]{DFGTV18b}  (recall \eqref{eq:1101} and notice that the $\theta$-derivative of $\theta\mapsto(\om\mapsto\tilde{\mathcal L}_{\sg^{-1}\om}(e^{\theta f_{\sg^{-1}\om}} \varphi _{\sg^{-1}\om}))$  is $\theta\mapsto(\om\mapsto\tilde{\mathcal L}_{\sg^{-1}\om}(f_{\sg^{-1}\om}e^{\theta f_{\sg^{-1}\om}}\varphi _{\sg^{-1}\om}))$),
one can see that, for each $\theta \in U$, 
\begin{align}\label{eq:0207e3a}
\lambda '(\theta )(\om)  =  \nu_{\sg\om,\infty } \left( \tilde{\mathcal L}^\theta _{\om}(f_{\om}\phi^\theta _{\om}) + \tilde{\mathcal L} ^\theta_{\om}( [\phi '(\theta)]_{\om} ) \right)   
\end{align}
for $m$-a.e.~$\oio$.
Hence, by \eqref{eq:1205cb} and \eqref{eq:0207e3},
we conclude
\[
\lambda '(0 )(\om)  =  \eta_{\om,\infty} (f_{\om})
=0 
\quad \text{for $m$-a.e.~$\oio$.}
\]

Furthermore, by considering the $\theta$-derivative of \eqref{eq:0207e3a}, we get
\begin{align*}
\lambda ''(\theta )(\omega ) =   \nu_{\sigma \omega ,\infty}  \left( \tilde{\mathcal L}_{\omega}^\theta (f^2_\om\phi^\theta_\omega ) + 2\tilde{\mathcal L}_{\omega}^\theta( f_\om[\phi '(\theta) ]_\omega) 
+ \tilde{\mathcal L}_{\omega}^\theta( [\phi ''(\theta) ]_\omega) \right).
\end{align*}
By \eqref{eq:1205cb}, \eqref{eq;0822}, and \eqref{eq:0207e3b},
\[
\lambda''(0)(\omega)
=
\nu_{\omega,\infty}(f_\omega^2\phi_\omega)
+
2\sum_{n=1}^\infty
\nu_{\omega,\infty}\!\left(
f_\omega\widetilde{\mathcal L}_{\sigma^{-n}\omega}^{n}
(f_{\sigma^{-n}\omega}\phi_{\sigma^{-n}\omega})
\right).
\]
By the $\sigma$-invariance of $m$ and \eqref{eq:20260705c},
this gives
\[
\int \lambda''(0)\,dm=\Sigma^2(f).
\]
Here the $n$th summand in the series is $a_n(\sigma^{-n}\omega)$,
so termwise integration is justified by the estimates in
Section~\ref{ss:variance}; see \eqref{eq:uniform-ak-bound}.

For item~(2) of Theorem~\ref{conj:1}, 
one can find a neighborhood $U_0$ of $0$ such that the twisted cocycle $(\tilde{\mathcal L}^\theta,\sigma)$ satisfies a Lasota--Yorke estimate with constants uniform in $\theta\in U_0$ (cf.~\cite[Lemma 3.13]{DFGTV18b}).
The proof of \cite[Theorem~3.12]{DFGTV18b} then gives quasicompactness of $(\tilde{\mathcal L}^\theta,\sigma)$ for $\theta\in U_0$, and also shows that the second Lyapunov exponent is upper semicontinuous in $\theta$. Since the untwisted cocycle has a one-dimensional top Oseledets space, there are a neighborhood $U_1$ of $0$ and $\kappa_1\in(0,1)$ such that the second Lyapunov exponent of $(\tilde{\mathcal L}^\theta,\sigma)$ is at most $\log \kappa_1$ for every $\theta\in U_1$.

On the other hand,
$
\tilde{\mathcal L}^\theta_\omega\phi^\theta_\omega=\lambda^\theta_{\omega}\phi^\theta_{\sigma\omega},
$
so $\mathbb C\phi^\theta_\omega$ is an equivariant one-dimensional family. Its Lyapunov exponent is $0$ at $\theta=0$ and depends continuously on $\theta$, so there is a neighborhood $U_2$ of $0$ on which this exponent is still strictly larger than $\log \kappa_1$. Replacing $U$ by $U_0\cap U_1\cap U_2$, we obtain
\[
\mathcal B=\mathbb C\phi^\theta_\omega\oplus \mathrm{Ker}(\nu^\theta_{\omega,\infty})
\]
for $m$-a.e.~$\omega$, and the first two claims of item~(2) follow.

It remains to prove the uniform estimate on $\mathrm{Ker}(\nu^\theta_{\omega,\infty})$. For $\theta\in U$ and $h\in\mathcal B$, set
\[
\Pi^\theta_\omega h:=h-\nu^\theta_{\omega,\infty}(h)\phi^\theta_\omega,
\qquad
Q^\theta_\omega h:=\tilde{\mathcal L}^\theta_\omega\Pi^\theta_\omega h.
\]
Then
$
\Pi^\theta_{\sigma\omega}\tilde{\mathcal L}^\theta_\omega
=
\tilde{\mathcal L}^\theta_\omega\Pi^\theta_\omega,
$
hence
\[
Q^{\theta,n}_\omega h=\tilde{\mathcal L}^{\theta,n}_\omega\Pi^\theta_\omega h
\]
for every $n\ge1$, where $Q^{\theta,n}_\omega:=Q^\theta_{\sigma^{n-1}\omega}\circ\cdots\circ Q^\theta_\omega$.

At $\theta=0$, \eqref{LL} gives
\[
\|Q^{0,n}_\omega h\|
=
\left\|\tilde{\mathcal L}^n_\omega\big(h-\nu_{\omega,\infty}(h)\phi_\omega\big)\right\|
\le
D\kappa^n \left\|h-\nu_{\omega,\infty}(h)\phi_\omega\right\|.
\]
Since $\phi\in L^\infty(\Omega,\mathcal B)$ and $\nu_\infty\in L^\infty(\Omega,\mathcal B^*)$, the constant
\[
M_0:=1+\|\phi\|_{L^\infty(\Omega,\mathcal B)}\|\nu_\infty\|_{L^\infty(\Omega,\mathcal B^*)}
\]
is finite, and therefore
\[
\|Q^{0,n}_\omega\|\le DM_0\kappa^n
\]
for every $n\ge1$ and $m$-a.e.~$\omega$.

Fix $r\in(\kappa,1)$ and choose $n_0\ge1$ so that $DM_0\kappa^{n_0}\le \frac12 r^{n_0}$. The maps $\theta\mapsto \tilde{\mathcal L}^\theta$, $\theta\mapsto \phi^\theta$, and $\theta\mapsto \nu^\theta_\infty$ are continuous into $L^\infty$, hence $\theta\mapsto Q^\theta$ is continuous as a map into $L^\infty(\Omega,L(\mathcal B))$. Therefore there is a neighborhood $U_3$ of $0$ such that
\[
\esssup_{\omega\in\Omega}\|Q^{\theta,n_0}_\omega-Q^{0,n_0}_\omega\|
\le \frac12 r^{n_0}
\]
for every $\theta\in U_3$. Replacing $U$ by $U\cap U_3$, we obtain
\[
\|Q^{\theta,n_0}_\omega\|\le r^{n_0}
\]
for every $\theta\in U$ and $m$-a.e.~$\omega$.

Set
\[
M_1:=\sup_{\theta\in U}\esssup_{\omega\in\Omega}\|Q^\theta_\omega\|<\infty.
\]
If $n=kn_0+j$ with $k\ge0$ and $0\le j<n_0$, then
\[
\|Q^{\theta,n}_\omega\|
\le M_1^j r^{kn_0}
\le C_2 r^n,
\qquad
C_2:=\max_{0\le j<n_0} M_1^j r^{-j}.
\]
If $h\in \mathrm{Ker}(\nu^\theta_{\omega,\infty})$, then $\Pi^\theta_\omega h=h$, hence
\[
\|\tilde{\mathcal L}^{\theta,n}_\omega h\|
=
\|Q^{\theta,n}_\omega h\|
\le C_2 r^n\|h\|.
\]
For general $h\in\mathcal B$,
\[
\tilde{\mathcal L}^{\theta,n}_\omega h-\lambda^{\theta,n}_{\omega}\nu^\theta_{\omega,\infty}(h)\phi^\theta_{\sigma^n\omega}
=
\tilde{\mathcal L}^{\theta,n}_\omega\Pi^\theta_\omega h
=
Q^{\theta,n}_\omega h,
\]
so the last estimate in item~(2) follows as well.

For the second claim of item~(3), we use the fact that $\mathbf F$ is of class~$\mathscr C^3$ (which follows from the analyticity of $\theta \mapsto \tilde{\mathcal L}^\theta$), so the implicit function theorem yields $\mathscr C^3$ regularity of $\theta \mapsto \phi^\theta$, and hence the $\mathscr C^3$ regularity of $\widehat{\Lambda}(\theta)$. We note that real-analyticity of $\widehat{\Lambda}$ was shown in \cite{DH-Limit20} in the closed random case.

For the former claim of item (4),
this is an immediate consequence of the fact that $\phi^\theta$ and $\nu_\infty^\theta$ are continuous with respect to $\theta$ and \eqref{eq:1101}.
The latter claim of item (4) follows from the $\mathscr C^2$ regularity of $\phi^\theta$ and $\nu_\infty^\theta$, Lemma \ref{eq:20240618a}, and \eqref{eq:LL4-uniform-Cphi}:
there is a constant $C>0$, independent of $\om$ and $k$, such that
\[
\left\vert \frac{\nu_{\om,c}(\hat X_{\om,k} \phi^\theta_\om ) - \nu_{\om,c}(\hat X_{\om,k} \phi_\om) }{ \nu_{\om, c}(\hat X_{\om,k} \phi _\om)} \right\vert 
\le
 \frac{\nu_{\om,c}(\hat X_{\om,k} )||\phi^\theta_\om -\phi _\om||_\infty }{ C_{\om,\phi}^{-1}\nu_{\om, c}(\hat X_{\om,k}  )}\le C_{\om,\phi} C_0||\phi^\theta_\om -\phi_\om||\le C|\theta|  
\]
and
\[
\left| \frac{\nu_{\om,\infty}^\theta (\psi _\om) -\nu_{\om,\infty}(\psi _\om) }{ \nu_{ \om,\infty}(\psi  _\om)} \right|
\le 
||\nu_{\om,\infty}^\theta  -\nu_{\om,\infty} ||_{\mathcal B^*}\frac{||\psi _\om||}{\nu_{\om, \infty}(\psi  _\om)} \le C|\theta| \frac{||\psi _\om||}{\nu_{\om, \infty}(\psi  _\om)} 
\]
for $m$-a.e.~$\oio$.

\section{Proof of Theorems \ref{clt}, \ref{ldp}, \ref{be}}\label{s:9}
\label{sec: proof of E, F, G}

Let $\psi$ be the random density of $\zeta\in\mathcal P_{\mathcal B,+}(\nu_c)$.
\subsection{Central limit theorem}\label{ss:CLT}
In this section, we establish Theorem~\ref{clt}.
Now we can calculate the asymptotics of the characteristic function of $\frac{S_nf }{\sqrt n}$ over 
$\zt_{N}$. 
It follows from \eqref{eq:1001a-noomega} that \forevery $t\in\mathbb R$,
\begin{align*}
\zt_{\om,N}(e^{it  \frac{(S_nf)_\om }{\sqrt n}})=
 \frac{\nu_{\sg^n\om,c}(\hat X_{\sg^n\om,N-n} \tilde{\mathcal L}^{\frac{it}{\sqrt n},n}_{\om}  \psi _{\om})}{\nu_{\sg^n\om,c}( \hat X_{\sg^n\om,N-n} \tilde{\mathcal L}^n_{\om}  \psi_{\om} ) } 
 = \frac{L_1'+R_1'}{L_2'+R_2'},
  \end{align*}
  where
\begin{align*}
&L_1'=
\lambda^{\frac{it}{\sqrt n}, n} _{\om}\, \nu_{\om,\infty}^{\frac{it}{\sqrt n}} (\psi _{\om})   \,  \nu_{\sg^n\om,c}\big( \hat X_{\sg^n\om,N-n}\, \phi ^{\frac{it}{\sqrt n}}_{\sg^n\om} \big),\\
&R_1' =\nu_{\sg^n\om,c}\big( \hat X_{\sg^n\om,N-n}\, \big(\tilde{\mathcal L}^{\frac{it}{\sqrt n},n} _\om \psi _\om - \lambda^{\frac{it}{\sqrt n}, n} _{\om}\, \nu_{\om,\infty}^{\frac{it}{\sqrt n}} (\psi _{\om})  \, \phi ^{\frac{it}{\sqrt n}}_{\sg^n\om} \big)\big),\\
 &L_2'=\nu_{\om,\infty} (\psi _\om) \, \nu_{\sg^n\om,c}\big( \hat X_{\sg^n\om,N-n}  \, \phi_{\sg^n\om} \big),\\
&R_2'= \nu_{\sg^n\om,c}\big( \hat X_{\sg^n\om,N-n}\,\big(\tilde{\mathcal L}^{n}_\om   \psi _\om - \nu_{\om,\infty} (\psi _\om) \, \phi_{\sg^n\om}\big) \big) .
 \end{align*}
Let $C_2>0$ and $r\in(0,1)$ be the constants from Theorem~\ref{conj:1}(2). By \eqref{B2} and Theorem~\ref{conj:1}(2),
\begin{align*}
|R_1'|
&\le C_0 \nu_{\sg^n\om,c}(\hat X_{\sg^n\om,N-n})
\left\|
\tilde{\mathcal L}^{\frac{it}{\sqrt n},n}_\om\psi_\om
-
\lambda^{\frac{it}{\sqrt n},n}_{\om}
\nu_{\om,\infty}^{\frac{it}{\sqrt n}}(\psi_\om)\phi^{\frac{it}{\sqrt n}}_{\sg^n\om}
\right\| \\
&\le C_0 C_2\|\psi_\om\| r^n \nu_{\sg^n\om,c}(\hat X_{\sg^n\om,N-n}),
\end{align*}
while \eqref{LL} gives
\[
|R_2'|\le C_0 D\|\psi_\om\|\kappa^n \nu_{\sg^n\om,c}(\hat X_{\sg^n\om,N-n}).
\]
Moreover, Lemma~\ref{eq:20240618a} yields
\[
|L_2'|
=
\nu_{\om,\infty}(\psi_\om)\nu_{\sg^n\om,c}(\hat X_{\sg^n\om,N-n}\phi_{\sg^n\om})
\ge
\frac{\nu_{\om,\infty}(\psi_\om)}{C_{\sg^n\om,\phi}}\nu_{\sg^n\om,c}(\hat X_{\sg^n\om,N-n}).
\]
Hence, for every $n\ge n_{\om,\zt}$,
\[
|L_2'+R_2'|
\ge
\frac12\frac{\nu_{\om,\infty}(\psi_\om)}{C_{\sg^n\om,\phi}}\nu_{\sg^n\om,c}(\hat X_{\sg^n\om,N-n}).
\]

Fix $t\in\mathbb R$. By item~(3) of Theorem~\ref{conj:1}, there are $\delta_\Lambda>0$ and $C_\Lambda>0$ such that
\[
\|\widehat{\Lambda}(\theta)\|_{L^\infty}\le C_\Lambda |\theta|^2
\]
whenever $|\theta|<\delta_\Lambda$. Therefore, whenever $|t|/\sqrt n<\delta_\Lambda$,
\[
\left|\lambda^{\frac{it}{\sqrt n},n}_{\om}\right|
=
\exp\left(\sum_{j=0}^{n-1}\Re \widehat{\Lambda}\left(\frac{it}{\sqrt n}\right)(\sigma^j\omega)\right)
\le e^{C_\Lambda t^2}.
\]
Also, Theorem~\ref{conj:1}(4) gives
\[
\left|
\frac{\nu_{\sg^n\om,c}(\hat X_{\sg^n\om,N-n}\phi^{\frac{it}{\sqrt n}}_{\sg^n\om})}{\nu_{\sg^n\om,c}(\hat X_{\sg^n\om,N-n}\phi_{\sg^n\om})}-1
\right|
\le C\frac{|t|}{\sqrt n}
\]
and
\[
\left|
\frac{\nu_{\om,\infty}^{\frac{it}{\sqrt n}}(\psi_\om)}{\nu_{\om,\infty}(\psi_\om)}-1
\right|
\le
C \frac{\|\psi_\om\|}{\nu_{\om,\infty}(\psi_\om)}\frac{|t|}{\sqrt n}
\]
for $m$-a.e.~$\oio$, every $N\ge n\ge1$, and every $t\in\mathbb R$ such that $\frac{it}{\sqrt n}\in U$.
Thus, for each fixed $t$ and $m$-a.e.~$\oio$ there exists $n_1(\omega,t)$ such that
\[
\frac{|L_1'|}{|L_2'|}\le 4e^{C_\Lambda t^2}
\]
for every $n\ge n_1(\omega,t)$ and every $N\ge n$.

Using
\[
\frac{L_1'+R_1'}{L_2'+R_2'}-\frac{L_1'}{L_2'}
=
\frac{R_1'}{L_2'+R_2'}-\frac{R_2'}{L_2'+R_2'}\cdot \frac{L_1'}{L_2'},
\]
we obtain
\begin{align*}
\left|
\zt_{\om,N}\left(e^{it  \frac{(S_nf)_\om }{\sqrt n}}\right)-\frac{L_1'}{L_2'}
\right|
&\le  \frac{|R_1'|}{|L_2' + R_2'|} + \frac{|R_2'|}{|L_2'+ R_2'|} \cdot \frac{|L_1'|}{|L_2'|} \\
&\le
\frac{2C_0C_{\sg^n\om,\phi}\|\psi_\om\|}{\nu_{\om,\infty}(\psi_\om)}
\left(C_2 r^n+4e^{C_\Lambda t^2}D\kappa^n\right)
\end{align*}
for every $n\ge \max\{n_{\om,\zt},n_1(\omega,t)\}$ and every $N\ge n$. Since $C_\phi$ is essentially bounded by \eqref{eq:LL4-uniform-Cphi} and both $r$ and $\kappa$ belong to $(0,1)$, the right-hand side tends to $0$ as $n\to\infty$.
Therefore,
\[
\lim_{n\to\infty}\sup_{N\ge n}\left|
\zt_{\om,N}\left(e^{it  \frac{(S_nf)_\om }{\sqrt n}}\right)-
\lambda^{\frac{it}{\sqrt n}, n} _{\om}
 \frac{\nu_{\sg^n\om,c}(\hat X_{\sg^n\om,N-n}\phi^{\frac{it}{\sqrt n}}_{\sg^n\om}) \nu_{\om,\infty}^{\frac{it}{\sqrt n}}(\psi _\om)}{\nu_{\sg^n\om,c}(\hat X_{\sg^n\om,N-n}\phi_{\sg^n\om}) \nu_{\om,\infty}(\psi _\om)}
\right|=0
\]
for $m$-a.e.~$\oio$.

Theorem~\ref{conj:1}(4) also implies that
\[
\lim_{n\to\infty}\sup_{N\ge n}\left|
\frac{\nu_{\sg^n\om,c}(\hat X_{\sg^n\om,N-n}\phi^{\frac{it}{\sqrt n}}_{\sg^n\om}) \nu_{\om,\infty}^{\frac{it}{\sqrt n}}(\psi _\om)}{\nu_{\sg^n\om,c}(\hat X_{\sg^n\om,N-n}\phi_{\sg^n\om}) \nu_{\om,\infty}(\psi _\om)}-1
\right|=0
\]
for $m$-a.e.~$\oio$. Consequently,
\[
\lim_{n\to\infty}\sup_{N\ge n}\left|
\zt_{\omega ,N}\left(e^{it  \frac{(S_nf)_\omega }{\sqrt n}}\right)-
\lambda^{\frac{it}{\sqrt n}, n} _{\omega}
\right|=0
\]
for $m$-a.e.~$\oio$.

Since
\[
\lambda^{\frac{it}{\sqrt n}, n} _{\omega}
=
\exp\left(\sum_{j=0}^{n-1}  \widehat \Lambda \left(\frac{it}{\sqrt n}\right)(\sigma ^j\omega) \right),
\]
item~(3) of Theorem~\ref{conj:1} gives
\[
\sum_{j=0}^{n-1}  \widehat \Lambda \left(\frac{it}{\sqrt n}\right)(\sigma ^j\omega)
=
-\frac{t^2}{2n}\sum_{j=0}^{n-1}\widehat \Lambda ''(0)(\sigma ^j\omega)+o(1)
\]
for $m$-a.e.~$\oio$. Birkhoff's ergodic theorem and Theorem~\ref{conj:1}(3) then yield
\[
\lambda^{\frac{it}{\sqrt n}, n} _{\omega}\to e^{-\frac{t^2\Sigma ^2(f)}{2}}
\]
for $m$-a.e.~$\oio$. Therefore, for every $t\in\mathbb R$ and $m$-a.e.~$\oio$,
\[
\sup_{N\ge n}\left|
\zt_{\omega ,N}\left(e^{it  \frac{(S_nf)_\omega }{\sqrt n}}\right)-e^{-\frac{t^2\Sigma ^2(f)}{2}}
\right|\to0
\qquad (n\to\infty).
\]
Now fix $a\in\mathbb R$ and let $(N_n)_{n\ge1}$ be any sequence with $N_n\ge n$. The characteristic functions of the laws of $\frac{(S_nf)_\omega}{\sqrt n}$ under $\zeta_{\omega,N_n}$ converge pointwise to $t\mapsto e^{-t^2\Sigma^2(f)/2}$. Since the Gaussian distribution is continuous, L\'evy's continuity theorem gives
\[
\zeta_{\omega,N_n}\left(\frac{(S_nf)_\omega}{\sqrt n}\le a\right)\to
\frac{1}{\sqrt{2\pi \Sigma ^2(f)}} \int ^a _{-\infty} \exp\left(-\frac{z^2}{2\Sigma ^2(f)}\right) dz.
\]
As this holds for every sequence $N_n\ge n$, the stated $\sup_{N\ge n}$ convergence follows. 

\subsection{Large deviation principle}

In this section we show Theorem~\ref{ldp}. 
We first recall the G\"artner--Ellis theorem (cf.~\cite[lemma XIII.2]{HH01}):
\begin{theorem}[G\"artner--Ellis]
Let $S_n$ be a random variable on a probability space $(Z,\mathscr  G,\xi_n)$ for each $n\ge 0$.
Assume that there is $\theta_0>0$ such that \forevery $\theta\in[-\theta_0,\theta_0]$,
\[
\lim _{n\to\infty} \frac{1}{n} \log \xi_n(e^{\theta S_n} ) = \Lambda (\theta ), 
\]
where $\Lambda $ is a strictly convex $\mathcal C^1$ function satisfying $\frac{d}{d\theta}\Lambda (0)=0$.
Then there exists $\delta_0>0$ such that,
for every $\delta\in[0,\delta_0)$,
\[
\lim_{n\to\infty}\frac{1}{n}
\log\xi_n\bigl(\{x\in Z:S_n(x)>n\delta\}\bigr)
=-c(\delta),
\]
where
\[
c(\delta)=
\sup_{\theta\in[-\theta_0,\theta_0]}
\{\theta\delta-\Lambda(\theta)\}.
\]
The function $c$ is nonnegative, continuous,
strictly convex on a sufficiently small
neighborhood of $0$, and satisfies $c(0)=0$.
\end{theorem}

Let $C$ be the constant in Theorem~\ref{conj:1}(4), and let
$C_2,r$ be as in item~(2) of that theorem. Choose $\theta_0>0$
so small that $[-\theta_0,\theta_0]\subset U$ and
\[
 C\theta_0\le\frac12,\qquad
 r\exp\!\left(\theta_0C_0\esssup_{\oio}\|f_\om\|\right)<1,
 \qquad
 \inf_{|\theta|\le\theta_0}\essinf_\omega\lambda_\omega^\theta>r.
\]
For $N\ge n$, use the same four terms as in Section~\ref{ss:CLT}, now with
$\theta$ in place of $it/\sqrt n$:
\begin{align*}
L_1'&=\lambda_\omega^{\theta,n}\nu_{\omega,\infty}^\theta(\psi_\omega)
 \nu_{\sigma^n\omega,c}
       (\hat X_{\sigma^n\omega,N-n}\phi_{\sigma^n\omega}^\theta),\\
R_1'&=\nu_{\sigma^n\omega,c}\!\left(\hat X_{\sigma^n\omega,N-n}
 [\widetilde{\mathcal L}_\omega^{\theta,n}\psi_\omega
   -\lambda_\omega^{\theta,n}\nu_{\omega,\infty}^\theta(\psi_\omega)
         \phi_{\sigma^n\omega}^\theta]\right),\\
L_2'&=\nu_{\omega,\infty}(\psi_\omega)
 \nu_{\sigma^n\omega,c}(\hat X_{\sigma^n\omega,N-n}\phi_{\sigma^n\omega}),\\
R_2'&=\nu_{\sigma^n\omega,c}\!\left(\hat X_{\sigma^n\omega,N-n}
 [\widetilde{\mathcal L}_\omega^n\psi_\omega
   -\nu_{\omega,\infty}(\psi_\omega)\phi_{\sigma^n\omega}]\right).
\end{align*}
Then
\begin{equation*}
 \zeta_{\omega,N}(e^{\theta(S_nf)_\omega})
 =\frac{L_1'+R_1'}{L_2'+R_2'}.
\end{equation*}
As in Section~\ref{ss:CLT}, for $n\ge n_{\omega,\zeta}$,
\begin{align*}
 &|R_1'|\le C_0C_2\|\psi_\omega\|r^n
               \nu_{\sigma^n\omega,c}(\hat X_{\sigma^n\omega,N-n}),\qquad
 |R_2'|\le C_0D\|\psi_\omega\|\kappa^n
               \nu_{\sigma^n\omega,c}(\hat X_{\sigma^n\omega,N-n}),\\
 &L_2'+R_2'\ge
 \frac{\nu_{\omega,\infty}(\psi_\omega)}{2C_{\sigma^n\omega,\phi}}
               \nu_{\sigma^n\omega,c}(\hat X_{\sigma^n\omega,N-n}).
\end{align*}
Also, Theorem~\ref{conj:1}(4) gives
\begin{equation}\label{eq:ldp-terminal-ratio}
 \frac12\le
 \frac{\nu_{\sigma^n\omega,c}
       (\hat X_{\sigma^n\omega,N-n}\phi_{\sigma^n\omega}^\theta)}
      {\nu_{\sigma^n\omega,c}
       (\hat X_{\sigma^n\omega,N-n}\phi_{\sigma^n\omega})}
 \le\frac32.
\end{equation}

We next verify that $\nu_{\omega,\infty}^\theta(\psi_\omega)>0$ for $\theta\in[-\theta_0,\theta_0]$.
Otherwise, continuity in $\theta$ and positivity at $\theta=0$ give some $\theta\in[-\theta_0,\theta_0]$ with $\nu_{\omega,\infty}^\theta(\psi_\omega)=0$, hence $L_1'=0$.
Taking $N=n$, we obtain
\begin{align*}
 \exp\!\left(-n|\theta|C_0\esssup_{\oio}\|f_\om\|\right)
 \le\zeta_{\omega,n}(e^{\theta(S_nf)_\omega})
 \le
 \frac{2C_0C_2\,\esssup_\omega C_{\omega,\phi}\,\|\psi_\omega\|}
      {\nu_{\omega,\infty}(\psi_\omega)}r^n.
\end{align*}
The first inequality uses only the boundedness of the real observable.
These bounds contradict the choice of $\theta_0$ as $n\to\infty$, due to \eqref{eq:LL4-uniform-Cphi}.
Thus $\nu_{\omega,\infty}^\theta(\psi_\omega)>0$
throughout $[-\theta_0,\theta_0]$.

Consequently $L_1',L_2'>0$, and
\begin{align*}
 \left|\frac{R_1'}{L_1'}\right|
 &\le
 \frac{2C_0C_2\,\esssup_\omega C_{\omega,\phi}\,\|\psi_\omega\|}
      {\nu_{\omega,\infty}^\theta(\psi_\omega)}
 \left(\frac{r}{\essinf_\omega\lambda_\omega^\theta}\right)^n,\\
 \left|\frac{R_2'}{L_2'}\right|
 &\le
 \frac{C_0D\,\esssup_\omega C_{\omega,\phi}\,\|\psi_\omega\|}
      {\nu_{\omega,\infty}(\psi_\omega)}\kappa^n.
\end{align*}
For each fixed $\omega,\psi_\omega,\theta$, both bounds tend to zero
uniformly for $N\ge n$. For all sufficiently large $n$, the logarithms
below are therefore real and
\begin{align*}
&\frac1n\log\zeta_{\omega,N}(e^{\theta(S_nf)_\omega})
       -\frac1n\sum_{j=0}^{n-1}\log\lambda_{\sigma^j\omega}^\theta\\
&\quad=\frac1n\log
       \frac{\nu_{\omega,\infty}^\theta(\psi_\omega)}
            {\nu_{\omega,\infty}(\psi_\omega)}
 +\frac1n\log
       \frac{\nu_{\sigma^n\omega,c}
          (\hat X_{\sigma^n\omega,N-n}\phi_{\sigma^n\omega}^\theta)}
            {\nu_{\sigma^n\omega,c}
          (\hat X_{\sigma^n\omega,N-n}\phi_{\sigma^n\omega})}
          +\frac1n\log\left(1+\frac{R_1'}{L_1'}\right)
       -\frac1n\log\left(1+\frac{R_2'}{L_2'}\right).
\end{align*}
The first logarithm on the right-hand side is independent of $n$ and $N$, the second is bounded by \eqref{eq:ldp-terminal-ratio}, and the last two terms tend to zero by the relative remainder bounds.
Thus the left-hand side tends to zero as $n\to\infty$, uniformly in $N\ge n$.
Birkhoff's theorem gives, for every sequence $N_n\ge n$,
\[
 \lim_{n\to\infty}\frac1n
 \log\zeta_{\omega,N_n}(e^{\theta(S_nf)_\omega})
 =\int_\Omega\widehat\Lambda(\theta)\,dm
 =\frac{\Sigma^2(f)}2\theta^2+o(\theta^2).
\]
The full-measure set on which the first equality holds may depend on $\theta$.
However, a common full-measure set for all $\theta\in[-\theta_0,\theta_0]$ can be obtained by first considering a countable dense subset of this interval and then extending the equality to the whole interval:
the functions $n^{-1}\log\zeta_{\omega,N}(e^{\theta(S_nf)_\omega})$ are Lipschitz in $\theta$ with Lipschitz constant $C_0\esssup_{\oio}\|f_\om\|$, uniformly in $N\ge n$.

The function $\theta\mapsto\int\widehat\Lambda(\theta)\,dm$ is $C^3$,
centered, and strictly convex on a smaller real neighbourhood because
its second derivative at $0$ is $\Sigma^2(f)>0$.
Apply the G\"artner--Ellis theorem with $Z=X$,
$\xi_n=\zeta_{\omega,N_n}$, and $S_n=(S_nf)_\omega$.
The resulting rate function is independent of $\omega$ and the
initial density. Since the argument applies to every sequence
$N_n\ge n$, the convergence is uniform in $N\ge n$ as in
Theorem~\ref{ldp}.

\subsection{Berry--Esseen theorem}\label{ss:BE}

We follow the proof of \cite[Section~4.4]{DH-Limit20} in the closed case.

Fix $\zt\in\cP_{\cB,+}(\nu_c)$, and let $\psi$ be its random density.
Recall $\sg_{\om,n}^2=\sg_{\om,n}^2(f)$ from \eqref{eq:0805b}.
Let $U$ and $(\lambda^\theta,\phi^\theta,\nu_\infty^\theta)$ be as in Theorem~\ref{conj:1}, and let $C_*>0$ and $r\in(0,1)$ be the constants in Theorem~\ref{conj:1}(2).
For $\delta>0$, write
\[
M(\delta):=\sup_{|\vartheta|\le \delta}\|\widehat\Lambda'''(\vartheta)\|_{L^\infty(m)}.
\]
Choose $\dl_0>0$ so that $\{\theta\in\mathbb C:|\theta|\le \dl_0\}\subset U$ and
\[
\frac{\dl_0\,M(\dl_0)}{3\Sigma^2(f)}\le \frac{1}{4}.
\]
Such a choice is possible because $\widehat\Lambda$ is $\mathscr C^3$ on $U$.
Shrinking $\delta_0$ if necessary, also require that $\{|z|\le2\delta_0\}\subset U$ and put
\[
q:=r\exp\Big(\sup_{|z|\le2\delta_0}\|\widehat\Lambda(z)\|_{L^\infty(m)}\Big)<1.
\]
This is possible since $\widehat\Lambda(0)=0$.

For $n\ge 1$ and $N\ge n$, set
\[
F_{\om,n,N}(a):=\zt_{\om,N}\!\left(\frac{S_nf_\om}{\sg_{\om,n}}\le a\right),
\qquad
\Psi_{\om,n,N}(t):=\zt_{\om,N}\!\left(e^{itS_nf_\om/\sg_{\om,n}}\right),
\]
\[
\Phi(a):=\frac{1}{\sqrt{2\pi}}\int_{-\infty}^a e^{-z^2/2}\,dz,
\qquad
T_{\om,n}:=\dl_0\sg_{\om,n}.
\]
By Esseen's smoothing inequality (see, for example, \cite[Section~4.4]{DH-Limit20}), for $m$-a.e.~$\oio$ we have
\begin{equation}\label{eq:20260409}
\sup_{a\in\mathbb R}\bigl|F_{\om,n,N}(a)-\Phi(a)\bigr|
\le \frac{1}{\pi}\int_{-T_{\om,n}}^{T_{\om,n}}
\left|\frac{\Psi_{\om,n,N}(t)-e^{-t^2/2}}{t}\right|\,dt
+\frac{24}{\pi\sqrt{2\pi}\,T_{\om,n}}.
\end{equation}

Since $\Sigma^2(f)>0$ and $\sg_{\om,n}^2/n\to\Sigma^2(f)$ $m$-a.s., define
\[
n_1(\om):=\inf\left\{n\ge 1:
\sigma_{\omega,k}^2=\sum_{j=0}^{k-1}\widehat\Lambda''(0)(\sigma^j\omega),\,
\frac{\Sigma^2(f)}{2}\le \frac{\sg_{\om,k}^2}{k}\le 2\Sigma^2(f)
\text{ for all }k\ge n\right\}.
\]
Then $n_1(\om)<\infty$ for $m$-a.e.~$\om$, and for every $n\ge n_1(\om)$,
\[
\dl_0\sqrt{\frac{\Sigma^2(f)}{2}}\,\sqrt n
\le T_{\om,n}\le
\dl_0\sqrt{2\Sigma^2(f)}\,\sqrt n.
\]
Hence
\[
\frac{24}{\pi\sqrt{2\pi}\,T_{\om,n}}
\le
\frac{24}{\pi\sqrt{\pi}\,\dl_0\sqrt{\Sigma^2(f)}}\,\frac{1}{\sqrt n}
\qquad\text{for every }n\ge n_1(\om).
\]
Let
\[
K_0:=\frac{24}{\pi\sqrt{\pi}\,\dl_0\sqrt{\Sigma^2(f)}}.
\]

Fix $\om$ in a full-measure set on which the conclusions of Theorem~\ref{conj:1}, Lemma~\ref{eq:20240618a}, and the convergence $\sg_{\om,n}^2/n\to\Sigma^2(f)$ all hold.
For $|t|\le T_{\om,n}$, set
\[
\theta(t)=\theta_{\om,n}(t):=\frac{it}{\sg_{\om,n}}.
\]
Then $|\theta(t)|\le \dl_0$.

Write
\[
\Lambda_{\om,n}(\theta):=\sum_{j=0}^{n-1}\widehat\Lambda(\theta)(\sg^j\om),
\]
\[
A_{\om,n,N}(t):=
\frac{\nu_{\sg^n\om,c}(\hat X_{\sg^n\om,N-n}\phi^{\theta(t)}_{\sg^n\om})}
     {\nu_{\sg^n\om,c}(\hat X_{\sg^n\om,N-n}\phi_{\sg^n\om})},
\qquad
B_{\om,n}(t):=
\frac{\nu^{\theta(t)}_{\om,\infty}(\psi_\om)}{\nu_{\om,\infty}(\psi_\om)}.
\]
We also set
\[
E_{\om,n,N}(t):=
\Psi_{\om,n,N}(t)-e^{\Lambda_{\om,n}(\theta(t))}A_{\om,n,N}(t)B_{\om,n}(t).
\]

By Theorem~\ref{conj:1}(4), there are a constant $K_A$ and a random variable $K_B(\om)>0$ such that
\[
|A_{\om,n,N}(t)-1|\le K_A|\theta(t)|,
\qquad
|B_{\om,n}(t)-1|\le K_B(\om)|\theta(t)|
\]
for every $n\ge 1$, $N\ge n$, and $|t|\le T_{\om,n}$.
Set
\[
M_A:=1+\dl_0K_A,
\qquad
M_B(\om):=1+\dl_0K_B(\om).
\]
Then
\[
|A_{\om,n,N}(t)|\le M_A,
\qquad
|B_{\om,n}(t)|\le M_B(\om).
\]

To estimate $E_{\om,n,N}(t)$, write
\[
R_{\om,n}^\theta:=
\tilde{\mathcal L}_{\om}^{\theta,n}\psi_\om
-
\lambda_{\om}^{\theta,n}\nu_{\om,\infty}^\theta(\psi_\om)\phi_{\sg^n\om}^\theta.
\]
By Theorem~\ref{conj:1}(2),
\[
\|R_{\om,n}^\theta\|\le C_*r^n\|\psi_\om\|
\qquad (|\theta|\le \dl_0,\, n\ge 0, \, \text{$m$-a.e.~$\oio$}).
\]
Moreover, $R_{\om,n}^0$ is obtained by setting $\theta=0$, and \eqref{eq:1001a-noomega} gives
\begin{multline*}
E_{\om,n,N}(t)
=
\frac{\nu_{\sg^n\om,c}\!\left(\hat X_{\sg^n\om,N-n}(R_{\om,n}^{\theta(t)}-R_{\om,n}^0)\right)}
{\nu_{\sg^n\om,c}\!\left(\hat X_{\sg^n\om,N-n}\tilde{\mathcal L}_{\om}^{n}\psi_\om\right)}\\
-
\bigl(e^{\Lambda_{\om,n}(\theta(t))}A_{\om,n,N}(t)B_{\om,n}(t)-1\bigr)
\frac{\nu_{\sg^n\om,c}\!\left(\hat X_{\sg^n\om,N-n}R_{\om,n}^0\right)}
{\nu_{\sg^n\om,c}\!\left(\hat X_{\sg^n\om,N-n}\tilde{\mathcal L}_{\om}^{n}\psi_\om\right)}.
\end{multline*}
Regard the expression on the right-hand side above as a function
$\mathcal E_{\omega,n,N}(z)$ of the complex parameter $z$.
The primal and dual fixed-point constructions in
Subsection~\ref{ss:implicit} are holomorphic: the twisted
operators depend holomorphically on $z$ by \eqref{B2},
and the normalizing denominators are bounded away from zero
on a neighbourhood of the unperturbed data.
Thus $\mathcal E_{\omega,n,N}$ is holomorphic on
$\{|z|\le2\delta_0\}$ and
$\mathcal E_{\omega,n,N}(0)=0$.

The denominator estimate \eqref{eq:20260215c} and
Lemma~\ref{eq:20240618a}, together with \eqref{LL4},
show that the ratio of
$\nu_{\sigma^n\omega,c}(\widehat X_{\sigma^n\omega,N-n}R)$
to the denominator in the preceding display is bounded by a
random constant times $\|R\|$, uniformly for
$N\ge n\ge n_{\omega,\zeta}$.
Theorem~\ref{conj:1}(2), the uniform bounds in
Theorem~\ref{conj:1}(4), and the definition of $q$ therefore
give a finite random variable $K_*(\omega)$, independent of
$n$ and $N$, such that
\[
\sup_{|z|\le2\delta_0}
|\mathcal E_{\omega,n,N}(z)|
\le K_*(\omega)q^n
\qquad(N\ge n\ge n_{\omega,\zeta}).
\]
By Cauchy's derivative estimate on the concentric discs,
\[
\sup_{|z|\le\delta_0}
|\mathcal E'_{\omega,n,N}(z)|
\le\frac{K_*(\omega)}{\delta_0}q^n.
\]
Since $\mathcal E_{\omega,n,N}(0)=0$, it follows that
\[
|E_{\omega,n,N}(t)|
\le K_E(\omega)q^n|\theta(t)|,
\qquad K_E(\omega):=K_*(\omega)/\delta_0,
\]
for $N\ge n\ge n_{\omega,\zeta}$ and
$|t|\le T_{\omega,n}$.

We now estimate the leading exponential factor.
Taylor's theorem gives
\[
\widehat\Lambda(\theta(t))(\sg^j\om)
=
\frac{\theta(t)^2}{2}\widehat\Lambda''(0)(\sg^j\om)
+
\frac{\theta(t)^3}{2}\int_0^1(1-s)^2\widehat\Lambda'''(s\theta(t))(\sg^j\om)\,ds.
\]
Summing over $j=0,\dots,n-1$ and using \eqref{eq:0805b}, we obtain
\[
\Lambda_{\om,n}(\theta(t))=-\frac{t^2}{2}+W_{\om,n}(t),
\]
where
\[
W_{\om,n}(t):=
\frac{\theta(t)^3}{2}
\sum_{j=0}^{n-1}\int_0^1(1-s)^2\widehat\Lambda'''(s\theta(t))(\sg^j\om)\,ds.
\]

To estimate the first term of the right-hand side of \eqref{eq:20260409}, we now write
\[
\Psi_{\om,n,N}(t)-e^{-t^2/2}
=e^{-t^2/2}\bigl(e^{W_{\om,n}(t)}A_{\om,n,N}(t)B_{\om,n}(t)-1\bigr)+E_{\om,n,N}(t),
\]
and
\begin{align*}
&e^{W_{\om,n}(t)}A_{\om,n,N}(t)B_{\om,n}(t)-1\\
&=
\bigl(e^{W_{\om,n}(t)}-1\bigr)A_{\om,n,N}(t)B_{\om,n}(t)
+\bigl(A_{\om,n,N}(t)-1\bigr)B_{\om,n}(t)
+\bigl(B_{\om,n}(t)-1\bigr).
\end{align*}

Since $\int_0^1(1-s)^2\,ds=1/3$, for $n\ge n_1(\om)$ we have
\[
|W_{\om,n}(t)|
\le
\frac{M(\dl_0)}{6}\,n\,|\theta(t)|^3
\le
\frac{\sqrt{2}M(\dl_0)}{3\Sigma^2(f)^{3/2}}\frac{|t|^3}{\sqrt n}.
\]
Set
\[
C_W:=\frac{\sqrt{2}M(\dl_0)}{3\Sigma^2(f)^{3/2}}.
\]
Using also $|t|\le T_{\om,n}=\dl_0\sg_{\om,n}$ and $\sg_{\om,n}^2\ge \Sigma^2(f)n/2$, we get
\[
|W_{\om,n}(t)|
\le
\frac{M(\dl_0)\dl_0}{3\Sigma^2(f)}\,t^2
\le \frac{t^2}{4}.
\]
Therefore, since $|e^z-1|\le |z|e^{|z|}$,
\[
e^{-t^2/2}|e^{W_{\om,n}(t)}-1|
\le
e^{-t^2/4}|W_{\om,n}(t)|
\le
C_W\frac{|t|^3}{\sqrt n}e^{-t^2/4}.
\]

Set
\[
I_0:=\int_{\mathbb R} e^{-t^2/2}\,dt,
\qquad
I_2:=\int_{\mathbb R} t^2e^{-t^2/4}\,dt.
\]
Then,
we obtain, for every $n\ge n_1(\om)$,
\begin{align*}
&\int_{-T_{\om,n}}^{T_{\om,n}}
\frac{e^{-t^2/2}\bigl(|A_{\om,n,N}(t)-1|\,|B_{\om,n}(t)|+|B_{\om,n}(t)-1|\bigr)}{|t|}\,dt \\
&\le
\frac{\sqrt{2}}{\sqrt{\Sigma^2(f)}}\,
\frac{K_AM_B(\om)+K_B(\om)}{\sqrt n}\,I_0,
\end{align*}
\[
\int_{-T_{\om,n}}^{T_{\om,n}}
\frac{e^{-t^2/2}|e^{W_{\om,n}(t)}-1|\,|A_{\om,n,N}(t)|\,|B_{\om,n}(t)|}{|t|}\,dt
\le
\frac{C_WM_AM_B(\om)}{\sqrt n}\,I_2,
\]
and
\[
\int_{-T_{\om,n}}^{T_{\om,n}}
\frac{|E_{\om,n,N}(t)|}{|t|}\,dt
\le 2\dl_0 K_E(\om)q^n.
\]

Let
\[
n_2:=\min\{n\ge1:q^k\le k^{-1/2}\text{ for every }k\ge n\},
\qquad
n_0(\om):=\max\{n_1(\om),n_2,n_{\om,\zeta}\},
\]
and define
\[
K_{\Psi}(\om):=
\frac{\sqrt{2}}{\sqrt{\Sigma^2(f)}}\bigl(K_AM_B(\om)+K_B(\om)\bigr)I_0
+
C_WM_AM_B(\om)I_2
+
2\dl_0K_E(\om).
\]
Then, for every $n\ge n_0(\om)$,
\[
\int_{-T_{\om,n}}^{T_{\om,n}}
\left|\frac{\Psi_{\om,n,N}(t)-e^{-t^2/2}}{t}\right|\,dt
\le
\frac{K_{\Psi}(\om)}{\sqrt n}.
\]
Substituting this into Esseen's inequality \eqref{eq:20260409} yields
\[
\sup_{a\in\mathbb R,\ N\ge n}\bigl|F_{\om,n,N}(a)-\Phi(a)\bigr|
\le
\frac{1}{\pi}\frac{K_{\Psi}(\om)}{\sqrt n}
+
\frac{K_0}{\sqrt n}
\qquad\text{for every }n\ge n_0(\om).
\]

Finally, set
\[
C(\om):=
\sqrt{n_0(\om)}
+
\frac{K_{\Psi}(\om)}{\pi}
+
K_0.
\]
If $1\le n<n_0(\om)$, then
\[
\sup_{a\in\mathbb R,\ N\ge n}\bigl|F_{\om,n,N}(a)-\Phi(a)\bigr|
\le 1
\le \frac{\sqrt{n_0(\om)}}{\sqrt n}
\le \frac{C(\om)}{\sqrt n}.
\]
If $n\ge n_0(\om)$, the previous inequality gives the same bound.
Hence
\[
\sup_{a\in\mathbb R,\ N\ge n}\bigl|F_{\om,n,N}(a)-\Phi(a)\bigr|
\le \frac{C(\om)}{\sqrt n}
\]
for every $n\ge 1$ and for $m$-a.e.~$\oio$.

\appendix

\section{Proof of claims in Remark \ref{rmk:1229}}\label{app:B}

\subsection{Equivalence of (\ref{eq:20240614}) and (\ref{eq:20240626})}
For a bounded measurable random variable $f$, write
\[
A_{n} f(\om,x)=A_{\om,n} f(x):=\frac{1}{n}\sum_{j=0}^{n-1}f_{\sigma^j\omega}\circ T_\omega^j(x)
\qquad \text{for }(\omega,x)\in X,\ n\ge 1.
\]
By the definition of the pushforward random measure,
\begin{equation}\label{eq:appB}
\frac{1}{n}\sum_{j=0}^{n-1}T_*^j(\zeta)_{\sigma^j\omega}(f_{\sigma^j\omega})
=
\zeta_\omega(A_{\om,n} f)
\qquad \text{for $m$-a.e.~$\oio$.}
\end{equation}

Assume first that \eqref{eq:20240614} holds.
Let $f$ be bounded and let $\zeta$ be a random probability measure such that $\zeta_\omega$ is absolutely continuous with respect to $\varrho_\omega$ for $m$-a.e.~$\omega$.
Then $A_n f\to \int f\,d\varrho$ $\varrho$-a.e. by \eqref{eq:20240614}.
Since $\zeta_\omega$ is absolutely continuous with respect to $\varrho_\omega$ and $|A_n f|\leq \|f\|_\infty$, dominated convergence together with \eqref{eq:appB} yields
\[
\lim_{n\to\infty}\frac{1}{n}\sum_{j=0}^{n-1}T_*^j(\zeta)_{\sigma^j\omega}(f_{\sigma^j\omega})
=
\int_\Omega\varrho_\omega(f_\omega)\,dm(\omega)
\qquad \text{for $m$-a.e.~$\oio$.}
\]
Taking $\zeta=\varrho$ and using $T_*(\varrho)=\varrho$ (or, simply by Birkhoff's ergodic theorem), we obtain
\[
\lim_{n\to\infty}\frac{1}{n}\sum_{j=0}^{n-1}\varrho_{\sigma^j\omega}(f_{\sigma^j\omega})
=
\int_\Omega\varrho_\omega(f_\omega)\,dm(\omega)
\qquad \text{for $m$-a.e.~$\oio$.}
\]
Hence \eqref{eq:20240626} follows.

Conversely, assume \eqref{eq:20240626}.
Let $f$ be bounded.
Since $\varrho$ is $T$-invariant, Birkhoff's ergodic theorem gives a $T$-invariant random variable $h_f\in L^1(\varrho)$ such that
\[
A_nf \to h_f \qquad \varrho\text{-a.e.}
\]
Because $|A_nf|\leq \|f\|_\infty$, dominated convergence and \eqref{eq:appB} imply that for every $\zeta$ such that $\zeta_\omega$ is absolutely continuous with respect to $\varrho_\omega$ for $m$-a.e.~$\omega$,
\[
\zeta_\omega(h_{\omega,f})
=
\lim_{n\to\infty}\zeta_\omega(A_{\om,n}f)
=
\int_\Omega\varrho_\omega(f_\omega)\,dm(\omega)
\qquad \text{for $m$-a.e.~$\oio$.}
\]
Fix $\varepsilon>0$ and consider
\[
A^+:=\Big\{(\omega,x)\in X : h_{\omega,f}(x)>\int_\Omega\varrho_\omega(f_\omega)\,dm(\omega)+\varepsilon\Big\}.
\]
If $\varrho(A^+)>0$, then the set
\[
B^+:=\{\omega\in\Omega : \varrho_\omega(A^+_\omega)>0\}
\]
has positive $m$-measure.
Define a random probability measure $\zeta^+$ by
\[
\zeta^+_\omega :=
\begin{cases}
\dfrac{\varrho_\omega(\,\cdot\,\cap A^+_\omega)}{\varrho_\omega(A^+_\omega)} & \text{if }\omega\in B^+,\\[1.2ex]
\varrho_\omega & \text{if }\omega\notin B^+.
\end{cases}
\]
Then $\zeta^+_\omega$ is absolutely continuous with respect to $\varrho_\omega$ for $m$-a.e.~$\omega$, while 
\[
\zeta^+_\omega(h_{f,\omega})\ge \int_\Omega\varrho_\omega(f_\omega)\,dm(\omega) +\varepsilon
\]
for every $\omega\in B^+$, contradicting the previous display.
Hence $\varrho(A^+)=0$.
Applying the same argument to
\[
A^-:=\Big\{(\omega,x)\in X : h_{\omega,f}(x)<\int_\Omega\varrho_\omega(f_\omega)\,dm(\omega)-\varepsilon\Big\}
\]
gives $\varrho(A^-)=0$.
Since $\varepsilon>0$ was arbitrary, we conclude that $h_f=\int_\Omega\varrho_\omega(f_\omega)\,dm(\omega)$ $\varrho$-a.e., that is,
\[
A_nf(\omega,x)\to \int_\Omega\varrho_\omega(f_\omega)\,dm(\omega)
\qquad \varrho\text{-a.e. }(\omega,x)
\]
for every bounded measurable random variable $f$.

In particular, if $A\in\mathscr F\otimes\mathscr B$ is $T$-invariant, then applying the previous conclusion to $f=1_A$ yields
\[
1_A(\omega,x)=\varrho(A)
\qquad \varrho\text{-a.e. }(\omega,x),
\]
so $\varrho(A)\in\{0,1\}$.
Thus $\varrho$ is ergodic.
The full statement \eqref{eq:20240614} for every $f\in L^1(\varrho)$ now follows from the usual Birkhoff ergodic theorem.

\subsection{Forward fibrewise mixing and pullback convergence}
Let $d\zeta_\omega=\psi_\omega\,d\varrho_\omega$, where
$\zeta_\omega$ is absolutely continuous with respect to $\varrho_\omega$.
For each $f\in C(\Xi)$,
\[
(T_*^n\zeta)_{\sigma^n\omega}(f)
=\varrho_\omega((f\circ T_\omega^n)\psi_\omega).
\]
Thus \eqref{eq:equivalence_mixing_2} implies forward convergence from each fixed initial fibre. In contrast, \eqref{eq:equivalence_mixing} is convergence at a fixed terminal fibre, with initial fibre $\sigma^{-n}\omega$.
When $\Omega$ is a singleton these coincide, and the converse follows by applying \eqref{eq:equivalence_mixing} to normalized nonnegative $L^1(\varrho)$ densities and taking linear combinations.

\section{Proof of Proposition \ref{cor:0201a}}
\label{app:survivor-invariance}

\subsection{Proof of Proposition \ref{cor:0201a}(1)}
Fix $\omega$ in a full-measure set on which the conclusions of
Theorem~\ref{main thm qstat and qerg measures} hold and every function
in $\mathcal B$ is
$\widehat\nu_{\omega,\infty}$-a.e.~continuous.

We first prove that
\begin{equation}\label{eq:nu-functional-measure-identification}
\nu_{\omega,\infty}(f)
=
\widehat\nu_{\omega,\infty}(f),
\qquad f\in\mathcal B.
\end{equation}
Let $f\in\mathcal B$.
By Theorem~\ref{main thm qstat and qerg measures}(4), $\nu_{\omega,n}$ converges weakly to $\widehat\nu_{\omega,\infty}$, whereas $f$ is $\nu_{\om,\infty}$-a.e.~continuous by the assumption.
Hence, Billingsley's mapping theorem gives weak convergence of the pushforward measures under $f$.
Their supports are contained in the closed disk of radius $C_0\|f\|$ by \eqref{B2};
hence integration of a bounded continuous function on $\mathbb C$ that agrees with the identity on that disk gives
$\nu_{\omega,n}(f)
\to
\widehat\nu_{\omega,\infty}(f)$.
On the other hand, Theorem~\ref{main thm qstat and qerg measures}(1),
applied with $\zeta=\nu_c$, gives
$\nu_{\omega,n}(f)\to\nu_{\omega,\infty}(f)$. 
The uniqueness of the limit proves
\eqref{eq:nu-functional-measure-identification}.

Now let
$\zeta\in\mathcal P_{\mathcal B,+}(\nu_c)$
have random density $\psi$.  By \eqref{B2}, $\psi_\omega$ is bounded,
and
$\nu_{\omega,\infty}(\psi_\omega)>0$,
so
$\frac{1}
     {\nu_{\omega,\infty}(\psi_\omega)}
\,\widehat\nu_{\omega,\infty}(\cdot \, \psi_\omega)$
is a finite Radon measure on $\Xi$.
Hence, Theorem~\ref{main thm qstat and qerg measures}(2), the definition of
$\zeta_{\omega,\infty}$, and
\eqref{eq:nu-functional-measure-identification} give
\begin{equation}\label{eq:representing-measure-densities}
\widehat\zeta_{\omega,\infty}(v)
=
\zeta_{\omega,\infty}(v)=
\frac{\nu_{\omega,\infty}(v\psi_\omega)}
     {\nu_{\omega,\infty}(\psi_\omega)}=
\frac{\widehat\nu_{\omega,\infty}(v\psi_\omega)}
     {\nu_{\omega,\infty}(\psi_\omega)},
     \qquad v\in\mathcal B\cap C(\Xi)
\end{equation}
By \eqref{B1}, this equality can be extended to that on $C(\Xi)$. 
Consequently, two finite Radon measures
$\widehat\zeta_{\omega,\infty}$ and
$\frac{1}
     {\nu_{\omega,\infty}(\psi_\omega)}
\,\widehat\nu_{\omega,\infty}(\cdot \, \psi_\omega)$
coincide. Hence, for every measurable function $f$,
\begin{equation*}
\widehat\zeta_{\omega,\infty}(f)
=
\frac{\widehat\nu_{\omega,\infty}(f\psi_\omega)}
     {\nu_{\omega,\infty}(\psi_\omega)}.
\end{equation*}
In particular, for every $f\in\mathcal B$, 
%since $f\psi_\omega\in\mathcal B$ by \eqref{B2},
\eqref{eq:nu-functional-measure-identification} implies
\begin{align*}
\widehat\zeta_{\omega,\infty}(f)
=
\frac{\widehat\nu_{\omega,\infty}(f\psi_\omega)}
     {\nu_{\omega,\infty}(\psi_\omega)}=
\frac{\nu_{\omega,\infty}(f\psi_\omega)}
     {\nu_{\omega,\infty}(\psi_\omega)}=
\zeta_{\omega,\infty}(f).
\end{align*}
This proves item~(1).

\subsection{\texorpdfstring{$\widetilde{\mathcal L}^n$}{L-tilde to the n} as a random closed transfer operator on \texorpdfstring{$X_\infty$}{X-infinity}}
We need the following duality relation to show items (2) and (3) of Proposition \ref{cor:0201a}.
\begin{lemma}\label{lem:borel-formulas-on-survivor}
Assume that \eqref{G} holds and that for $m$-a.e.~$\omega$, every function in $\mathcal B$ is $\widehat \nu_{\om,\infty}$-a.e.~continuous. 
Then, for $m$-a.e.~$\omega$, every $n\ge1$, and every $u\in\mathcal B$, it holds that
\begin{equation*}
(T_\omega^n)_*\bigl(u\,\widehat\nu_{\omega,\infty}\bigr)
=
\bigl(\widetilde{\mathcal L}_\omega^n u\bigr)
\widehat\nu_{\sigma^n\omega,\infty}.
\end{equation*}
Here $u\,\widehat\nu_{\omega,\infty}$ denotes the measure $A\mapsto\widehat\nu_{\omega,\infty}(u\hat A)$.
\end{lemma}

\begin{proof}
Fix $\omega$, $n\ge1$, and $u\in\mathcal B$.  
We first see that a function
\begin{equation}\label{eq:borel-set-function-on-source}
A\longmapsto
\widehat\nu_{\sigma^n\omega,\infty}
\left(
\widetilde{\mathcal L}_\omega^n(u\hat A)
\right)
\end{equation}
on the set of Borel sets on $\Xi$ is a finite complex Borel measure.
The expression in \eqref{eq:borel-set-function-on-source} is defined by
iterating the preimage sum in \eqref{eq: L_c def}, including the survivor
indicators in the definition of $\widetilde{\mathcal L}_\omega^n$.
Namely,
for every bounded Borel function $q$ and
$x\in \Xi$, 
\begin{equation}\label{eq:n-step-preimage-sum-borel}
\widetilde{\mathcal L}_\omega^n q(x)
:=
\frac{1}
{\lambda_{\omega}^n}
\sum_{
        y\in X_{\omega,n-1},\;
        T_\omega^n(y)=x}
\left(
 \prod_{j=0}^{n-1}
 g_{\sigma^j\omega}(T_\omega^j y)
\right)q(y),
\end{equation}
which is compatible to the definition of $\widetilde{\mathcal L}_\omega^n$ on $\cB$.

Let $A_1,A_2,\ldots$ be pairwise disjoint Borel subsets of $\Xi$.
The series in \eqref{eq:n-step-preimage-sum-borel} is absolutely
convergent when $q=u\hat A$.  Indeed,
\[
\begin{aligned}
\frac{1}
{\lambda_{\omega}^n}
\sum_{
        y\in X_{\omega,n-1},\;
        T_\omega^n(y)=x}
\left(
 \prod_{j=0}^{n-1}
 g_{\sigma^j\omega}(T_\omega^j y)
\right)
|u(y)|\hat A(y)
\le
\|u\|_\infty\,
\widetilde{\mathcal L}_\omega^n\hat X_\omega(x)
<\infty.
\end{aligned}
\]
Since the sets $A_r$ are pairwise disjoint,
$
\widehat{\bigcup_{r=1}^{\infty}A_r}(y)
=
\sum_{r=1}^{\infty}\hat A_r(y)
$
for every $y\in\Xi$.  Hence, Tonelli's theorem and the preceding 
estimate give
\[
\begin{aligned}
&\sum_{r=1}^{\infty}
\frac{1}
{\lambda_{\omega}^n}
\sum_{
        y\in X_{\omega,n-1},\;
        T_\omega^n(y)=x}
\left(
 \prod_{j=0}^{n-1}
 g_{\sigma^j\omega}(T_\omega^j y)
\right)
|u(y)|\hat A_r(y)\\
&=
\frac{1}
{\lambda_{\omega}^n}
\sum_{
        y\in X_{\omega,n-1},\;
        T_\omega^n(y)=x}
\left(
 \prod_{j=0}^{n-1}
 g_{\sigma^j\omega}(T_\omega^j y)
\right)
|u(y)|
\widehat{\bigcup_{r=1}^{\infty}A_r}(y)
\le
\|u\|_\infty\,
\widetilde{\mathcal L}_\omega^n\hat X_\omega(x).
\end{aligned}
\]
We may therefore interchange the two absolutely convergent sums in
\eqref{eq:n-step-preimage-sum-borel} with $q=u \widehat{\bigcup_{r=1}^{\infty}A_r}$.  It follows that
\[
\begin{aligned}
\widetilde{\mathcal L}_\omega^n
\Big(
 u\widehat{\bigcup_{r=1}^{\infty}A_r}
\Big)(x)
=
\sum_{r=1}^{\infty}
\widetilde{\mathcal L}_\omega^n(u\hat A_r)(x)
=
\lim_{N\to\infty}
\sum_{r=1}^{N}
\widetilde{\mathcal L}_\omega^n(u\hat A_r)(x).
\end{aligned}
\]
On the other hand, positivity of the preimage sum also gives, for every $N\ge 1$,
\begin{equation}\label{eq:borel-input-domination}
\Big|
\sum_{r=1}^{N}
\widetilde{\mathcal L}_\omega^n(u\hat A_r)(x)
\Big|
\le\|u\|_\infty\,\widetilde{\mathcal L}_\omega^n\Big(\hat X_\omega\widehat{\bigcup_{r=1}^{N}A_r}\Big)
\le
\|u\|_\infty\,
\widetilde{\mathcal L}_\omega^n\hat X_\omega(x).
\end{equation}
The last function is integrable with respect to
$\widehat\nu_{\sigma^n\omega,\infty}$.  Indeed,
\eqref{eq:zeta-functional-measure-identification}, applied with
$\zeta=\nu_c$, and \eqref{G1} imply
\[
\widehat\nu_{\sigma^n\omega,\infty}
\left(\widetilde{\mathcal L}_\omega^n\hat X_\omega\right)
=
\nu_{\sigma^n\omega,\infty}
\left(\widetilde{\mathcal L}_\omega^n\hat X_\omega\right)
=
\nu_{\omega,\infty}(\hat X_\omega)
=1.
\]
Dominated convergence therefore proves that
\eqref{eq:borel-set-function-on-source} is a finite complex Borel measure.

For every $v\in\mathcal B\cap C(\Xi)$, integration of $v$ with respect
to the measure in \eqref{eq:borel-set-function-on-source}, first for simple
functions and then by uniform approximation, equals
\[
\widehat\nu_{\sigma^n\omega,\infty}
\left(\widetilde{\mathcal L}_\omega^n(uv)\right).
\]
Since $uv\in\mathcal B$ by \eqref{B2},
\eqref{eq:zeta-functional-measure-identification}, applied with
$\zeta=\nu_c$, and \eqref{G1} give
\begin{align}\label{eq:202620720}
\widehat\nu_{\sigma^n\omega,\infty}
\left(\widetilde{\mathcal L}_\omega^n(uv)\right)
=
\nu_{\sigma^n\omega,\infty}
\left(\widetilde{\mathcal L}_\omega^n(uv)\right)
=
\nu_{\omega,\infty}(uv)
=
\widehat\nu_{\omega,\infty}(uv).
\end{align}
Thus the measure in \eqref{eq:borel-set-function-on-source} and
$u\,\widehat\nu_{\omega,\infty}$ have the same integrals on
$\mathcal B\cap C(\Xi)$, and hence on $C(\Xi)$ by \eqref{B1}.

Furthermore, the measure in \eqref{eq:borel-set-function-on-source} and $u\,\widehat\nu_{\omega,\infty}$ are Radon. 
Indeed, when $u=\hat X_\omega$, the measure in \eqref{eq:borel-set-function-on-source} is a positive probability measure.
Given $\varepsilon>0$, choose $q\in C(\Xi)$ with $0\le q\le1$ and $\widehat\nu_{\omega,\infty}(q)>1-\varepsilon$. 
Its integral with respect to the measure in \eqref{eq:borel-set-function-on-source} has the same value by \eqref{eq:202620720}, so that measure assigns mass greater than $1-\varepsilon$ to the compact set $\operatorname{supp}q$.  
It is therefore tight and hence Radon because $\Xi$ is a metric space.  
For general $u$, inequality \eqref{eq:borel-input-domination} implies that the total variation of the measure in \eqref{eq:borel-set-function-on-source} is bounded by $\|u\|_\infty$ times the Radon measure obtained when $u=\hat X_\omega$; it is therefore Radon as well.

In conclusion, by the Riesz--Markov representation theorem, it follows from \eqref{eq:202620720} that
\begin{equation*}
\widehat\nu_{\sigma^n\omega,\infty}
\left(\widetilde{\mathcal L}_\omega^n(u\hat A)\right)
=
\widehat\nu_{\omega,\infty}(u\hat A)
\end{equation*}
for every Borel set $A\subset\Xi$.
Therefore, for any Borel set $A\subset\Xi$, the pointwise equality
\begin{equation}\label{eq:borel-pullout-on-survivor}
\widetilde{\mathcal L}_\omega^n
\left(u\widehat{(T_\omega^n)^{-1}A}\right)
=
\hat A\,\widetilde{\mathcal L}_\omega^n u,
\end{equation}
which follows directly from the iterated preimage sum in
\eqref{eq: L_c def}, gives
\begin{align*}
(T_\omega^n)_*\bigl(u\,\widehat\nu_{\omega,\infty}\bigr)(A)
=
\widehat\nu_{\omega,\infty}
\left(u\widehat{(T_\omega^n)^{-1}A}\right)
=
\widehat\nu_{\sigma^n\omega,\infty}
\left(
\widetilde{\mathcal L}_\omega^n
\left(u\widehat{(T_\omega^n)^{-1}A}\right)
\right)
=
\widehat\nu_{\sigma^n\omega,\infty}
\left(\hat A\,\widetilde{\mathcal L}_\omega^n u\right).
\end{align*}
This proves the desired equality.
\end{proof}

\begin{remark}\label{rem:borel-extension-scope}
The proof of Lemma~\ref{lem:borel-formulas-on-survivor} evaluates
$\widetilde{\mathcal L}_\omega^n(u\hat A)$ and
$\widetilde{\mathcal L}_\omega^n
(u\widehat{(T_\omega^n)^{-1}A})$ even when the indicated bounded Borel
functions do not belong to $\mathcal B$.  These expressions are meaningful
in the present setting because \eqref{eq: L_c def} defines the operator by a
pointwise preimage sum.  The proof uses the two pointwise relations
\eqref{eq:borel-input-domination} and
\eqref{eq:borel-pullout-on-survivor}.

If $\widetilde{\mathcal L}_\omega$ is instead given only as an abstract
bounded operator $\mathcal B\to\mathcal B$, one should add either an extension
of $\widetilde{\mathcal L}_\omega^n$ to bounded Borel functions for which
\eqref{eq:borel-input-domination} and
\eqref{eq:borel-pullout-on-survivor} hold, or an algebra of Borel subsets
which generates the Borel $\sigma$-algebra and whose indicators and
required pullbacks admit representatives in $\mathcal B$, followed by a
monotone class argument.
\end{remark}

\subsection{Proof of Proposition \ref{cor:0201a}(2) and (3)}
By Lemma \ref{lem:borel-formulas-on-survivor} with $u=\phi_\omega$, followed by
\eqref{G2}, gives
\[
(T_\omega^n)_*\widehat\eta_{\omega,\infty}
=
\left(\widetilde{\mathcal L}_\omega^n\phi_\omega\right)
\widehat\nu_{\sigma^n\omega,\infty}
=
\phi_{\sigma^n\omega}\widehat\nu_{\sigma^n\omega,\infty}
=
\widehat\eta_{\sigma^n\omega,\infty},
\]
which is the $T$-invariance of $\widehat\eta_\infty$.

We next prove uniqueness.  Let $\rho$ be a $T$-invariant random probability
measure on $X_\infty$ such that
$\rho_\omega\ll\widehat\nu_{\omega,\infty}$.  Set
\[
   a(\omega):=
   \|\rho_\omega-\widehat\eta_{\omega,\infty}\|_{\mathrm{TV}}
   =
   \left\|
   \frac{d\rho_\omega}{d\widehat\nu_{\omega,\infty}}
   -\phi_\omega
   \right\|_{L^1(\widehat\nu_{\omega,\infty})}.
\]
Since pushforwards decrease total variation and both $\rho$ and
$\widehat\eta_\infty$ are $T$-invariant,
\[
   a(\sigma\omega)\le a(\omega)
   \qquad\text{for $m$-a.e.~$\omega$.}
\]
As $\sigma$ preserves $m$, this implies $a\circ\sigma=a$ $m$-a.e.; by
ergodicity of $(\sigma,m)$, $a$ is constant $m$-a.e.  Denote this constant
by $a_0$.

Fix $\omega$ in a forward-invariant full-measure set on which
\eqref{eq:stronggap}, Lemma \ref{lem:borel-formulas-on-survivor}, and
$a(\sigma^n\omega)=a_0$ hold for every $n\ge0$.  Let $\varepsilon>0$.
Since $\widehat\nu_{\omega,\infty}$ is a finite Radon measure,
$C(\Xi)$ is dense in $L^1(\widehat\nu_{\omega,\infty})$.
Condition~\eqref{B1} says that $\mathcal B\cap C(\Xi)$ is uniformly dense
in $C(\Xi)$; hence $\mathcal B$ is dense in
$L^1(\widehat\nu_{\omega,\infty})$.  Choose $u_\omega\in\mathcal B$ such
that
\[
   \left\|
   \frac{d\rho_\omega}{d\widehat\nu_{\omega,\infty}}
   -u_\omega
   \right\|_{L^1(\widehat\nu_{\omega,\infty})}
   <\varepsilon.
\]
For every $n\ge1$,
\begin{align*}
&a_0=a(\sigma^n\omega)\\
&\le
\left\|(T_\omega^n)_*
\left(\rho_\omega-u_\omega\widehat\nu_{\omega,\infty}\right)
\right\|_{\mathrm{TV}}+
\left\|(T_\omega^n)_*
\left(u_\omega\widehat\nu_{\omega,\infty}\right)
      -\nu_{\omega,\infty}(u_\omega)
       \widehat\eta_{\sigma^n\omega,\infty}\right\|_{\mathrm{TV}}+
|\nu_{\omega,\infty}(u_\omega)-1|.
\end{align*}
The first term is at most $\varepsilon$, because pushforward does not
increase total variation.  The third term is also at most $\varepsilon$.
Indeed, \eqref{eq:zeta-functional-measure-identification}, applied with
$\zeta=\nu_c$, identifies $\nu_{\omega,\infty}(u_\omega)$ with
$\widehat\nu_{\omega,\infty}(u_\omega)$, while
$\widehat\nu_{\omega,\infty}
\left(\frac{d\rho_\omega}{d\widehat\nu_{\omega,\infty}}\right)=1$.
Consequently,
\[
\begin{aligned}
|\nu_{\omega,\infty}(u_\omega)-1|
=
\left|
\widehat\nu_{\omega,\infty}
\left(
 u_\omega-
 \frac{d\rho_\omega}{d\widehat\nu_{\omega,\infty}}
\right)
\right|
\le
\left\|
 u_\omega-
 \frac{d\rho_\omega}{d\widehat\nu_{\omega,\infty}}
\right\|_{L^1(\widehat\nu_{\omega,\infty})}
<\varepsilon.
\end{aligned}
\]
For the second term, Lemma \ref{lem:borel-formulas-on-survivor} and
\eqref{eq:stronggap} give
\begin{align*}
\left\|
   (T_\omega^n)_*
   \left(u_\omega\widehat\nu_{\omega,\infty}\right)
   -
   \nu_{\omega,\infty}(u_\omega)
   \widehat\eta_{\sigma^n\omega,\infty}
\right\|_{\mathrm{TV}}
&=
\left\|
   \widetilde{\mathcal L}_\omega^n u_\omega
   -
   \nu_{\omega,\infty}(u_\omega)\phi_{\sigma^n\omega}
\right\|_{L^1(\widehat\nu_{\sigma^n\omega,\infty})}\\
&\le
C_0D_{\sigma^n\omega}\|u_\omega\|\kappa^n.
\end{align*}
Letting $n\to\infty$ gives $a_0\le2\varepsilon$ since $D$ is tempered,
and then letting $\varepsilon\to0$ gives $a_0=0$.  Hence
$\rho=\widehat\eta_\infty$.

It remains to prove ergodicity.  Let $A\subset X_\infty$ be a
$T$-invariant measurable set, meaning that
\[
   \hat A_{\sigma\omega}\circ T_\omega
   =
   \hat A_\omega
   \qquad
   \widehat\eta_{\omega,\infty}\text{-a.e.}
\]
Then the invariance of $\widehat\eta_\infty$ gives
$\widehat\eta_{\sg\omega,\infty}(A_{\sg\omega})=\widehat\eta_{\omega,\infty}(A_\omega)$ for $m$-a.e.~$\oio$; hence $\widehat\eta_{\omega,\infty}(A_\omega)$ is $m$-a.e.~constant, which we denote by $c_0$.
If $0<c_0<1$, then
\[
   \rho_\omega(B)
   :=
   \frac{\widehat\eta_{\omega,\infty}(B\cap A_\omega)}{c_0}
\]
defines a $T$-invariant random probability measure on $X_\infty$ which is
absolutely continuous with respect to $\widehat\nu_\infty$.  By uniqueness,
$\rho=\widehat\eta_\infty$, and therefore
\[
   1=\rho_\omega(A_\omega)
   =\widehat\eta_{\omega,\infty}(A_\omega)=c_0,
\]
a contradiction.  Thus $c_0\in\{0,1\}$, and
$\widehat\eta_\infty$ is ergodic.
This proves item~(2).

Applying Lemma \ref{lem:borel-formulas-on-survivor} with $u=v\psi_\omega$ to the bounded Borel function $f$ gives
\begin{align*}
&\widehat\zeta_{\omega,\infty}
\bigl((f\circ T_\omega^n)v\bigr)
-
\widehat\eta_{\sigma^n\omega,\infty}(f)\,
\widehat\zeta_{\omega,\infty}(v)=
\frac{1}{\nu_{\omega,\infty}(\psi_\omega)}
\widehat\nu_{\sigma^n\omega,\infty}
\left(
 f\left[
 \widetilde{\mathcal L}_\omega^n(v\psi_\omega)
 -\nu_{\omega,\infty}(v\psi_\omega)
  \phi_{\sigma^n\omega}
 \right]
\right).
\end{align*}
Equation~\eqref{eq:stronggap}, applied to $v\psi_\omega$, gives
\[
\left\|
\widetilde{\mathcal L}_\omega^n(v\psi_\omega)
-\nu_{\omega,\infty}(v\psi_\omega)\phi_{\sigma^n\omega}
\right\|
\le D_{\sigma^n\omega}\|v\psi_\omega\|\kappa^n.
\]
Using this inequality and \eqref{B2} gives
\begin{align*}
&\left|
\widehat\zeta_{\omega,\infty}
\bigl((f\circ T_\omega^n)v\bigr)
-
\widehat\eta_{\sigma^n\omega,\infty}(f)\,
\widehat\zeta_{\omega,\infty}(v)
\right|\le
\frac{C_0\|f\|_\infty\|v\psi_\omega\|}
     {\nu_{\omega,\infty}(\psi_\omega)}
D_{\sigma^n\omega}\kappa^n\le
C_{\omega,\zeta,f,v}
\bigl(D_{\sigma^n\omega}C_{\sigma^n\omega,\phi}\bigr)\kappa^n,
\end{align*}
where the last inequality follows from
$C_{\sigma^n\omega,\phi}\ge1$ and the definition of
$C_{\omega,\zeta,f,v}$.  This proves item~(3).

\section{\texorpdfstring{$\widehat\nu_\infty$}{nu-hat-infinity}-a.e.~continuity for multidimensional expanding maps}\label{app:multidimensional-ae-continuity}

This appendix gives the proof of $\widehat\nu_\infty$-a.e.~continuity of the representative $f^\#$ in Remark~\ref{rem:eqclass} for every equivalence class $[f]$ of $V_\alpha$ in the multidimensional expanding setting.
By this fact and Proposition \ref{cor:0201a}(1), we can identify $\nu_{\omega,\infty}$ with integration against $\widehat\nu_{\omega,\infty}$ on $V_\alpha(X)$ in the sense of Remark~\ref{rem:eqclass}.
%with the representative $f^\#$.
%fixed in Remark~\ref{rem:eqclass}.
The same argument gives the almost everywhere continuity required in Corollary~\ref{cor:CLTLDPXinfty} and Proposition~\ref{prop:0625a}(3),(4).
For $E\subset X$ and $r>0$, write
\[
   B_r(E):=\{x\in X:d(x,E)<r\}.
\]

\begin{lemma}
\label{lem:neighborhoodNull}
Assume~\eqref{G}, and fix $\omega$ in a set of full measure on which the
conclusions of Theorem~\ref{main thm qstat and qerg measures} hold.  Let
$E\subset X$ be closed, and suppose that, for some $C<\infty$ and $r_0>0$,
\[
   \Leb(B_r(E))\le Cr^\alpha,
   \qquad 0<r<r_0.
\]
Then $\widehat\nu_{\omega,\infty}(E)=0$.
\end{lemma}

\begin{proof}
Choose a Lipschitz function $\vartheta:[0,\infty)\to[0,1]$ such that $\vartheta=1$ on $[0,1]$ and $\vartheta=0$ on $[2,\infty)$, and, for $0<r<r_0/3$, put
\[
   h_r(x):=\vartheta\!\left(\frac{d(x,E)}r\right).
\]
Then $h_r=1$ on $E$ and $h_r=0$ on $X\setminus B_{2r}(E)$.
If $0<\varepsilon\le r$, the $\varepsilon$-oscillation function $\operatorname{osc}(h_r,B_\varepsilon(x)\cap X)$ vanishes outside $B_{3r}(E)$ and is at most $2\mathrm{Lip}(\vartheta)\varepsilon/r$ on this set.
Hence, it follows from the assumption that
\[
   \varepsilon^{-\alpha}
   \int_X\operatorname{osc}(h_r,B_\varepsilon(x)\cap X)\,d\Leb(x)
   \le 2\mathrm{Lip}(\vartheta)   \varepsilon^{-\alpha}\frac\varepsilon r\Leb(B_{3r}(E))
   \le 2C3^\alpha\mathrm{Lip}(\vartheta)\left(\frac\varepsilon r\right)^{1-\alpha}.
\]
If $r<\varepsilon\le\min\{\varepsilon_0,r_0/3\}$, the $\varepsilon$-oscillation is at most one and vanishes unless $x\in B_{3\varepsilon}(E)$ because $d(y,E)\ge 3\varepsilon -\varepsilon>2r$ for any $x\not\in B_{3\varepsilon}(E)$ and $y\in B_\varepsilon(x)$, so
\[
   \varepsilon^{-\alpha}
   \int_X\operatorname{osc}(h_r,B_\varepsilon(x)\cap X)\,d\Leb(x)
   \le \varepsilon^{-\alpha} \Leb(B_{3\varepsilon}(E))
   \le C3^\alpha.
\]
If $\min\{\varepsilon_0,r_0/3\}<\varepsilon\le\varepsilon_0$,
the same quantity is at most $\Leb(X)\min\{\varepsilon_0,r_0/3\}^{-\alpha}$ due to the trivial bound of the $\varepsilon$-oscillation  by $1$.
For the zero extension of $h_r$,
\[
\begin{aligned}
\int_{\mathbb R^d}\operatorname{osc}(h_r,B_\varepsilon(x))
\,d\Leb(x)
&\le
\int_X\operatorname{osc}(h_r,B_\varepsilon(x)\cap X)
\,d\Leb(x)\\
&\quad+
\Leb(\{x\in\mathbb R^d:
\operatorname{dist}(x,\partial X)<\varepsilon\}).
\end{aligned}
\]
By \eqref{P3}, the last term is at most $C_X\varepsilon$.
Since $\alpha\le1$, its contribution to $|h_r|_\alpha$
is bounded uniformly in $r$.
Since $\|h_r\|_1\le\Leb(X)$, it follows that
\begin{equation}\label{eq:20260802}
   \sup_{0<r<r_0/3}\|h_r\|_\alpha<\infty.
\end{equation}

Theorem~\ref{main thm qstat and qerg measures}(2) gives
\[
   \nu_{\omega,\infty}(h_r)
   =\widehat\nu_{\omega,\infty}(h_r),
\]
while Theorem~\ref{main thm qstat and qerg measures}(1) and the uniform bound \eqref{eq:20260802} give
\[
   \sup_{0<r<r_0/3}
   |\nu_{\omega,N}(h_r)-\nu_{\omega,\infty}(h_r)|\longrightarrow0.
\]
Suppose that $\widehat\nu_{\omega,\infty}(E)>0$.  Choose $N$ so large
that the last supremum is less than
$\widehat\nu_{\omega,\infty}(E)/4$. 
The assumed estimate implies $\Leb(E)=0$, and hence $\nu_{\omega,N}(E)=\Leb(E\cap X_{\omega,N})/\Leb(X_{\omega,N})=0$.
Since $E$ is closed,
$B_{2r}(E)\downarrow E$ as $r\downarrow0$, so $r$ may be chosen such that
\[
   \nu_{\omega,N}(B_{2r}(E))
   <\frac14\widehat\nu_{\omega,\infty}(E).
\]
Hence, since $\widehat E(x)\le h_r(x)$ for every $x$,
%For this $r$,
\[
   \widehat\nu_{\omega,\infty}(E)
   \le \widehat\nu_{\omega,\infty}(h_r)
    =\nu_{\omega,\infty}(h_r)
   \le
   |\nu_{\omega,\infty}(h_r)-\nu_{\omega,N}(h_r)|
   +\nu_{\omega,N}(B_{2r}(E))
   <\frac12\widehat\nu_{\omega,\infty}(E),
\]
which is a contradiction.
\end{proof}

\begin{proposition}
\label{prop:ValphaIdentification}
Assume~\eqref{G}.  In the setting of
Section~\ref{s:multidimensional-expanding}, for every
$f\in L^\infty(\Omega,V_\alpha(X))$ and $m$-a.e.~$\omega$,
\[
   \widehat\nu_{\omega,\infty}
   \bigl(\operatorname{Disc}(f_\omega^\#)\bigr)=0.
\]
If the boundary of every element of $\mathcal U$ is contained in a finite
union of $C^1$ hypersurfaces, then, for every $n\ge0$,
\[
   \widehat\nu_{\omega,\infty}
   \bigl(\operatorname{Disc}
      (f_{\sigma^n\omega}^\#\circ T^n)\bigr)=0.
\]
Here $\operatorname{Disc}(g)$ denotes the discontinuity set of $g$.
\end{proposition}

\begin{proof}
Fix $f\in L^\infty(\Omega,V_\alpha(X))$ and $\omega\in\Omega$.
For $k,\ell\ge1$, put
\[
   E_{\om,k,\ell}:=
   \left\{
      x\in X:
      \operatorname{osc}(f_\om,B_s(x))\ge k^{-1}
      \text{ for every }0<s<\ell^{-1}
   \right\}.
\]
These sets are closed.
Indeed, if $x_j\to x$, $x_j\in E_{\om,k,\ell}$, and
$0<s<\ell^{-1}$, then, for all sufficiently large $j$,
$B_{s/2}(x_j)\subset B_s(x)$, and hence
$\operatorname{osc}(f,B_s(x))\ge k^{-1}$.
The defining property of $f^\#_\om$ in Remark~\ref{rem:eqclass} gives
\[
   \operatorname{Disc}(f^\#_\om)
   \subset
   \bigcup_{k,\ell\ge1}E_{\om,k,\ell}.
\]
If $0<r<\min\{\ell^{-1},\varepsilon_0/2\}$ and
$y\in B_r(E_{\om,k,\ell})$, then
$\operatorname{osc}(f_\om,B_{2r}(y))\ge k^{-1}$.
Therefore
\begin{equation}\label{eq:20260802b}
   \Leb(B_r(E_{\om,k,\ell}))
   \le
   2^\alpha k\|f_\om\|_\alpha r^\alpha.
\end{equation}
Lemma~\ref{lem:neighborhoodNull}, applied to each
$E_{\om,k,\ell}$, yields
\[
   \widehat\nu_{\omega,\infty}
   \bigl(\operatorname{Disc}(f^\#_\om)\bigr)=0.
\]

It remains to consider the compositions.
The assertion for $n=0$ has already been proved.
Let $n\ge1$, and let $\mathcal U^{(n)}$ denote the partition into
cylinders of length $n$.
The discontinuity set of $T^n$ is contained in a finite union of
$C^1$ hypersurfaces.
Each such hypersurface has an $r$ neighborhood of Lebesgue measure
$O(r)$, and hence $O(r^\alpha)$.
Moreover,
\begin{equation}\label{eq:20260802c}
\begin{split}
   \operatorname{Disc}(f^\#_{\sg^n\om}\circ T^n)
   &\subset
   \operatorname{Disc}(T^n)
   \cup
   (T^n)^{-1}
   \bigl(\operatorname{Disc}(f^\#_{\sg^n\om})\bigr)
   \\
   &\subset
   \operatorname{Disc}(T^n)
   \cup
   \bigcup_{V\in\mathcal U^{(n)}}
   \bigcup_{k,\ell\ge1}
   \left[
      \overline V
      \cap
      (\widetilde T_V^n)^{-1}
      (E_{\sg^n\om,k,\ell})
   \right],
\end{split}
\end{equation}
where $\widetilde T_V^n$ is a $C^1$ diffeomorphic extension of
$T^n|_V$ to a neighborhood of $\overline V$.
After restricting each extension to a relatively compact neighborhood
of $\overline V$, the finiteness of $\mathcal U^{(n)}$ gives a constant
$C_n\ge1$ such that
\begin{equation}\label{eq:20260802d}
   d\bigl(\widetilde T_V^n(x),\widetilde T_V^n(z)\bigr)
   \le C_n d(x,z)
\qquad
\text{and}
\qquad
   \left|
      \det D_y(\widetilde T_V^n)^{-1}
   \right|
   \le C_n
\end{equation}
on the corresponding neighborhoods, for every
$V\in\mathcal U^{(n)}$.

Fix $V\in\mathcal U^{(n)}$ and $k,\ell\ge1$.
For all sufficiently small $r>0$, the $r$-neighborhood of $
   \overline V
   \cap
   (\widetilde T_V^n)^{-1}
   (E_{\sg^n\om,k,\ell})
$ is contained in
the neighborhood on which \eqref{eq:20260802d} holds.
If
$
   x\in
   B_r\!(
      \overline V
      \cap
      (\widetilde T_V^n)^{-1}
      (E_{\sg^n\om,k,\ell})
   ),
$
choose
$
   z\in
   \overline V
   \cap
   (\widetilde T_V^n)^{-1}
   (E_{\sg^n\om,k,\ell})
$
such that $d(x,z)<r$.
Since
$\widetilde T_V^n(z)\in E_{\sg^n\om,k,\ell}$, we have
\[
   d\bigl(
      \widetilde T_V^n(x),
      E_{\sg^n\om,k,\ell}
   \bigr)
   \le
   d\bigl(
      \widetilde T_V^n(x),
      \widetilde T_V^n(z)
   \bigr)
   \le C_n r.
\]
Consequently,
\[
   \widetilde T_V^n
   \left(
      B_r\!\left(
         \overline V
         \cap
         (\widetilde T_V^n)^{-1}
         (E_{\sg^n\om,k,\ell})
      \right)
   \right)
   \subset
   B_{C_n r}(E_{\sg^n\om,k,\ell}).
\]
The change of variables formula therefore gives
\begin{align*}
   \Leb\!\left(
      B_r\!\left(
         \overline V
         \cap
         (\widetilde T_V^n)^{-1}
         (E_{\sg^n\om,k,\ell})
      \right)
   \right)
   &=
   \int_{
      \widetilde T_V^n
      \left(
         B_r\left(
            \overline V
            \cap
            (\widetilde T_V^n)^{-1}
            (E_{\sg^n\om,k,\ell})
         \right)
      \right)
   }
   \left|
      \det D_y(\widetilde T_V^n)^{-1}
   \right|
   \,d\Leb(y)
   \\
   & \le
   C_n
   \Leb\bigl(B_{C_n r}(E_{\sg^n\om,k,\ell})\bigr).
\end{align*}
If $r$ is small enough that
$
   C_n r
   <
   \min\{\ell^{-1},\varepsilon_0/2\},
$
the inequality \eqref{eq:20260802b}, applied with $\sg^n\om$ in place of $\om$,
gives
\[
\Leb\!\left(
      B_r\!\left(
         \overline V
         \cap
         (\widetilde T_V^n)^{-1}
         (E_{\sg^n\om,k,\ell})
      \right)
   \right)
    \le
   2^\alpha k
   \|f_{\sg^n\om}\|_\alpha
   C_n^{1+\alpha}r^\alpha.
\]

Lemma~\ref{lem:neighborhoodNull} applies to the finitely many
$C^1$ hypersurfaces containing $\operatorname{Disc}(T^n)$ and to each
of the closed sets
$
   \overline V
   \cap
   (\widetilde T_V^n)^{-1}
   (E_{\sg^n\om,k,\ell}).
$
Since $\mathcal U^{(n)}$ is finite and $k,\ell$ range over a countable
set, the inclusion \eqref{eq:20260802c} yields
$
   \widehat\nu_{\omega,\infty}
   \bigl(
      \operatorname{Disc}
      (f^\#_{\sg^n\om}\circ T^n)
   \bigr)
   =0.
$
\end{proof}

\section*{AI statement}
The core of this paper, including its main ideas, theoretical framework, and much of its mathematical development, was developed by the authors in 2023–2025 without AI assistance. The authors subsequently used ChatGPT for literature review, drafting and refining parts of technical proofs, and checking mathematical arguments. All AI-assisted content has been critically reviewed and revised by the authors, who take full responsibility for the results and their proofs.

\section*{Acknowledgements}
JA was supported by an Australian Research Council Discovery Project
(DP210100357).
GF was supported by an Australian Research Council Laureate Fellowship
(FL230100088) and in part by an Einstein Visiting Fellowship
(Einstein Foundation Berlin).
CGT acknowledges support from the Australian Research Council.
YN was supported by JSPS KAKENHI Grant Number~JP23K03188 and
JST PRESTO Grant Number~JPMJPR25K8.
SV was supported by the project Dynamics and Information Research
Institute under the agreement between UniCredit Bank and the
Scuola Normale Superiore di Pisa, and by the Laboratoire International
Associ\'e (LIA) LYSM of the French CNRS and INdAM (Italy).
SV was also supported by the French CNRS and MEAE through the
MATHAmSud project TOMCAT (22-Math-10, No.~49958WH) and by the CNRS International Research Laboratory IRL 2032 FAMSI, France Australia, at the Australian National University in Canberra,
where this work was completed.

The authors thank MATRIX, the University of Pisa, Waseda University,
the Centre de Physique Th\'eorique in Marseille, Hokkaido University,
the University of New South Wales, the University of Queensland,
and the Isaac Newton Institute for Mathematical Sciences
for their support and hospitality.
Parts of this research were carried out at these institutions.

\bibliographystyle{abbrv}
\bibliography{RQE}
\end{document}